\documentclass[reqno,11pt]{amsart}
\usepackage{amsfonts,amsmath,amssymb,fullpage,array,rotate}
\theoremstyle{plain} \newtheorem{Thm}{Theorem}[section]
\theoremstyle{plain} \newtheorem{Cor}[Thm]{Corollary}
\theoremstyle{plain} \newtheorem{Prop}[Thm]{Proposition}
\theoremstyle{plain} \newtheorem{Lemma}[Thm]{Lemma}
\theoremstyle{plain} 
\theoremstyle{definition} \newtheorem{Def}[Thm]{Definition}
\theoremstyle{remark} \newtheorem{Rem}[Thm]{Remark} 
\theoremstyle{definition} \newtheorem{Ex}[Thm]{Example} 

\def\pf{\noindent {\it Proof.}\hskip 8pt}
\def\qed{\hfill{\setlength{\fboxrule}{0.8pt}\setlength{\fboxsep}{1mm}
\fbox{\null}} \vskip 10pt}
\def\nqed{\hfill{\setlength{\fboxrule}{0.8pt}\setlength{\fboxsep}{1mm}
\fbox{\null}}}

\newcommand{\clearemptydoublepage}
             {\newpage{\pagestyle{empty}\cleardoublepage}}

\newcommand{\thmlist}%
          {\renewcommand{\theenumi}{\alph{enumi}}
                  \renewcommand{\labelenumi}{{\rm (\theenumi)}}
                            }
\newcommand{\factlist}%
          {\renewcommand{\theenumi}{\arabic{enumi}}
                  \renewcommand{\labelenumi}{{\rm (\theenumi)}}
                            }
\newcommand{\factlisti}%
          {
                  
                            }
\newcommand{\condlisti}%
          {
                  \renewcommand{\labelenumi}{{\rm(\theenumii)}}
                            }
\newcommand{\condlist}%
          {\renewcommand{\theenumi}{\roman{enumi}}
                  \renewcommand{\labelenumi}{{\rm (\theenumi)}}
                            }

\def\frak{\mathfrak}
\renewcommand{\Re}{\mathop{\rm{Re}}}
\renewcommand{\Im}{\mathop{\rm{Im}}}

\newcommand{\conv}{\mathop{\rm{conv}}}
\newcommand{\id}{\mathop{\rm{id}}}

\newcommand{\supp}{\mathop{\rm{supp}}}

\newcommand{\Span}{\mathop{\rm{span}}}

\newcommand{\Pav}{\mathop{\rm{P^{av}}}}
\newcommand{\PW}{\mathop{\rm{PW}}}
\newcommand{\abs}[1]{\left|\/#1\/\right|}

\newcommand{\inner}[2]{\langle#1,#2\rangle}

\newcommand{\C}{\ensuremath{\mathbb C}}

\newcommand{\D}{\ensuremath{\mathbb D}}

\newcommand{\R}{\ensuremath{\mathbb R}}
\newcommand{\Z}{\ensuremath{\mathbb Z}}
\newcommand{\N}{\ensuremath{\mathbb N}}

\renewcommand{\l}{\lambda}
\renewcommand{\a}{\alpha}
\renewcommand{\b}{\beta}
\newcommand{\e}{\varepsilon}

\newcommand{\la}{\l_\a}

\newcommand\SO{\mathop{\rm{SO}}}

\newcommand{\frakacs}{\frak a_{\scriptscriptstyle{\C}}^*}

\newcommand{\liecomplex}[1]{{\frak #1}_{\scriptscriptstyle{\C}}}
\newcommand{\complex}[1]{{#1}_{\scriptscriptstyle{\C}}}

\newcommand{\cmax}{C_{\rm max}}

\newcommand{\Hreg}{\complex{A}^{\rm reg}}

\newcommand{\rootstheta}{\langle\Theta\rangle}

\newcommand{\polya}{{\rm S}(\liecomplex{a})}
\newcommand{\polyAc}{\C[\complex{A}]}
\newcommand{\ma}{m_\a}
\newcommand{\mduea}{m_{2\a}}

\newcommand{\pedtheta}[1]{{#1}_{\scriptscriptstyle{\Theta}}}
\newcommand{\aptheta}[1]{{#1}^{\scriptscriptstyle{\Theta}}}

\newcommand{\Pavtheta}{\pedtheta{\rm P}^{\rm av}}
\newcommand{\pedempty}[1]{{#1}_{\scriptscriptstyle{\emptyset}}}
\newcommand{\ped}[2]{{#2}_{\scriptscriptstyle{#1}}}

\newcommand{\pedR}[1]{{#1}_{\scriptscriptstyle{R}}}

\newcommand{\pedPi}[1]{{#1}_{\scriptscriptstyle{\Pi}}}

\newcommand{\cthetampl}{\pedtheta{c}^+(m;\l)}
\newcommand{\cthetamml}{\pedtheta{c}^-(m;\l)}

\newcommand{\cCa}{C_c^\infty(C)}
\newcommand{\cCA}{C_c^\infty(C)}

\newcommand{\pwthetamC}{\pedtheta{\rm PW}(m;C)}

\newcommand{\cctheta}{C_c(\pedtheta{A})^{\pedtheta{W}}}
\newcommand{\ccitheta}{C_c^\infty(\pedtheta{A})^{\pedtheta{W}}}

\newcommand{\smC}{\scriptscriptstyle{\mathbb{C}}}

\begin{document}
\makeatletter
\title{A Paley-Wiener theorem for the $\Theta$-spherical transform: \\ the even 
multiplicity case}
\author{Gestur \'{O}lafsson}
\address{Department of Mathematics, Louisiana State University, Baton Rouge, 
LA 70803, U.S.A.}
\email{olafsson@math.lsu.edu}
\thanks{The first author was supported by NSF grant DMS-0070607 and DMS-0139783}
\thanks{Both authors were partially supported by DFG-Schwerpunkt ``Global 
methods in complex geometry'' of the Ruhr-Universit\"at Bochum and by the Lorentz Center
at the Rijksuniversiteit Leiden.}
\author{Angela Pasquale}
\address{Institut f\"ur Mathematik, TU-Clausthal,
38678 Clausthal-Zellerfeld, Germany.}
\email{mapa@math.tu-clausthal.de}
\date{}

\begin{abstract}
The $\Theta$-spherical functions generalize the
spherical functions on Riemannian symmetric
spaces and the spherical functions on non-compactly
causal symmetric spaces. In this article we consider
the case of even multiplicity functions. We construct
a differential shift operator $D_m$ with smooth
coefficients which generates the $\Theta$-spherical
functions from finite sums of exponential functions.
We then use this fact to prove a Paley-Wiener theorem
for the $\Theta$-spherical transfrom.     
\end{abstract}

\maketitle
\makeatother

\section{Introduction}

Let $\frak a$ be an $l$-dimensional Euclidean real vector space with complexification
$\liecomplex{a}$, and let $\frak a^*$ and $\frakacs$ respectively denote the real and 
complex dual vector spaces of $\frak a$.  
For a compact subset  $E$ of $\frak a$ let $\conv(E)$ denote
its closed convex hull 
(i.e. the intersection of all closed half-spaces in $\frak a$ containing $E$). 
The support function of $E$ is the function $q_E:\frak a^* \to \R$ defined by
\begin{equation}
  \label{eq:support}
  q_E(\l):=\sup_{H \in E} \l(H)=\sup_{H \in \conv(E)} \l(H).
\end{equation}
Let $C$ be a compact convex subset of $\frak a$, and let $\cCa$ denote the
space of smooth functions on $\frak a$ with support contained in $C$.
The Paley-Wiener space ${\rm PW}(C)$ consists of the entire functions $g:\frakacs\to \C$ 
which are of exponential type $C$ and rapidly decreasing, i.e., for every $N \in \N$  
there is a constant $C_N\geq 0$ such that
$$
\abs{g(\l)} \leq C_N (1+\abs{\l})^{-N} e^{q_C(\Re\l)}
$$
for all $\l \in \frakacs$. 
The Euclidean Fourier transform of a  sufficiently regular
function $f: \frak a \to \C$ is the function $\mathcal Ff:\frakacs \to \C$ defined by 
\begin{equation*}
 \mathcal Ff(\l):=\int_{\frak a} f(H) e^{\l(H)} \; dH, 
\end{equation*}
where $dH$ denotes the Lebesgue measure on $\frak a$. 
A classical theorem due to Paley and Wiener 
characterizes ${\rm PW}(C)$ as 
the image of $\cCa$ under the Euclidean Fourier transform 
(see e.g. \cite{HoerLPDO}, Theorem 7.3.1, or \cite{JL}, Theorem 8.3 and Proposition
8.6).

\begin{Thm}[Paley-Wiener] \label{thm:classicalPW}
Let $C$ be a compact convex subset of $\frak a$. Then the Euclidean Fourier transform
maps  $\cCa$ bijectively onto ${\rm PW}(C)$.
Moreover, if $C$ is stable under the action of a finite group $W$
of linear automorphisms of $\frak a$, then the Fourier transform
maps the subspace ${\cCa}^W$ of $W$-invariant elements in $\cCa$
onto the subspace ${\rm PW}(C)^W$ of $W$-invariant elements in ${\rm PW}(C)$.
\end{Thm}

Suppose that $A$ is a connected simply connected  abelian Lie  group with Lie algebra $\frak a$.
Then $\exp:\frak a\to A$ is a diffeomorphism. Denote by $\log$  the inverse of $\exp$.
The Euclidean Fourier transform of a sufficiently regular
functions $f: A \to \C$  is the function  $\mathcal F_Af: \frakacs \to \C$ defined by
\begin{equation*}
 \mathcal F_Af(\l):=\int_{A} f(a) e^{\l(\log a)} \; da,
\end{equation*}
where the Haar measure $da$ on $A$ is the pullback under the exponential map of the
Haar measure $dH$ on $\frak a$.
Let $W$ be a finite group acting on $\frak a$ by linear automorphism.
Then we define an action by $W$ on $A$ by $w (\exp H)=\exp w(H)$.
Denote by $\cCA$  the space of smooth functions on $A$ with
support in the compact set $\exp C$, and let $\cCA^W$ be the subspace of $W$-invariant
elements. Composition with $\exp$ and Theorem \ref{thm:classicalPW} prove that
the Euclidean Fourier transform $\mathcal F_A$ is a bijection of
$\cCA$ onto ${\rm PW}(C)$ and restricts to a bijection of $\cCa^W$ onto ${\rm PW}(C)^W$.

In this paper we generalize Theorem \ref{thm:classicalPW}  to the case of the
$\Theta$-spherical transform of \cite{PHab} and \cite{P2} corresponding to a triple
$(\frak a, \Sigma, m)$. Here $\frak a$ is a Euclidean space, $\Sigma$ an irreducible 
reduced root system in the dual space $\frak a^*$ of $\frak a$, and $m$ 
a positive even multiplicity function on $\Sigma$ 
(see Section \ref{section:prelim} for the precise definitions). 
Furthermore, $\Theta$ denotes a subset in a given fundamental system 
$\Pi$ of a fixed set $\Sigma^+$ of positive simple roots in $\Sigma$. 

Let $\pedtheta{W}$ denote the parabolic subgroup of $W$ generated by the 
reflections $r_\a$ with $\a \in \Theta$.
With each choice of $\Theta$ and $\l \in \frakacs$ is associated a function 
$\pedtheta{\varphi}(m;\l,a)$, which is defined for $a=\exp H$ with $H$ 
in a certain $\pedtheta{W}$-invariant open cone $\pedtheta{\frak a}$ in $\frak a$. 
We set $\pedtheta{A}:=\exp(\pedtheta{\frak a})$.
The function $\pedtheta{\varphi}(m;\l)$
is called the {\em $\Theta$-spherical function of spectral parameter $\l$}. 
It is a common eigenfunction
of the hypergeometric system of differential operators of spectral
parameter $\l$ constructed by Heckman and Opdam. It is real analytic and
$\pedtheta{W}$-invariant. As a function of $\l \in \frakacs$,  
$\pedtheta{\varphi}(m;\l,a)$ is meromorphic with simple poles located 
along a specific
finite family of affine complex hyperplanes. 
Indeed, there is a polynomial function 
$\pedtheta{e}^-(m;\l)$, which is a finite product of affine functions of $\l$, and 
there is a tubular neighborhood 
$\pedtheta{U}$ of $\pedtheta{A}$ in the complexification $\complex{A}$ of 
$A$ so that
$\pedtheta{e}^-(m;\l) \pedtheta{\varphi}(m;\l,a)$ extends as a holomorphic function 
of $(\l,a) \in \frakacs \times \pedtheta{A}$.
We refer the reader to Theorem  \ref{thm:sphentire} below
for detailed information on the regularity properties of the $\Theta$-spherical functions.

The $\Theta$-spherical transform of a $\pedtheta{W}$-invariant function 
$f$ on  $\pedtheta{A}$ is the $\pedtheta{W}$-invariant 
function $\pedtheta{\mathcal F}f(m)$ on $\frakacs$ defined for 
$\l \in \frakacs$ by
\begin{equation} \label{eq:thetasphtr}
\pedtheta{\mathcal F}f(m;\l):=
\frac{1}{\abs{\pedtheta{W}}} \, \int_{\pedtheta{A}} 
f(a) \, \pedtheta{\varphi}(m;\l,a) \Delta(m;a)\; da,
\end{equation}
provided the integral converges. 
Here 
\begin{equation}\label{eq:delta}
\Delta(m):=\prod_{\a \in \Sigma^+} (e^\a -e^{-\a})^{m_\a}.
 \end{equation}

The $\Theta$-spherical transform cointains important special cases:
\begin{enumerate}
\thmlist
\item
For $m=0$ the $\Theta$-spherical spherical function of spectral parameter 
$\l$ is 
$\pedtheta{\varphi}(m;\l,a)=\sum_{w \in \pedtheta{W}} e^{w\l(\log a)}$.  
The $\Theta$-spherical transform of a 
$\pedtheta{W}$-invariant function on $\pedtheta{A}$  
therefore coincides with its Euclidean Fourier transform. 
\item
If $\Theta=\Pi$, then the $\Theta$-spherical transform
agrees with the Opdam transform \cite{OpdActa}, and hence with the 
spherical transform of Harish-Chandra when the triple
$(\frak a,\Sigma,m)$ corresponds to a Riemannian symmetric space $G/K$
of noncompact type. The condition on the even multiplicities 
singles out the spaces $G/K$ with the property that all Cartan 
subalgebras in the Lie algebra $\frak g$ of $G$ are conjugate 
under the adjoint group of $\frak g$ (cf. \cite{He1}, p. 429; 
see also Section \ref{section:prelim} below).
\item 
When $(\frak a,\Sigma,m)$ corresponds to a Hermitian symmetric space $G/K$ 
and $\Theta$ is the set $\Pi_0$ of simple positive compact roots, then the 
$\Theta$-spherical transform coincides with the spherical Laplace transform on the
non-compactly causal space $G/H$ having $G/K$ as Riemannian dual. 
See Section \ref{section:prelim}. 
\end{enumerate}

\begin{Def}[Paley-Wiener space] \label{def:pwspacetheta}
Let $m \in \mathcal M^+$ be an even multiplicity function, 
and let $\Theta \subset \Pi$ be a fixed
set of positive simple roots. Let $C$ be a compact, convex and 
$\pedtheta{W}$-invariant subset of $\pedtheta{\frak a}$. 
The {\em Paley-Wiener space} $\pwthetamC$ is the space of all 
$\pedtheta{W}$-invariant meromorphic functions $g:\frakacs \to \C$ 
satisfying the following properties:
\begin{enumerate}
\renewcommand{\labelenumi}{\theenumi.}
\item  \label{item:pwspacethetauno}
$\pedtheta{e}^-(m;\l) g(\l)$ is a rapidly decreasing entire function of exponential type $C$,
that is 
for every $N \in \N$  
there is a constant $C_N\geq 0$ such that
$$
\abs{\pedtheta{e}^-(m;\l) g(\l)} \leq C_N (1+\abs{\l})^{-N} e^{q_C(\Re\l)}
$$
for all $\l \in \frakacs$.
\item \label{item:pwspacethetadue}
The function 
\begin{equation}\label{eq:Pav}
\Pavtheta g(\l):=\sum_{w \in \pedtheta{W}\setminus W} g(w\l)
\end{equation}
 extends to an entire function on $\frakacs$.
\end{enumerate}
\end{Def}

By $\pedtheta{W}$-invariance of $g$, Condition \ref{item:pwspacethetadue} in 
Definition \ref{def:pwspacetheta} is equivalent to
\begin{enumerate}
\setcounter{enumi}{1}
\renewcommand{\labelenumi}{\theenumi'.}
\item
The function 
\begin{equation}\label{eq:Pavtheta}
\Pav g(\l):=\frac{1}{\abs{W}}\; \sum_{w \in W} g(w\l)
\end{equation}
extends to an entire function on $\frakacs$.
\end{enumerate}

Condition \ref{item:pwspacethetadue} is automatically 
satisfied in the Euclidean case 
$m=0$, in the complex case $m=2$, and 
when $\Theta=\Pi$ (see Proposition \ref{prop:PWspace}). 
Notice that  $\pedtheta{e}^-\equiv 1$ in the Euclidean case and when
$\Theta=\Pi$. 
Therefore $\pedtheta{\PW}(C;0)\equiv \PW(C)^{\pedtheta{W}}$ and 
$\pedPi{\PW}(C;m) \equiv \PW(C)^W$.

A generalization of the classical Paley-Wiener theorem to the 
$\Theta$-spherical transform holds under certain restrictions
on $\Theta$ and $m$, which we indicate as Condition A. 
It requires either that there are ``not too many'' $\l$-singular 
hyperplanes (condition on $m$) or that the cone $\pedtheta{\frak a}$ is ``wide enough'' 
in $\frak a$ (condition on $\Theta$). 
We refer to Section \ref{section:condA}  for the 
precise statement of Condition A.  
Here we only observe that all pairs $(\Theta,m)$
with $\Theta=\Pi$ or $m=0,2$ or corresponding to a 
$K_\varepsilon$ symmetric space with even multiplicities satisfy 
this condition.  The $K_\e$ spaces are the relevant 
symmetric spaces for the geometric realization of the $\Theta$-spherical
transform, and include 
the non-compactly causal (NCC) symmetric spaces as special case. 
See Section \ref{section:prelim} for more information. 

The Paley-Wiener theorem for the $\Theta$-spherical
transform, which is the main result of this paper, is given by the following theorem. 
   
\begin{Thm}[Paley-Wiener theorem]\label{thm:pw}
Let $\Theta \subset \Pi$ and let $m \in \mathcal M^+$ be an even multiplicity 
function. Suppose Condition A is satisfied.
Let $C$ be a compact, convex and $\pedtheta{W}$-invariant subset of
$\pedtheta{\frak a}$.
Then the $\Theta$-spherical transform $\pedtheta{\mathcal F}(m)$ 
is a bijection of $C^\infty_c(C)^{\pedtheta{W}}$ onto $\pwthetamC$.
\end{Thm}

It is important to remark that Condition A only plays a role in the proof that the 
$\Theta$-spherical transform is surjective.
The 
injectivity will be proven without any 
restriction on the pair $(m,\Theta)$. See Theorem \ref{lemma:pwinto}.

In the case $m=0$, Theorem \ref{thm:pw} reduces to the classical Paley-Wiener
theorem for the Euclidean spherical transform. 
For $\Theta=\Pi$, it is a slight generalization of
the Paley-Wiener theorem for the Odpam transform 
as stated in Theorems 8.6 and 9.13(4) in \cite{OpdActa}, 
where it has been proven for $W$-invariant 
convex sets of $\frak a$ of the form 
$\conv (W( H))$ with $H \in \frak a$. Indeed in this case the Theorem
holds without assuming that the multiplicities must be even.  
For $\Theta=\Pi$ and $m$ geometric, it reduces to the celebrated theorem of 
Helgason-Gangolli-Rosenberg (see e.g. \cite{He2}, Ch. IV, Theorem 7.1, 
or \cite{GV}, Theorem 6.6.8). An elementary proof of Theorem \ref{thm:pw} 
for $\Theta=\Pi$ and $m$ even is given on p. \pageref{secondproof}.

If $(\frak a, \Sigma,m)$ corresponds to a NCC symmetric space 
$G/H$ and $\Pi=\Theta_0$, then 
Theorem \ref{thm:pw} yields a Paley-Wiener theorem for the spherical Laplace 
transform on $G/H$. For a very special exausting family of compact convex 
subsets in the maximal cone 
$\cmax:=\frak a_{\Pi_0}$, a Paley-Wiener type theorem in this context 
was also proven in \cite{AOS}. 
However, our description of the Paley-Wiener space in Definition 
\ref{def:pwspacetheta} is more explicit, 
and very close to the Euclidian one, since
the exponential growth condition is only modified by 
multiplication by  the polynomial 
$\pedtheta{e}^-$
giving the location of the possible singularities. 
The simplification of the Paley-Wiener space in 
the complex NCC case was also not remarked before. Finally, the proof 
presented in \cite{AOS} depends heavily on the causal structure of 
the symmetric space, 
whereas our proof here applies to a much more general context, and 
is based only on the root structure and the fact that the multiplicities 
are even.

More precisely, the proof of Theorem \ref{thm:pw} 
is based on a reduction to the 
Paley-Wiener theorem for the Euclidean Fourier transform. 
The tools allowing this reductions are $W$-invariant differential operators 
$D_m$ on $A$ linking the $\Theta$-spherical functions for even multiplicities 
to the exponential functions. 
These differential operators are essentially special cases 
of the shift operators of Opdam \cite{Opd88a}. 
Shift operators are differential operators
with smooth coefficients on the positive Weyl chamber $A^+$ but
generally singular on the ``walls'' of $A^+$. 
The crucial observation is that the modified shift operators $D_m$ 
we are  using have
smooth coefficients on all of $A$. 
In fact its coefficients are even holomorphic
on a certain torus $\complex{A}$ complexifying $A$
(see Theorem \ref{thm:Dm} and Corollary \ref{cor:shiftnonsing}). 

An important intermediate step in the proof of the Paley-Wiener theorem 
is a new inversion
formula for the $\Theta$-spherical transform in the even multiplicity case. 
The formula
resembles the inversion formula for the Euclidean Fourier transform, 
with only a correction 
due to the differential operator $D_m$. 
See Theorem \ref{cor:newinv} and Corollary 
\ref{cor:newinv}.

The study of the modified shift operators $D_m$ gives as a parallel result explicit formulas 
for the $\Theta$-spherical functions corresponding to even multiplicites. 
These formulas are of independent interest. 
They yield as special instances new explicit formulas for the
spherical functions on Riemannian and NCC symmetric spaces with 
even multiplicities. See Theorem \ref{thm:formulas}, 
Corollary \ref{cor:formulas}, 
and Examples \ref{ex:formularone} and 
\ref{ex:phithetacomplex}. We determine two types of formulas for the
$\Theta$-spherical functions. The first type involves 
the differential operator $D_m$ applied to a sum over 
the Weyl group $\pedtheta{W}$ of exponential functions. 
This type is used in the proof of the Paley-Wiener theorem.  
The second type presents a differential operator applied to an alternating
sum of exponential functions over the Weyl group $\pedtheta{W}$. 
The formulas of the second type resemble the classical formulas 
by Harish-Chandra for the
spherical functions on Riemannian symmetric spaces of noncompact type with 
a complex structure. 
The second type is obtained from the first one by splitting $D_m$
as a composition of a differential operator and the shift operator 
of constant shift 2.  

The proof that the operator $D_m$ have nonsingular coefficients 
depends on the behavior of the Harish-Chandra series 
and their derivatives in the spectral parameters on the 
walls of the positive Weyl chamber. 
This behavior can be deduced from new estimates, 
which hold for arbitrary positive 
multiplicity functions and which we have collected in Appendix A. 
The fact that the above mentioned splitting of $D_m$ returns a 
differential operator is stated in Corollary 
\ref{cor:pieceofDm}, which is proven in Appendix \ref{app:proofcor}.

\tableofcontents

\section{Preliminaries} \label{section:prelim}
Let $\frak a$ be an $l$-dimensional real Euclidean vector space with 
inner product $\inner{\cdot}{\cdot}$. 
For every $\a$ in the dual space $\frak a^*$
of $\frak a$, let $A_\a \in \frak a$ be determined by
$\a(H)=\inner{H}{A_\a}$ for all $H \in \frak a$. Then
for $\a\not= 0$ the vector
$H_\a:=2A_\a/ \inner{A_\a}{A_\a}$
satisfies $\a(H_\a)=2$. The assignment 
$\inner{\a}{\b}:=\inner{A_\a}{A_\b}$ defines an inner product in 
$\frak a^*$. Let $\Sigma$ be a reduced 
root system in  $\frak a^*$ with associated Weyl group $W$. 
For every $\a \in \Sigma$, 
we denote by $r_\a$ the reflection in $\frak a^*$ defined by
$r_\a(\l):=\l-\l(H_\a)\a$ for all $\l \in \frak a^*$. 

Let $\Sigma^+$ be a choice of positive roots in $\Sigma$ and 
$\Pi=\{\a_1,\dots,\a_l\}$ the  system of simple roots 
associated with $\Sigma^+$. The positive Weyl chamber 
$\frak a^+$ consists of the elements  $H\in\frak a$ for which 
$\a(H)>0$ for all $\a \in \Sigma^+$.  

The complexification $\liecomplex{a}:=\frak a \otimes_\R \C$ of $\frak a$ can 
be viewed as the Lie algebra of the complex torus 
$\complex{A}:= \liecomplex{a} / \Z\{i\pi  
H_\a: \a \in \Sigma\}$. 
The exponential map $\exp: \liecomplex{a} \rightarrow \complex{A}$ 
is the canonical projection of $\liecomplex{a}$ onto $\complex{A}$. Its
multi-valued inverse is denoted by $\log$. 
The split real form $A:=\exp \frak a $ of $\complex{A}$ is an 
abelian subgroup of $\complex{A}$ with Lie algebra $\frak a$ and
$\exp: \frak a \rightarrow A$ is a diffeomorphism. 
Set $A^+:=\exp \frak a^+$.  
The polar decomposition of $\complex{A}$ is 
$\complex{A}=AT$, 
where $T:=\exp(i \frak a)$ is a compact torus with Lie algebra $i\frak a$. 
Let $\frakacs$ be the space of all $\C$-linear functionals on $\frak a$.
The action of $W$ extends to $\frak a$ by duality, 
and then to $\frakacs$ and $\liecomplex{a}$ by $\C$-linearity, 
and to $\complex{A}$ and $A$ by 
the exponential map. Moreover, $W$ acts on functions $f$ on any of these spaces
by $(wf)(x):=f(w^{-1}x)$, $w \in W$.
The $\C$-bilinear extension to $\frakacs$ and 
$\liecomplex{a}$ of the inner products $\inner{\cdot}{\cdot}$ on
$\frak a^*$ and $\frak a$ will also be denoted by $\inner{\cdot}{\cdot}$.

A multiplicity function on $\Sigma$ is a $W$-invariant function 
$m:\Sigma\rightarrow \C$. Setting $m_\a:=m(\a)$, 
we therefore have $m_{w\a}=m_\a$ for all  $\a \in \Sigma$ and $w \in W$. 
The set ${\mathcal M}$ of all multiplicity functions on $\Sigma$ is a 
subspace of the finite-dimensional $\C$-vector space $\C^\Sigma$. 
Given $k \in \C$, the multiplicity function with $(km)_\a=k\ma$
for all $\a \in \Sigma$ is denoted by $km$. The multiplicity 
function with constant value $\ma=k$ for all $\a \in \Sigma$ is denoted by $k$.  
The multiplicity function $m$ is said to be positive (resp., even) if $\ma\geq 0$ 
(resp., $\ma\in 2\Z$) for all $\a \in \Sigma$.  
The set of positive multiplicity functions is denoted by ${\mathcal M}^+$.
Finally, a multiplicity function $m$ is said to be \emph{geometric} if there is a 
Riemannian symmetric space of noncompact type $G/K$ 
with restricted root system $\Sigma$ such that $m_\a$ is the multiplicity
of the root $\a$ for all $\a\in \Sigma$. Otherwise, $m$ is said to be 
\emph{non-geometric}. \footnote{
We adopt the multiplicity notation commonly used in the theory of 
symmetric spaces. 
It differs from the notation employed by Heckman and Opdam 
in the following ways.
The root system $R$ used by Heckman and Opdam
is related to our root system $\Sigma$ by the relation 
$R=\{2\a:\a \in \Sigma\}$; 
the multiplicity function $k$ in Heckman-Opdam's work is related to 
our $m$ by $k_{2\a}=\ma/2$.}

In the geometric case, even multiplicities correspond to Riemannian 
symmetric spaces $G/K$ with the property that all Cartan subalgebras 
in the Lie algebra $\frak g$  of $G$ are
conjugate under the adjoint group of $\frak g$ (see \cite{He1}, p. 429). 
The simplest examples correspond to spaces with a complex structure, 
for which all multiplicities are equal to $2$. 
Among the pseudo-Riemannian symmetric spaces, the material developed 
in this paper is particularly relevant for the study of spherical 
functions on the $K_{\varepsilon}$-spaces and the so-called
non-compactly causal (NCC) symmetric spaces (see e.g. \cite{HO}). 
We shall assume that $G$ is simple.  Let $\mathfrak{k}$ be
the Lie algebra of $K$ and denote by $\theta$ the corresponding Cartan 
involution. Denote by $\mathfrak{p}$ the $(-1)$-eigenspace of $\theta$. 
We assume that $\mathfrak{a}$ is a maximal
abelian subspace of $\mathfrak{p}$. 
Set $\mathfrak{m}=\mathfrak{z}_{\mathfrak{k}}(\mathfrak{a})$
and $\mathfrak{g}^{\alpha}=
\{X\in \mathfrak{g}^{\alpha}: \text{$[H,X]=\a(H) X$ for all 
$H \in \frak a$}\,\}$.
A map $\e:\Sigma \to \{1,-1\}$
is a {\em signature} if $\e(\alpha + \beta)=\e(\alpha)\e(\beta)$
whenever $\alpha,\beta, \alpha+\beta\in \Sigma$ and 
$\e(\a)=\e(-\a)$ for all $\a \in \Sigma$. 
Let $\e$ be a signature. Define an involution $\theta_\e$
on $\mathfrak{g}$ by
\[\theta_\e|_{\mathfrak{m}\oplus \mathfrak{a}}=
\theta|_{\mathfrak{m}\oplus \mathfrak{a}}\]
and for each $\alpha\in \Sigma$ by
\[\theta_\e|_{\mathfrak{g}^\alpha}=
\e(\alpha)\theta|_{\mathfrak{g}^\alpha}\, .\]
Denote also by $\theta_\e$ the corresponding involution on $G$, and
let $H=\{g\in G: \theta_\e(g)=g\}$.
Then $G/H$ is called a $K_\e$-space. 
If there exists a $X_0\in \mathfrak{a}$ such that 
$\alpha (X_0)\in \{0,1,-1\}$ for all $\alpha \in \Sigma$
and  the signature $\e$ is given by
\begin{equation} \label{eq:signatureNCC}
\e(\alpha)=(-1)^{\alpha (X_0)}
\end{equation}
then $G/H$ is a non-compactly causal symmetric space. 
We notice that the root system
and the multiplicity function corresponding to the
symmetric space $G/H$ is the same as that of $G/K$.
The infinitesimal classification of the 
$K_\e$-spaces $G/H$ with simple $G$ can be found in the
appendix of \cite{OS}.  
The root multiplicities can be read off from the list. We
ist all the  $K_\e$ symmetric pairs with even multiplicites in Appendix 
\ref{app:list}. An interesting fact, which however will not be used 
in the present paper, is that all multiplicities of a $K_\e$-space 
with even multiplicities are equal. 

For $\a \in \Sigma$ and $\l \in \frakacs$ we set 
\begin{equation}\label{eq:la}
  \la:=\frac{\l(H_\a)}{2}=\frac{\inner{\l}{\a}}{\inner{\a}{\a}}.
\end{equation}
The \emph{restricted weight lattice} of $\Sigma$ is the set $P$ of 
all $\l \in \frak a^*$ with $\la\in\Z$ 
for all $\a \in\Sigma$. Notice that $2\a \in P$ for all $\a \in \Sigma$.
The elements in $\frakacs \setminus P$ are said
to be \emph{generic}.

If $\l \in P$, then the exponential $e^\l$ defined by 
$e^\l(a):=e^{\l(\log a)}$ is single valued on $\complex{A}$.
 The $\C$-linear span of the $e^\l$, $\l\in P$,is the ring  
of regular functions $\C[\complex{A}]$ on the affine algebraic 
variety $\complex{A}$.
The lattice $P$ is $W$-invariant, and the Weyl group acts on
$\C[\complex{A}]$
according to $w(e^\l):=e^{w\l}$ for all $w \in W$. 
The set $\Hreg:=\{h \in \complex{A}: e^{2\a(\log h)}\neq 1 \ 
\text{for all $\a \in \Sigma$}\}$
consists of the regular points of $\complex{A}$ for the action of $W$. 
Notice that $A^+$ is a 
subset of $\Hreg$. The algebra $\C[\Hreg]$ of regular functions on $\Hreg$ 
is the subalgebra of the quotient field $\C(\complex{A})$ 
of $\C[\complex{A}]$ generated by
$\C[\complex{A}]$ and $1/(1-e^{-2\a})$, $\a \in \Sigma^+$.
Its $W$-invariant elements form the subalgebra $\C[\Hreg]^W$.

Given a multiplicity function $m \in \mathcal M$ on $\Sigma$, we define
\begin{align}
   \label{eq:rhom}
  \rho(m)&:=\frac{1}{2} \sum_{\a \in \Sigma^+} \ma \a, \\
  \label{eq:Deltam}
\Delta(m)&:=\prod_{\a \in \Sigma^+} (e^\a-e^{-\a})^{m_\a}
=e^{\rho(m)}\prod_{\a \in \Sigma^+} (1-e^{-2\a})^{m_\a}\, .
\end{align}
In particular, $\Delta:=\Delta(1)$ is the {\em Weyl denominator}.

\begin{Lemma} \label{lemma:deltainv}
Let $\Sigma$ be a reduced root system and let $m \in \mathcal M^+$ be an 
even multiplicity function. Then $\rho \in P$. Consequently, $\Delta(m/2) \in 
\polyAc$. Furthermore, $\Delta(m) \in\polyAc^W$.
\end{Lemma}
\pf
If $\a \in \Pi$ is a simple root, then $\rho_\a=\ma/2 \in \Z$. For general 
$\a \in\Sigma$ there are $\beta \in \Pi$ and $w \in W$ so that $\a=w\beta$.
The property that $\rho_\a\in \Z$ follows then by induction on the length of 
$w$. This proves $\rho \in P$. Since $-2\a \in P$ and $\ma/2 \in \Z$, we
also obtain
\begin{equation*}
  \Delta(m/2)=e^\rho \prod_{\a \in \Sigma^+} (1-e^{-2\a})^{m_\a/2} \in \polyAc.
\end{equation*}
Let $c=\frac{1}{2}\sum_{\alpha \in \Sigma^+}m_\alpha$. Then
\begin{equation*}
 \Delta(m)=(-1)^c \prod_{\a \in \Sigma^+} (e^\a-e^{-\a})^{m_\a/2}
   \prod_{\a \in \Sigma^+} (e^{-\a}-e^{\a})^{m_\a/2}
  =(-1)^{c}  \prod_{\a \in \Sigma} (e^\a-e^{-\a})^{m_\a/2},
\end{equation*}
which proves that $\Delta(m) \in \polyAc^W$.
\qed

Let $\polya$ denote the symmetric algebra over $\liecomplex{a}$ 
considered as the space of polynomial functions on $\frakacs$, 
and let $\polya^W$ be the subalgebra of $W$-invariant elements.
For $p\in \polya$ write $\partial (p)$ for the corresponding
 constant-coefficient differential operator
on $\complex{A}$
(or on $\liecomplex{a}$).
Let $\D(\Hreg):=\C[\Hreg]\otimes \polya$ 
denote the algebra of differential operators on 
$\complex{A}$ with coefficients 
in $\C[\Hreg]$. 
The Weyl group $W$ acts on $\D(\Hreg)$ according to 
\begin{equation*}
  w\big(\phi\otimes \partial(p)\big):=w\phi \otimes \partial(wp).
\end{equation*}
Let $\D(\Hreg)^W$ denote the subspace of $W$-invariant elements in
$\D(\Hreg)$. The set 
$\D(\Hreg) \otimes \C[W]$ 
of differential reflection operators on $\Hreg$ 
can be endowed with the structure of an 
associative algebra 
with respect to the product
\begin{equation*}
 (D_1 \otimes w_1)\cdot (D_2 \otimes w_2)=D_1w_1(D_2) \otimes w_1w_2,
\end{equation*}
where the action of $W$ on differential operators is defined by 
$(wD)(wf):=w(Df)$
for every sufficiently differentiable function $f$. 
It is also a left $\C[\Hreg]$-module.
Considering $D \in \D(\Hreg)$ as 
element of $\D(\Hreg) \otimes \C[W]$, we shall usually write 
$D$ instead of $D \otimes
1$. 
The differential-reflection operators act on
functions $f$ on $\Hreg$ according to $(D\otimes w)f:=D(wf)$. 

Define a linear map $\Upsilon:\D(\Hreg)\otimes \C[W] \rightarrow \D(\Hreg)$ by 
\begin{equation*}
  \Upsilon(\sum_j D_j \otimes w_j):=\sum_j D_j.
\end{equation*}
Then $\Upsilon(Q)f=Qf$ for all  $Q \in \D(\Hreg) \otimes \C[W]$ and all 
$W$-invariant $f$ on $\Hreg$.

 \begin{Def}  {\textrm (\cite{CherInv})} \label{def:cherednik}
Let $m \in {\mathcal M}$ and $H \in \liecomplex{a}$.
The \emph{Dunkl-Cherednik operator\/}  $T(m;H) \in \D(\Hreg) \otimes \C[W]$ 
is defined by 
\begin{equation*}
  T(m;H):=\partial(H)-\rho(m)(H)+\sum_{\a\in \Sigma^+} m_\a \a(H) 
(1-e^{-2\a})^{-1}
\otimes (1-r_\a).
\end{equation*}
 \end{Def}

The  Dunkl-Cherednik operators map $\C[\complex{A}]$ into itself, 
but they can also be considered 
as operators acting on other function spaces, 
for instance, on the spaces of smooth and compactly supported smooth
functions on $A$ and $\frak a$. 
Indeed, as can be seen from the Taylor formula, 
the term $1-r_\a$ cancels the apparent singularity on $A$ and $\frak a$
arising from the denominator $1-e^{-2\a}$.

The Dunkl-Cherednik operators
$\{T(m;H): H \in \liecomplex{a}\}\,$ form a 
commuting family of differential-reflection operators in 
$\D(\Hreg) \otimes \C[W]$
(cf. \cite{OpdActa}, Section 2).
Therefore the map $H \mapsto T(m;H)$ on $\liecomplex{a}$
extends uniquely to an algebra homomorphism of $\polya$ into
$\D(\Hreg) \otimes \C[W]$.
For $p \in \polya$ we set 
$D(m;p):=\Upsilon\big(T(m;p)\big)$.   
Then $\Upsilon$ establishes an algebra homomorphism of 
$\{T(m;p): p \in \polya^W\}$ into $\D(\Hreg)^W$ (see \cite{HS}, Lemma 1.2.2).
The differential operator $D(m;p)$ is the unique 
element of $\D(\Hreg)^W$ coinciding on $\C[\Hreg]^W$ with the
restriction of $T(m;p)$. The algebra
\begin{equation} \label{eq:diffops}
\D(\frak a, \Sigma,m):=\{D(m;p):p \in \polya^W\}
\end{equation}
is a commutative subalgebra of $\D(\Hreg)^W$.
It is called the algebra of
hypergeometric differential operators associated with the data 
$(\frak a,\Sigma,m)$.

Let $L_A$ denote the Laplace operator on $A$ and let $p_L \in  \polya^W$ 
be the polynomial defined by 
$p_L(\l):=\inner{\l}{\l}$ for $\l \in \frakacs$.
Then
\begin{equation} \label{eq:ML}
   ML(m):=D(m;p_L) =L(m)+ \inner{\rho(m)}{\rho(m)},
\end{equation}
  where 
  \begin{equation}
    \label{eq:Laplm}
  L(m):=L_A+\sum_{\a \in \Sigma^+} m_\a \,\frac{1+e^{-2\a}}{1-e^{-2\a}} \;
     \partial(A_\a)  
  \end{equation}
generalizes to arbitrary multiplicity functions $m$ 
the radial component on $A^+$ 
of the Laplace operator on a Riemannian symmetric space
$G/K$ of noncompact type.
The algebra $\D(\frak a, \Sigma,m)$
coincides with the commutant 
$\{Q \in \D(\Hreg)^W: L(m)Q=QL(m)\}$
of $L(m)$ in $\D(\Hreg)^W$.
It is therefore the analog, 
for arbitrary multiplicity functions, of the commutative 
algebra of the radial parts on $A$ of the 
invariant differential operators on a Riemannian symmetric space of  
noncompact type.
The map $\gamma(m): \D(\frak a, \Sigma,m)\rightarrow \polya^W$ defined by
\begin{equation}
  \label{eq:HChomo}
\gamma(m)\big(D(m;p)\big)(\l):=p(\l)  
\end{equation}
is an algebra isomorphism, called the Harish-Chandra homomorphism 
(see \cite{HS}, Theorem 1.3.12 and Remark 1.3.14). 
Chevalley's theorem implies that
$\D(\frak a,\Sigma,m)$ is generated by $l (=\dim \frak a)$ elements.

Let $\l \in \frakacs$ be fixed. The system of differential equations
\begin{equation}
  \label{eq:hypereq}
  D(m;p) \varphi=p(\l)\varphi, \qquad p \in \polya^W,
\end{equation}
is called the hypergeometric system of differential equations 
with spectral parameter $\l$ associated with the data $(\frak a,\Sigma,m)$. 
For geometric multiplicities, 
the hypergeometric system (\ref{eq:hypereq}) agrees
with the system of differential equations on $A$ defining 
Harish-Chandra's spherical function
of spectral parameter $\l$.

In the following we denote by $\N$ the set of positive integers and
set $\N_0=\N \cup \{0\}$.

\section{$\Theta$-spherical functions}
\label{section:HC}

As in the theory of spherical functions on Riemannian 
symmetric spaces of noncompact type, Heckman and Opdam \cite{HOpd1} 
looked for solutions of the hypergeometric system  
(\ref{eq:hypereq})
with spectral parameter $\l$ which are of the form
\begin{equation} \label{eq:HCexp}
\Phi(m;\l,a)=e^{(\l-\rho)(\log a)} 
\sum_{\mu \in 2\Lambda} \Gamma_\mu(m;\l) e^{-\mu(\log a)}
\end{equation}
for all $a$ in some tubular neighborhood of
$A^+$ in $\complex{A}$.
Here $\Lambda:=\left\{\sum_{j=1}^l n_j \a_j: n_j \in \N_0\right\}$ 
is the positive semigroup generated by the fundamental system of simple 
roots $\Pi:=\{\a_1,\dots, \a_l\}$ in $\Sigma^+$. 
For $\mu \in 2\Lambda \setminus \{0\}$, 
the coefficients $\Gamma_\mu(m;\l)$ are rational functions of 
$\l \in \frakacs$ determined from the recursion relations  
\begin{equation}\label{eq:recursion}
\inner{\mu}{\mu-2\l} \Gamma_\mu(m;\l) = 2 \sum_{\a\in \Sigma^+} m_\a 
 \sum_{\substack{k\in \N\\\mu -2k\a \in \Lambda}}
 \Gamma_{\mu-2k\a}(m;\l)  \inner{\mu+\rho-2k\a-\l}{\a}, 
\end{equation}
which are derived by formally inserting the series for $\Phi$ 
into the differential equation of (\ref{eq:hypereq}) corresponding to
$p=p_L$. With the initial condition $\Gamma_0(m;\l)=1$, the relations
(\ref{eq:recursion}) admit unique solutions $\Gamma_\mu(m;\l)$ provided
$\inner{\mu}{\mu-2\l}\neq 0$ for all $\mu \in  2\Lambda\setminus \{0\}$. 
The functions $\Phi(m;\l,a)$ are commonly known as the 
Harish-Chandra series.

\begin{Thm} \label{thm:HC}
\begin{enumerate}
\thmlist
\item {\rm (\cite{Opd88b}, Corollary 2.3; see also \cite{OpdGauss}, Lemma 2.1)}
There is a connected and simply connected open subset $U$ of $T$ containing 
the identity element $e$ such that $\Phi(m;\l,a)$ is a meromorphic functions of $(m,\l,a) \in 
\mathcal M \times \frakacs \times A^+U$ with at most simple poles along 
hyperplanes of the form $\mathcal M \times \mathcal H_{n,\a} \times A^+U$,
 where
$$
\mathcal H_{n,\a} :=\{\l \in \frakacs: \l_\a=n\}
$$
is a complex hyperplane in $\frakacs$ corresponding to some $\a \in \Sigma^+$
and $n \in \N$.
 \item {\rm (\cite{HS}, Corollary 4.2.6)}
For $\l \in \frakacs \setminus P$ 
the set $\{\Phi(m;w\l,a):w\in W\}$ is a basis for the
solution space of (\ref{eq:hypereq}) on $A^+U$.
\end{enumerate}
\end{Thm}

The set $U$ in Theorem \ref{thm:HC} is chosen so that the function $\log$ is 
single valued on it, i.e. so that $e^{(\l-\rho(m))(\log a)}$ 
is holomorphic on $U$ for all $\l \in \frakacs$. 
Observe that there is a neighborhood $V$ of $\l=0$ in $\frakacs$ such that
$\Phi(m;\l,a)$ is a holomorphic function of 
$(m,\l,a) \in \mathcal M \times V \times A^+U$. 

Let $\Theta \subset \Pi$ be an arbitrary set of positive simple roots, and let 
$\rootstheta=\Sigma \cap \,\Span \Theta$ be the subsystem of $\Sigma$ generated by 
$\Theta$. Then the Weyl group $\pedtheta{W}$ of $\rootstheta$ is the subgroup of
$W$ generated by the reflections $r_\a$ with $\a \in \Theta$.

Given an even multiplicity function $m\in\mathcal M^+$, we define
\begin{align*} 
\cthetampl:&=\prod_{\a \in \rootstheta^+} \prod_{k=0}^{\ma/2-1} 
\frac{1}{\la+k},\\
\cthetamml:&= \prod_{\a \in \Sigma^+\setminus \rootstheta^+} 
 \prod_{k=-\ma/2+1}^0 \frac{1}{-\la+k}\\
&=(-1)^{d(\Theta,m)}  \prod_{\a \in \Sigma^+\setminus \rootstheta^+} 
 \prod_{k=0}^{\ma/2-1} \frac{1}{\la+k}\,,
\end{align*}
where 
\begin{equation}
  \label{eq:dmtheta}
d(\Theta,m):=\frac{1}{2} \sum_{\a \in \Sigma^+\setminus \rootstheta^+} \ma\,.  
\end{equation}

\begin{Def} 
Suppose $\l \in \frakacs \setminus P$. Then the {\em $\Theta$-spherical function of 
spectral parameter $\l$} is the $\pedtheta{W}$-invariant 
solution of  (\ref{eq:hypereq}) defined on $A^+U$ by
\begin{equation}
  \label{eq:thetaphi}
  \pedtheta{\varphi}(m;\l,a):=\pedtheta{c}^-(m;\l) \; 
\sum_{w \in \pedtheta{W}}
\pedtheta{c}^+(m;w\l) \Phi(m;w\l,a).
\end{equation} 
\end{Def}
 
\begin{Rem}
The functions $\pedtheta{c}^+$ and $\pedtheta{c}^-$ are respectively the analog
for $\rootstheta$ of Harish-Chandra $c$-function for Riemannian symmetric 
spaces and of the $c$-function of \cite{KO} for NCC spaces.  
We refer the reader to  \cite{PHab} or \cite{P1}  
for their general definition and for further information on $\Theta$-spherical 
functions corresponding to arbitrary multiplicity functions. 
\end{Rem}

\begin{Ex}
When $\Theta=\Pi$, the function 
$\pedPi{\varphi}(m;\l,a)/\pedPi{c}^+(m;\rho(m))$ coincides with the
hypergeometric function of spectral parameter $\l$ as defined by Heckman and 
Opdam \cite{HOpd1}. 
For geometric multiplicity functions, it therefore agrees with the 
restriction of Harish-Chandra's spherical function  of spectral parameter $\l$
to a maximal flat subspace $A$ of the corresponding Riemannian symmetric 
space $G/K$.
\end{Ex}

\begin{Ex}
Suppose $\Sigma$ is the restricted root system of a NCC symmetric space $G/H$. 
If $\Theta=\Pi_0$
is the set of positive compact simple roots and $m$ is geometric, 
then the function
$\ped{\Pi_0}{\varphi}(m;\l,a)/\ped{\Pi_0}{c}^+(m;\rho(m))
\ped{\Pi_0}{c}^-(m;\rho(m))$ coincides 
with the spherical function on $G/H$ with spectral parameter $\l$ as 
defined by \cite{FHO} (see \cite{O}).
\end{Ex}

Set 
\begin{align}
 \pedtheta{\frak a}&:=\big(\pedtheta{W}( \overline{\frak a^+})\big)^0=
\{H \in \frak a: \text{$\a(H)>0$ for all $\a \in \Sigma^+ \setminus 
\rootstheta^+$}\}, 
\label{eq:atheta}\\
 \pedtheta{A}&:=\big(\pedtheta{W}( \overline{A^+})\big)^0=
\exp(\pedtheta{\frak a}).
\label{eq:Atheta} 
\end{align}
The regularity properties of the $\Theta$-spherical functions for even 
multiplicity functions are collected in the following theorem.

\begin{Thm}{\rm (\cite{OP2}, Theorem 10)}
\label{thm:sphentire} Define
\begin{equation} \label{eq:eminus} 
\pedtheta{e}^-(m;\l):=\prod_{\a \in \Sigma^+\setminus \rootstheta^+} \;
 \prod_{k=-\ma/2+1}^{\ma/2-1} (\la-k).
\end{equation}
Then there is a $\pedtheta{W}$-invariant tubular neighborhood 
$\pedtheta{U}$ in $\complex{A}$
of $\pedtheta{A}$ such that the function
\begin{equation*} 
\pedtheta{e}^-(m;\l)\;\pedtheta{\varphi}(m;\l,h)
\end{equation*}
extends as a $\pedtheta{W}$-invariant holomorphic function of
$(\l,h) \in \frakacs \times \pedtheta{U}$.
\end{Thm}

Estimates for the $\Theta$-spherical functions and their derivatives 
for $\Theta=\Pi$ and $m \in \mathcal M^+$ arbitrary 
have been proven by Opdam. 

\begin{Lemma}[cf. \cite{OpdActa}, Theorem 3.15, Proposition 6.1 
and Corollary 6.2] 
\label{lemma:estsphPi}
Let $m\in \mathcal M^+$ be a fixed multiplicity function. 
For $a\in \complex{A}$ write $\log a=a_{\rm R} + i a_{\rm I}$ with 
$a_{\rm R}, a_{\rm I} \in \frak a$.
 \begin{enumerate}
\thmlist
  \item
For all $\l \in \frakacs$ and all $a \in \complex{A}$ satisfying  
$\abs{\a(a_{\rm I})}<\pi/2$ for every $\a \in \Sigma$
\begin{equation*}
  \abs{\pedPi{\varphi}(m;\l,a)} \leq \abs{W}^{1/2} \pedPi{c}^+(m;\rho(m))\;
e^{-\min_{w\in W} \Im w\l(a_{\rm I})+\max_{w\in W} w\rho(a_{\rm I})+ 
\max_{w\in W} \Re w\l(a_{\rm R})} 
\end{equation*}
In particular, for all $a \in A$ und $\l \in \frakacs$ 
\begin{equation} \label{eq:estsphPi}
  \abs{\pedPi{\varphi}(m;\l,a)} \leq \abs{W}^{1/2} \pedPi{c}^+(m;\rho(m))\;
e^{\max_{w\in W} \Re w\l(\log a)} 
\end{equation}
\item
Let $I=(\iota_1,\dots,\iota_l)$ be a multi-index and let 
$\partial_a^I=\partial(H_1)^{\iota_1} \cdots \partial(H_l)^{\iota_l}$ 
be the corresponding partial
differential operator associated with an orthonormal basis 
$\{H_1,\dots,H_l\}$ of
$\frak a$. 
Then, for every $\e>0$, there is a constant 
$C_{I,\e}>0$ such that for all $a \in A$
and $\l \in i\frak a^*$ with $\abs{\l}>\e$
\begin{equation*}
\abs{\partial_a^I \pedPi{\varphi}(m;\l,a)} \leq C_{I,\e} \abs{\l}^{\abs{I}+l},
\end{equation*}
where $\abs{I}:=\sum_{k=1}^l \iota_k$. In particular, if $K \subset A$ is 
compact, then there is a constant $C_{I,K}>0$ so that
for all $a \in K$ and $\l \in i\frak a^*$
\begin{equation} \label{eq:estsphPider}
 \abs{\partial_a^I \pedPi{\varphi}(m;\l,a)} 
\leq C_{I,K} (1+\abs{\l})^{\abs{I}+l}.
\end{equation}
  \end{enumerate}
\end{Lemma}

The definition of the $\Theta$-spherical functions and 
analytic continuation yield the following important functional relation.

\begin{Lemma}  \label{lemma:functionaleq}
Let $m \in \mathcal M^+$ be even, and let $d(\Theta,m)$ be as in 
(\ref{eq:dmtheta}).
Then there is a $W$-invariant neighborhood
of $\pedtheta{A}$ in $\complex{A}$ on which the relation
  \begin{equation*}
  \pedPi{\varphi}(m;\l,a)=(-1)^{d(\Theta,m)} \sum_{\pedtheta{W}\setminus W} 
  \pedtheta{\varphi}(m;w\l,a)
    \end{equation*}
holds as equality of meromorphic functions on $\frakacs$.
\end{Lemma}

We shall prove in Section \ref{section:formulas} that, 
in the case of  even multiplicity functions,
explicit global formulas for the $\Theta$-spherical functions can be 
obtained by means of Opdam's shift operators.

\section{Harish-Chandra series and Opdam's shift operators}
\label{section:shift}

Let $\C_\Delta[\complex{A}]=\cup_{k\in\Z} \, \Delta^k  \C[\complex{A}]$ 
denote the localization of $\C[\complex{A}]$ along the Weyl denominator
$\Delta=\prod_{\a \in \Sigma^+} (e^\a-e^{-\a}) \in \C[\complex{A}]$.
Then $\C[\mathcal M]\otimes\C_\Delta[\complex{A}]\otimes \polya$ 
is the algebra of 
the differential operators on $\complex{A}$ (or $\liecomplex{a}$) 
with coefficients in $\C_\Delta[\complex{A}]$ 
depending also polynomially on the multiplicities. 
Recall from (\ref{eq:ML}) the notation 
$ML(m):=L(m)+\inner{\rho(m)}{\rho(m)}$.

\begin{Def}{\rm (\cite{Opd88a}; see also \cite{HS}, Chapter 3)}  
\label{def:shiftop}
Let $l \in \mathcal M$ be an even multiplicity function.
A \emph{shift operator with shift} $l$ is a 
 differential operator $D(l)\in 
\C[\mathcal M]\otimes\C_\Delta[\complex{A}]\otimes \polya$ 
satisfying the following properties:
\begin{enumerate}
\condlist
\item \label{item:commshift}
For all $m \in \mathcal M$
\begin{equation*}
D(l;m) \circ ML(m)=ML(m+l) \circ D(l;m).
\end{equation*}
\item \label{item:asshift}
The differential operator $D(l)$ admits on $A^+$ an expansion 
of the form
\begin{equation} \label{eq:aseq}
 D(l)=\sum_{\mu \in 2\Lambda} e^{-\rho(l)-\mu} \partial(p_\mu) 
\end{equation}
with $p_\mu \in \C[\mathcal M] \otimes \polya$

\end{enumerate}
\end{Def}

\begin{Rem} \label{rem:shiftcomp}
In \ref{def:shiftop} (\ref{item:commshift}) we have set $D(l;m):=D(l)(m)$, 
where $m$ denotes the variable in
$\mathcal M$ from which $D(l)$ depends polynomially. The asymptotics of
$D(l)$ in (\ref{item:asshift}) 
correspond to the power series expansion 
$(e^\a-e^{-\a})^{-1}=\sum_{k=1}^\infty e^{-2k\a}$ on $A^+$. 
Moreover, the elements 
of $\C[\mathcal M]$ act in  (\ref{eq:aseq}) as constants, 
i.e. $\partial(q\otimes p):=q \; \partial(p)$ for
$q \in  \C[\mathcal M]$ and $p \in \polya$.
\end{Rem}

All shift operators turn out to be $W$-invariant.  
See e.g. \cite{HS}, Corollary 3.1.4. 
The name ``shift operator of shift $l$'' reflects the property of $D(l)$ of
relating Harish-Chandra series corresponding to 
multiplicity functions differing by the shift $l$. 
This is the content of the next Theorem.
.

In the following we denote by $D^*$ the 
formal transpose of a differential operator $D$ on $\Hreg$ 
with respect to the Haar measure $da$ on $A$:
$$
\int_A (Df)g\; da = \int_A f (D^*g)\; da
$$ 
for all smooth functions $f,g$ on 
$A^{\rm reg}:=\Hreg \cap A$ with at least one
of $f,g$ with compact support.

\begin{Thm} \label{thm:shiftops}
 {\rm (\cite{Opd88a})}  
\label{thm:exshift}
Let $l \in \mathcal M^+$ be an even multiplicity function.
\begin{enumerate}
\thmlist
\item
There exists a unique shift operator $G_-(-l)$ of shift $-l$ such that  
\begin{equation*}
  G_-(-l;m)\Phi(m;\l,a)=
\frac{\pedPi{c}^+(m-l;\l)}{\pedPi{c}^+(m;\l)} \; \Phi(m-l;\l,a)
\end{equation*}
 for all $(m,\l,a) \in \mathcal M \times (\frakacs \setminus P) \times A^+$.
\item
Let $\Delta(m)$ be as in (\ref{eq:Deltam}). Then 
\begin{equation}
  \label{eq:gpluseminus}
 G_+(l;m):=\Delta(-l-m)\circ G^*_-(-l;m+l)\circ \Delta(m)
\end{equation}
defines the unique shift operator $G_+(l)$ of shift $l$ such that
\begin{equation*}
   G_+(l;m)\Phi(m;\l,a)
=\frac{\pedPi{c}^+(m;-\l)}{\pedPi{c}^+(m+l;-\l)} \; \Phi(m+l;\l,a)
\end{equation*}
 for all $(m,\l,a) \in \mathcal M \times (\frakacs \setminus P) \times A^+$.
\end{enumerate} 
\end{Thm}
\pf
The existence of $G_-(-l)$ and $G_+(l)$ is proven in 
\cite{HS}, Theorem 3.4.3 and Corollary 3.4. 
Their uniqueness depends on the injectivity of the Harish-Chandra mapping for
shift operators ({\em loc. cit.}, Proposition 3.1.6; see also Theorems
3.3.6 and 3.3.7, and the remark after the proof of Corollary 3.4.4).
\qed
 
\begin{Rem}
Obviously $G_{\pm}(0;m)=\id$ for all $m \in \mathcal M$.  
\end{Rem}

The multiplicity $m\equiv 0$ corresponds the Euclidean case with 
$\Phi(0;\l,a)=e^{\l(\log a)}$. Hence
Theorem  \ref{thm:exshift} provides a link between Harish-Chandra functions
for even multiplicities and exponential functions.

\begin{Cor}
  Let $m \in \mathcal M$ be even. Then the differential operators in 
$\C_\Delta[\complex{A}]\otimes \polya$
\begin{align}
D_-(m)&:=G_-(-m;m), \label{eq:Dminus}\\
D_+(m)&:=G_+(m;0)=\Delta(-m) \circ  G^*_-(-m;m)=\Delta(-m) \circ D^*_-(m) 
                               \label{eq:Dplus}
\end{align}
satisfy for all $(\l,a) \in (\frakacs\setminus P) \times A^+$ 
\begin{align}
D_-(m)\Phi(m;\l,a)&=\frac{1}{\pedPi{c}^+(m;\l)} \; e^{\l(\log a)},  
\label{eq:DminusHC}\\
D_+(m)e^{\l(\log a)}&= \frac{1}{\pedPi{c}^+(m;-\l)} \; \Phi(m;\l,a).    
\label{eq:DplusHC}
\end{align}
\end{Cor}

\begin{Rem}
\label{rem:compoG}
Suppose $m \in \mathcal M$ is even.
The uniqueness of the shift operators $G_\pm$
implies for even $l',l'' \in \mathcal M^+$ 
\begin{align*}
  G_-(-l'-l'';m)&=G_-(-l';m-l'') \circ G_-(-l'';m),\\
  G_+(l'+l'';m)&=G_+(l';m+l'') \circ G_+(l'';m).
\end{align*}
In particular, suppose $m \in \mathcal M$ is even with $m_\a\geq 2$
for all $\a$. Then
\begin{align}
 D_-(m)&:= G_-(-m;m)= G_-(-2;2)\circ G_-(2-m;m)=D_-(2)\circ G_-(-(m-2);m),
\label{eq:splitshiftminus} \\
D_+(m)&:= G_+(m;0) = G_+(m-2;2)\circ G_+(2;0) =G_+(m-2;2)\circ D_+(2).
\label{eq:splitshiftplus}  
\end{align}
Furthermore, shift operators of even constant shift $m \in \mathcal M$ 
can be obtained as composition of fundamental shift operators
corresponding to shifts $\pm 2$. For instance, when 
$\ma=m \in 2\N$ for all $\a \in \Sigma$, we have
\begin{align*}
G_-(-m;m)&=G_-(-2;2) \circ \dots \circ G_-(-2;m-2) \circ G_-(-2;m), \\
G_+(0;m)&=G_+(2;m-2) \circ \dots \circ G_+(2;2) \circ G_+(2;0). 
\end{align*}
Explicit formulas for $D_-(2)=G_-(-2;2)$ and $D_+(2)=G_+(2;0)$ 
are given in Example \ref{ex:shiftcomplex} below.
\end{Rem}

\begin{Ex}[The rank-one case]
The \emph{rank-one case} corresponds to triples $(\frak a,\Sigma,m)$ in which 
$\frak a$ is one dimensional. Suppose $\Sigma=\{\pm \a\}$ is of type $A_1$.
Notice that $P=\Z \a$ in this case. 
We identify $\liecomplex{a}$ and $\frakacs$ with $\C$ by setting $\l\a\equiv \l$  and 
$zH_\a/2\equiv z$ for $\l,z \in \C$. A multiplicity function can be identified with a complex 
number $m \in \C$. It is even if and only if $m \in 2\N_0$. We have then $\rho(m)= m/2$. 
The exponential function maps $\liecomplex{a}\equiv \C$ onto 
$\complex{A}\equiv \C^\times$, and 
$\Hreg=\C\setminus\{0,\pm1\}$. The Weyl chamber $\frak a^+$ coincides 
with the half-line $(0,\infty)$. 
The Weyl group $W$ reduces to $\{1,-1\}$ and acts on $\liecomplex{a}$ 
by multiplication. We normalize the 
inner product so that $\inner{\a}{\a}=1$.
With the identification of $\C[\complex{A}]=\C[u,u^{-1}]$, 
the Weyl denominator is given by $\Delta(u)=u-u^{-1}$. 
Set $\theta:=u \frac{d}{du}$.  Fix $u^{-1}du$ as 
Haar measure on $A\equiv  \R\setminus \{0\}$. 
Then (cf. \cite{HS}, \S 3.3)
\begin{align*}
  ML(m)&=\theta^2 + m\; \frac{1-u^{-2}}{1+u^{-2}}\;\theta+\rho(m)^2,\\
  G_-(-2;m)&=\Delta(u)\theta+(m-1)(u+u^{-1}),\\
  G_+(2;m)&=-\frac{1}{\Delta(u)}\;\theta \qquad \text{(independent of $m$)}. 
\end{align*}
Seen as differential operators on $\liecomplex{a}\equiv \C$, we have with 
$\Delta(z):=e^z-e^{-z}$,
\begin{align*}
ML(m)&=\frac{d^2}{dz^2}+ m\; \frac{1-e^{-2z}}{1+e^{-2z}}\;\frac{d}{dz}+
\rho(m)^2,\\
G_-(-2;m)&=\Delta(z)\frac{d}{dz}+(m-1)(e^z+e^{-z}), \\
G_+(2;m)&=-\frac{1}{\Delta(z)}\;\frac{d}{dz}.  
\end{align*}
Moreover $G_\pm(0;m)=\id$. For an arbitrary $l \in 2\N$, the shift operators 
$G_\pm(\pm l)$ can be obtained by iteration 
according to Remark (\ref{rem:shiftcomp}).
In particular, 
$$D_+(m)=G_+(m;0)=G_+(2;0)^{m/2}.$$
\end{Ex}

\begin{Ex}[The complex case] \label{ex:shiftcomplex}
The  complex case is characterized by the condition $\ma=2$ 
for all $\a \in \Sigma$. As differential operators on $\complex{A}$ 
or $\liecomplex{a}$ we have 
\begin{alignat*}{2}
D_-(2)&=G_-(-2;2)&&=\sigma \Big(\prod_{\a \in \Sigma^+} 
\partial(A_\a)\Big) \circ \Delta,\\
D_+(2)&=G_+(2;0)&&=\widetilde{\sigma} \Delta^{-1} 
\prod_{\a \in \Sigma^+} \partial(A_\a),  
\end{alignat*}
where ${\sigma}^{-1}:=\prod_{\a \in \Sigma^+}\inner{\a}{\a}$,
$\widetilde{\sigma}:=(-1)^{\abs{\Sigma^+}} {\sigma}$, and, 
as before, $\Delta=\Delta(1)$. 
\end{Ex}

As proven by Heckam and Opdam, arbitrary shift operators can be computed using 
Cherednik operators. 

\begin{Prop}[\cite{HeckBou}, Definition 4.2 and Proposition 4.4; 
\cite{Opd00}, Theorem 5.13]
Let $T(m;p)$ be the Cherednik operator associated with $p \in \polya$ 
and write
\begin{equation*}
  T(m;p)=\sum_{w\in W} D_w(m;p) \otimes w
\end{equation*}
with $D_w(m;p) \in \D(\Hreg)$. 
Let $q(m) \in \polya$ be defined by
\begin{equation*} \label{eq:pml}
q(m;\l):=\prod_{\a\in \Sigma^+} \Big(\l_\a +\frac{\ma}{2}\Big)  
\end{equation*} 
Then the fundamental shift operators $G_+(2;m), G_-(2;m) \in \D(\Hreg)$ are given by
\begin{align*}
  G_+(2;m)&=\Delta^{-1} \sum_{w\in W} D_w(m;q(m)), \\
  G_-(-2;m+2)&=\abs{W}^{-1} \sum_{w,v \in W} \varepsilon(w) \, v\big(D_w(m;q(m))\Delta\big),
\end{align*}
where $\varepsilon:W\to \{\pm 1\}$ is the sign character. 
In particular, for $f \in \C[\complex{A}]^W$ we have
\begin{align*}
  G_+(2;m)f&=\Delta^{-1} \sum_{w\in W} T(m;q(m))f, \\
  G_-(-2;m+2)f&=\abs{W}^{-1} \sum_{w,v \in W}  v\big(T(m;q(m))\Delta\big)f.
\end{align*}
\end{Prop}
 
The remaining of this section is devoted to the proof of the following theorem,
which is our first main result. It shows that multiplication of the
shift operators $D_\pm(m)$ by $\Delta(m)$ yields differential operators with 
holomorphic coefficients on the entire $\complex{A}$. 

\begin{Thm} \label{thm:Dm}
Set  
\begin{equation} \label{eq:Dm}
  D_m:=\Delta(m) D_+(m)  \qquad \big(=G^*_-(-m;m)=D_-(m)^*\big).
\end{equation}
Then $D_m$
is a $W$-invariant element of $\polyAc\otimes \polya$.
All the coefficients of $D_m$ vanish on $\bigcup_{\a \in \Sigma^+}
\{ a \in A: \a(\log a)=0\}$.  
\end{Thm}

For the proof of Theorem \ref{thm:Dm} we need to introduce some notation.
Let $\{H_1,\dots,H_l\}$ be a fixed orthonormal basis in $\frak a$ 
with dual basis $\{\xi_1,\dots, \xi_l\}$ for $\frak a^*$.
The coordinates of $\l \in \frakacs$ with respect to $\{\xi_1,\dots, \xi_l\}$ 
are $(\l_1,\dots,\l_l)$ with $\l_j:=\l(H_j)$. 
For a multi-index $I=(\iota_1,\dots,\iota_l) \in \N_0^l $ we adopt the 
common notation 
\begin{equation*}
\abs{I}=\sum_{j=1}^l \iota_j \qquad\text{and} \qquad I!=\prod_{j=1}^l \iota_j!
\end{equation*}
and set 
\begin{equation*}
  \l^I:=\l_1^{\iota_1} \dots \l_l^{\iota_l}.
\end{equation*}
If  $p:=H_1^{\iota_1} \otimes \dots \otimes H_l^{\iota_l} \in \polya$,
then $p(\l)=\l^I$.
The constant coefficient differential operator $\partial(p)$ 
on $\complex{A}$ (or $\liecomplex{a}$) will be denoted $\partial_a^I$, that is
\begin{equation*}
  \partial_a^I:=\partial(p)=\partial(H_1)^{\iota_1} \cdots 
                                    \partial(H_l)^{\iota_l}.
\end{equation*}
We set 
\begin{equation*}
  \partial_\l^I:=\partial(\xi_1)^{\iota_1} \cdots 
                                    \partial(\xi_l)^{\iota_l}=
\left(\frac{\partial}{\partial \l_1}\right)^{\iota_1} \dots  
             \left(\frac{\partial}{\partial \l_l}\right)^{\iota_l}
\end{equation*}
for the constant coefficient differential operator 
on $\frakacs$ corresponding to 
$\xi_1^{\iota_1} \otimes \dots \otimes \xi_l^{\iota_l}$.
Observe that, for multi-indices $I=(\iota_1,\dots,\iota_l)$ 
and $K=(\kappa_1,\dots,\kappa_l)$,
one has
\begin{equation} \label{eq:diffuno}
\partial^I_a e^\l=\l^I e^\l
\end{equation}
and
\begin{equation} \label{eq:diffdue}
\left.\big(\partial^K_\l \l^I\big) \right|_{\l=0}=K! \; \delta_{I,K} 
\end{equation}
with $\delta_{I,K}=1$ if $\iota_j=\kappa_j$ for all 
$j=1,\dots,l$ and equal to $0$ otherwise.

Since $\Delta \in \polyAc$, 
every $D \in  \C_\Delta[\complex{A}]\otimes \polya$ can be written as
\begin{equation*}
  D=\Delta^k \; {\sum}_I \; \omega_I \otimes \partial_a^I,
\end{equation*}
where $\sum_I$ denotes a sum over a finite set of multi-indices $I$,
$\omega_I \in \polyAc$ and $k \in \Z$. 
We can assume that this representation of $D$
has been chosen with $k$ as big as possible. 
Then $D \in \polyAc\otimes \polya$ if and only if $k \in \N_0$.
(In particular, the condition $k \in \N_0$ is independent of the 
choice of the basis 
$\{H_1,\dots,H_j\}$ of $\frak a$ fixed for the representation of $D$.)

\begin{Lemma} \label{lemma:relprime}
Suppose $\Sigma$ is a reduced root system and $\a, \beta \in \Sigma^+$ with $\a \neq \beta$.
Then $1-e^{-2\a}$ and $1-e^{-2\beta}$ are relatively prime in   
$\polyAc$.
\end{Lemma}
\pf
This follows from \cite{Bou}, Ch. VI, \S 3, Lemma 1(ii).
\qed

\begin{Lemma} \label{lemma:Dnonsingular}
  Let $D  \in  \C_\Delta[\complex{A}]\otimes \polya$. Suppose
 $D=\Delta^k \; {\sum}_I \; \omega_I \otimes \partial_a^I$ (with 
$\omega_I \in \polyAc$ and $k \in \Z$ maximal) 
is the representation of $D$ with 
respect to the fixed basis $\{H_1,\dots,H_l\}$ of $\frak a$.
\begin{enumerate}
\thmlist
\item
If all coefficients 
$\Delta^k \omega_I$ are nonsingular on $A$, then $D \in 
\polyAc\otimes \polya$. 

\item \label{lemma:DnonsingularonAplus}
If $D$ is $W$-invariant and all coefficients 
$\Delta^k \omega_I$ are nonsingular on $\overline{A^+}$, 
then $D \in \polyAc\otimes \polya$. 

\item \label{lemma:coeffsDzeroonbdryAplus}
If $D$ is $W$-invariant and all coefficients $\Delta^k \omega_I$
vanish on $\partial(A^+)$, then they also vanish on 
$\bigcup_{\a \in \Sigma^+} \{a \in A: 
\a(\log a)=0\}$.
\end{enumerate}
\end{Lemma}
\pf
By maximality of $k$, there is a multi-index $I$ so that 
$\Delta$ does not divide
$\omega_I$. Hence there is $\a \in \Sigma^+$ so that 
$1-e^{-2\a}$ does not divide
 $\omega_I$. If
$$\Delta^k \omega_I=e^{k \rho} (1-e^{-2\a})^k 
\prod_{\b\in \Sigma^+\setminus\{\a\}} (1-e^{-2\beta})^k \; \omega_I$$
is non-singular on $A$, then $k \in \N_0$ by Lemma \ref{lemma:relprime}.
This proves (a).

For (b), notice first that
the possible singularities of $\Delta^k \omega_I$ in $A$ lie in the zero 
set of $\Delta$, i.e. along the hypersurfaces 
$\mathcal H_\a:=\{a\in A:\a(\log a)=0\}$ with $\a \in \Sigma^+$.
By assumption, no singularities of the coefficients $\Delta^k \omega_I$ 
lie in $\overline{A^+}$, that is, no hypersurface $\mathcal H_{\beta}$
is singular when $\beta \in \Pi$ is a positive simple root. 
Suppose $\mathcal H_\a$ is 
a singular hypersurface of a coefficient $\Delta^k \omega_K$ of $D$. 
There exist $w \in W$ and $\beta \in\Pi$ such that $w\a=\beta$.  
Since $D$ is $W$-invariant, we have
\begin{equation} \label{eq:altroD}
D=wD={\sum}_I \; w(\Delta^k \omega_I) \otimes w\partial_a^I.
\end{equation}
From (\ref{eq:diffuno}) we obtain for arbitrarily fixed $\l \in \frakacs$
\begin{equation*}
  (De^\l)e^{-\l}=\Delta^k \; {\sum}_I \; \omega_I \big(\partial_a^I e^\l \big)e^{-\l}=
             {\sum}_I \;  \Delta^k\omega_I \, \l^I.
\end{equation*}
Observe that 
\begin{equation*}
   \big(w\partial_a^I\big)(e^\l):=w\big(\partial_a^I(w^{-1}e^\l)\big)=
   w\big(\partial_a^I e^{w^{-1}\l}\big)= w\big((w^{-1}\l)^I e^{w^{-1}\l}\big)=(w^{-1}\l)^I e^\l.
\end{equation*}
Hence (\ref{eq:altroD}) yields
\begin{equation*}
  (De^\l)e^{-\l}= {\sum}_I \;  w\big(\Delta^k\omega_I\big) (w\l)^I.
\end{equation*}
Therefore, for all $\l \in \frakacs$,
\begin{equation*}
  {\sum}_I \;  w\big(\Delta^k\omega_I\big) (w\l)^I = {\sum}_I \;  
\Delta^k\omega_I \,\l^I.
\end{equation*}
 Replacing $\l$ with $w\l$, we conclude
\begin{equation*}
   {\sum}_I \;  w\big(\Delta^k\omega_I\big) \l^I ={\sum}_I \;  
\Delta^k\omega_I \,(w\l)^I,
\end{equation*}
and thus, by (\ref{eq:diffdue}),
\begin{equation} \label{eq:womegaK}
  K! \; w\big(\Delta^k\omega_K\big)= 
{\sum}_I \; \Delta^k\omega_I \left.\Big(\partial_\l^K (w\l)^I\Big) 
\right|_{\l=0}.
\end{equation}
Since $\Delta^k\omega_K$ is singular along $\mathcal H_\a$, it follows that
$w\big(\Delta^k\omega_K\big)$ is singular along 
$w\mathcal H_\a=\mathcal H_\beta$.
According to (\ref{eq:womegaK}), $w\big(\Delta^k\omega_K\big)$ is a
linear combination of the coefficients  $\Delta^k\omega_I$. So there must be 
a multi-index $I$ for which $\Delta^k \omega_I$ is singular 
along $\mathcal H_\beta$, in contradiction to our assumption. 
Thus no coefficient $\Delta^k\omega_I$ 
can be singular in $A$, and the claim follows then from (a).

Part (c) follows immediately from (\ref{eq:womegaK}) and from the 
fact that $W$ acts transitively on the Weyl chambers.
\qed
  
The following key lemma is a consequence of a more general result, 
which will be proven in Appendix \ref{appendix:estPhi} (see Corollary 
\ref{cor:estonA}). 

\begin{Lemma} \label{lemma:partialPhi}
  For every $m \in \mathcal M^+$ and every multi-index $I$ the function 
$\left.\Delta(m;a) \partial_\l^I \Phi(m;\l,a)\right|_{\l=0}$ 
extends continuously on $\overline{A^+}$ 
by setting it equal to zero on the boundary $\partial(A^+)$ of $\overline{A^+}$.
\nqed
\end{Lemma}

\noindent {\em Proof of Theorem \ref{thm:Dm}.\; }
The $W$-invariance of $D_m$ is an immediate consequence of Lemma 
\ref{lemma:deltainv} 
and the fact that $D_+(m)=G_+(m;0)$ is $W$-invariant. 
Lemma  \ref{lemma:deltainv}  also
guarantees that $D_m \in  \C_\Delta[\complex{A}]\otimes \polya$. 
Hence we can write 
\begin{equation} \label{eq:Dmaspartial}
  D_m=\Delta^k \; {\sum}_I \; \omega_I \otimes \partial_a^I,
\end{equation}
where $I$ are multi-indices, 
$\omega_I \in \polyAc$ and $k \in \Z$ is maximal. 
Because of (\ref{eq:diffuno}), 
\begin{equation} \label{eq:Dmaspartialbis}
 \big( D_m e^\l \big) e^{-\l}=\Delta^k \; {\sum}_I \; \omega_I \l^I.
\end{equation}
From (\ref{eq:diffdue}), (\ref{eq:Dmaspartialbis}), (\ref{eq:DplusHC}) 
and from the product rule for differentiation we obtain
\begin{align*}
\Delta(a)^k \omega_I(a)&=
    \frac{1}{I!} \;\partial_\l^I 
     \left.\Big( (D_m e^{\l(\log a)}) e^{-\l(\log a)} \Big) 
        \right|_{\l=0}\\
       &=\frac{1}{I!} \; \partial_\l^I  
  \left. \Big( \frac{e^{-\l(\log a)}}{\pedPi{c}^+(m;-\l)} \; 
     \Delta(m;a) \Phi(m;\l,a) \Big)\right|_{\l=0}\\
   &=\sum_{J+K=I} \frac{1}{J! K!}  \;  \partial_\l^J   
   \left. \Big(  \frac{e^{-\l(\log a)}}{\pedPi{c}^+(m;-\l)} \Big) 
        \right|_{\l=0} 
   \left. \Big(\Delta(m;a) \partial_\l^K \Phi(m;\l,a) \Big)\right|_{\l=0}.
\end{align*}
Lemma \ref{lemma:partialPhi} implies then that each coefficient 
$\Delta^k\omega_I$ 
in the representation (\ref{eq:Dmaspartial}) of $D_m$ extends continuously
to $\overline{A^+}$ by setting it equal to $0$ on $\partial(A^+)$. 
Since $D_m$ is $W$-invariant, the theorem thus follows from 
Lemma \ref{lemma:Dnonsingular} (\ref{lemma:DnonsingularonAplus}) 
and (\ref{lemma:coeffsDzeroonbdryAplus}).
\qed 

In the following corollaries we shall always assume that $m\in\mathcal M^+$
is a fixed even multiplicity function.

\begin{Cor} \label{cor:shiftnonsing}
$D_-(m):=G_-(-m;m) \in \polyAc \otimes \polya$. 
\end{Cor}

\begin{Rem}
  Corollary \ref{cor:shiftnonsing} is consistent with the idea that, since 
$G_-(-m;m)$ ``removes'' the singularities of $\Phi(m;\l,a)$ 
along $\partial(A^+)$,
it should not have singular coefficients there.
\end{Rem}

The following corollary will be important to determine formulas for the $\Theta$-spherical
functions in which, as in the complex case, an alternating sum appears (see Corollary 
\ref{cor:formulas}). Its proof, which 
depends on Theorem \ref{thm:Dm} and also on some elementary algebraic
properties of $\polyAc$, can be found in Appendix \ref{app:proofcor}.

\begin{Cor} \label{cor:pieceofDm}
$\Delta(m)G_+(m-2;2)\circ \Delta(-1)\in \polyAc \otimes \polya$. Consequently
$\Delta(m)G_+(m-2;2)$ is a  $W$-invariant element of $\polyAc \otimes \polya$.
\end{Cor}
\pf 
See Appendix \ref{app:proofcor}.
\qed

\begin{Cor} \label{cor:division}
Let $U$ denote the tubular neighborhood of $A$ from Theorem \ref{thm:HC},
and let $D_m=\Delta(m)D_+(m)$ be as in Theorem \ref{thm:Dm}. 
Set
\begin{equation}
  \label{eq:pi}
  \pi(\l):=\prod_{\a\in \Sigma^+} \la.
\end{equation}
Then the polynomial $\pi(\l)$ divides 
$D_me^{\l(\log a)}$ for all $a \in U$.  
\end{Cor}
\pf
By (\ref{eq:splitshiftplus}) 
 and Example \ref{ex:shiftcomplex}, we have
\begin{align*}
  D_m e^{\l(\log a)}&=\Delta(m)D_+(m)e^{\l(\log a)}\\
&=\Delta(m)G_+(m-2;2)D_+(2)e^{\l(\log a)}\\
&=\Delta(m)G_+(m-2;2)\Big[ (-1)^{\abs{\Sigma^+}} \Delta^{-1}
\prod_{\a\in \Sigma^+} \frac{\partial(A_\a)}{\inner{\a}{\a}} e^{\l(\log a)}
\Big]\\
&=\pi(\l)\,(-1)^{\abs{\Sigma^+}}\Big(\Delta(m)G_+(m-2;2) \circ \Delta(-1)\Big)
e^{\l(\log a)}, 
\end{align*}
and $\Delta(m)G_+(m-2;2) \circ \Delta(-1) \in \polyAc \otimes \polya$ by 
Corollary \ref{cor:pieceofDm}.
\qed

\section{Explicit formulas for the $\Theta$-spherical functions} 
\label{section:formulas}

In this section we prove explicit formulas for the $\Theta$-spherical functions
corresponding to even multiplicity functions. They generalize the well-known
formula by Harish-Chandra for the spherical functions on Riemannian symmetric 
spaces of the noncompact type with a complex structure and the formula 
by Faraut, Hilgert and \'Olafsson for the spherical functions on NCC symmetric
spaces with complex structure. 

We introduce the polynomial
\begin{equation*}
  \pedtheta{e}^+(m;\l)
           =(-1)^{\frac{1}{2}\, \sum_{\a\in \rootstheta^+} \ma} 
        \prod_{\a\in \rootstheta^+} \prod_{k=-\ma/2+1}^{\ma/2-1} (\l_\a-k)
  \end{equation*}
(with empty products equal to $1$).
Recall the polynomial $\pedtheta{e}^-(m;\l)$ and $\pi(\l)$ from from (\ref{eq:eminus}) 
and (\ref{eq:pi}), respectively. 
Notice that 
\begin{equation}
  \label{eq:epluseminus}
  \pedtheta{e}^-(m;\l)\pedtheta{e}^+(m;\l)=(-1)^{d(\Theta,m)} \pedPi{e}^+(m;\l).
\end{equation}
Finally, let 
\begin{equation*}
  \pedtheta{c}^{+,c}(m;\l):=\frac{\pedPi{c}^+(m;\l)}{\pedtheta{c}^+(m;\l)}=
\prod_{\a \in \Sigma^+ \setminus \rootstheta^+} \prod_{k=0}^{\ma/2-1} \frac{1}{\l_\a+k} =
(-1)^{d(\Theta,m)} \pedtheta{c}^-(m;\l).
\end{equation*}

\begin{Thm} \label{thm:formulas}
\begin{enumerate}
\thmlist
\item
Let $m \in \mathcal M^+$ be even, and let $D_m=\Delta(m)D_+(m) 
\in \polyAc\times \polya$ be as in
(\ref{eq:Dm}). 
Then the $\Theta$-spherical function $\pedtheta{\varphi}(m;\l,a)$ 
is determined by the formula
\begin{equation}
  \label{eq:formula}
 \pedtheta{e}^-(m;\l) \Delta(m;a) \pedtheta{\varphi}(m;\l,a)=
 \frac{1}{\pi(\l)\pedtheta{e}^+(m;\l)} D_m \Big(\sum_{w \in \pedtheta{W}}e^{w\l(\log a)}\Big).
\end{equation}
The left-hand side of (\ref{eq:formula}) is a holomorphic function of $(\l,a)\in \frakacs \times
\pedtheta{A}U$. For fixed $\l \in \frakacs$ the function 
$$D_m \Big(\sum_{w \in \pedtheta{W}}e^{w\l(\log a)}\Big)$$
is single valued and 
holomorphic in $A\pedtheta{U}^\l$, where $\pedtheta{U}^\l$ ($\supset U$)
is a largest open neighborhood of $e$ in $T$ in  which $e^{w\l(\log a)}$ 
is single-valued for all $w \in 
\pedtheta{W}$. As a function of $(\l,a)$, it is single valued and holomorphic in $\frakacs 
\times \pedtheta{U}$ with $\pedtheta{U}:=\cap_{\l \in \frakacs} \pedtheta{U}^\l$.
We can take $\pedtheta{U}^\l=T$ when $\l \in P$. 
\item
For all $(\l,a)\in \frakacs \times \pedtheta{A}U$, we have the following equality 
(as meromorphic functions
of $\l$ with singularities located along the zero set of $\pedtheta{e}^-(m;\l)$):
\begin{equation}
  \label{eq:formulabis}
  \Delta(m;a) \pedtheta{\varphi}(m;\l,a)=
(-1)^{d(\Theta,m)} 
\Big[{\prod_{\a\in\Sigma^+} \prod_{k=0}^{\ma/2-1} (k^2-\la^2)}\Big]^{-1} \;
D_m \Big(\sum_{w \in \pedtheta{W}}e^{w\l(\log a)}\Big)
\end{equation}
with $d(\Theta,m)$ as in (\ref{eq:dmtheta}).
\item
If $\Theta=\Pi$, then $\Delta(m;a) \pedPi{\varphi}(m;\l,a)$ extends as a 
holomorphic function on $\frakacs \times A\pedPi{U}$ by means of the formula
\begin{equation}\label{eq:formulaPi} 
 \Delta(m;a) \pedPi{\varphi}(m;\l,a)= 
\Big[{\prod_{\a\in\Sigma^+} \prod_{k=0}^{\ma/2-1} (k^2-\la^2)}\Big]^{-1} \;
D_m \Big(\sum_{w \in W}e^{w\l(\log a)}\Big).
\end{equation}
If, moreover, $\l \in P$, then $\Delta(m;a) \pedPi{\varphi}(m;\l,a)$ extends
by (\ref{eq:formulaPi}) as $W$-invariant entire function on $\complex{A}$.
\end{enumerate}
\end{Thm}
\pf
To prove (a), observe first that,
since $\pedtheta{c}^{+,c}(m;\l)$ and 
$\pedtheta{c}^{+}(m;\l)\pedtheta{c}^{+}(m;-\l)$ are 
$\pedtheta{W}$-invariant, we have for all $w \in \pedtheta{W}$
\begin{align}
 \pedtheta{c}^-(m;\l)\pedtheta{c}^+(m;w\l)\pedPi{c}^+(m;-w\l)&=
 \pedtheta{c}^-(m;\l)\pedtheta{c}^+(m;w\l)\pedtheta{c}^+(m;-w\l)\pedtheta{c}^{+,c}(m;-w\l)
\notag \\
&=\pedtheta{c}^-(m;\l)\pedtheta{c}^+(m;\l)\pedtheta{c}^+(m;-\l)\pedtheta{c}^{+,c}(m;-\l)
\notag \\
&=\Big[\Big(\prod_{\a \in \Sigma^+ \setminus \rootstheta^+} \la\Big) \,
\pedtheta{e}^-(m;\l)
\Big(\prod_{\a \in \rootstheta^+} \la\Big) \,  \pedtheta{e}^+(m;\l)\Big]^{-1}
\notag \\
&=\Big[\pi(\l) \pedtheta{e}^-(m;\l)\pedtheta{e}^+(m;\l)\Big]^{-1}  \label{eq:cande}
\end{align}
Recall from Theorem \ref{thm:sphentire} that the function 
 $\pedtheta{e}^-(m;\l)\pedtheta{\varphi}(m;\l,a)$ is entire as a function of
$(\l,a) \in \frakacs \times \pedtheta{U}$, where $\pedtheta{U}$ denotes a 
$\pedtheta{W}$-invariant tubular neighborhood of $\pedtheta{A}$ in 
$\complex{A}$.

The definition of $\Theta$-spherical functions, (\ref{eq:DplusHC}) 
and (\ref{eq:cande}) 
give for all $a \in A^+U$ and $\l \in \frakacs \setminus P$
\begin{align*}
  \pedtheta{\varphi}(m;\l,a)&=\pedtheta{c}^-(m;\l) \; \sum_{w \in \pedtheta{W}}
                                                         \pedtheta{c}^+(m;w\l) \Phi(m;w\l,a)\\
&=\pedtheta{c}^-(m;\l) \; \sum_{w \in \pedtheta{W}}
            \pedtheta{c}^+(m;w\l) \pedPi{c}^+(m;-w\l) D_+(m)e^{w\l(\log a)}\\
&=\frac{1}{\pi(\l) \pedtheta{e}^-(m;\l)\pedtheta{e}^+(m;\l)} \sum_{w \in \pedtheta{W}}
D_+(m)e^{w\l(\log a)}.
\end{align*} 
Hence
\begin{equation*}
 \pedtheta{e}^-(m;\l) \Delta(m;a) \pedtheta{\varphi}(m;\l,a)=
 \frac{1}{\pi(\l)\pedtheta{e}^+(m;\l)} D_m \Big(\sum_{w \in \pedtheta{W}}e^{w\l(\log a)}\Big).
\end{equation*}
Since $D_m \in \polyAc \otimes \polya$ is $W$-invariant, 
the formula extends by $\pedtheta{W}$-invariance
to $\frakacs \times \pedtheta{A}U$.

Parts (b) and (c) follow from (a) with easy computations. 
\qed  

\begin{Cor} \label{cor:formulas}
Let $m \in \mathcal M^+$ be even. Then there exist differential operators
$G_m,\widetilde{G}_m \in \polyAc\otimes \polya$ with $\widetilde{G}_m=
G_m \circ \Delta$ and $\widetilde{G}_m$ $W$-invariant such that
for all $(\l,a) \in \frakacs \times \pedtheta{A}U$ the following equality of 
meromorphic functions of $\l$ holds: 
\begin{align*}
 \Delta(m;a)\pedtheta{\varphi}(m;\l,a)&=
\frac{(-1)^{d(\Theta,m)}}{\pi(\l)} \Big[\prod_{\a\in\Sigma^+} 
   \prod_{k=1}^{m_\a/2-1} (k^2-\l_\a^2)\Big]^{-1} \;
  G_m \Big( \sum_{w\in \pedtheta{W}} \e(w) e^{w\l(\log a)} \Big)\\
&= (-1)^{d(\Theta,m-2)} \Big[\prod_{\a\in\Sigma^+} 
   \prod_{k=1}^{m_\a/2-1} (k^2-\l_\a^2)\Big]^{-1} \;
  \widetilde{G}_m \, \pedtheta{\varphi}(2;\l,a).
\end{align*}
If $m=2$, then $G_m=\id$. 
If $\Theta=\Pi$, then  
\begin{align*}
 \Delta(m;a)\pedPi{\varphi}(m;\l,a)&=
\Big[\pi(\l) \prod_{\a\in\Sigma^+} 
   \prod_{k=1}^{m_\a/2-1} (k^2-\l_\a^2)\Big]^{-1} \;
  G_m \Big( \sum_{w\in W} \e(w) e^{w\l(\log a)} \Big)\\
&= \Big[\prod_{\a\in\Sigma^+} 
   \prod_{k=1}^{m_\a/2-1} (k^2-\l_\a^2)\Big]^{-1} \;
  \widetilde{G}_m \, \pedPi{\varphi}(2;\l,a),
\end{align*}
in which the right hand sides extend as holomorphic functions of 
$(a,\l) \in AU \times \frakacs$.
If $\l\in P$ is fixed, then the right-hand sides of the previous equalities
extend as holomorphic functions of $a\in\complex{A}$. 
\end{Cor}
\pf
Set $G_m:=\Delta(m)G_+(m-2;2)\circ \Delta(-1)$. Then the formulas 
follow immediately from Theorem \ref{thm:formulas}, Corollary \ref{cor:pieceofDm}
and the equalities in the proof of Corollary \ref{cor:division}
(see also Example \ref{ex:phithetacomplex} below).
\qed

\begin{Ex}[The rank-one case] 
\label{ex:formularone} 
In the rank-one case the only two possiblities are
$\Theta=\emptyset$ and $\Theta=\Pi$. If $\Theta=\emptyset$, we have a multiple (with a function
of $\l$ as a constant) of the Harish-Chandra series
\begin{equation}
  \label{eq:HCseriesrankone}
  \Phi(m;\l,z)=2^{-m/2} \prod_{k=0}^{\frac{m}{2}-1} (\l-k) 
     \Big( \frac{1}{\sinh z} \; \frac{d}{dz}\Big)^{m/2} e^{\l z}.
\end{equation}
In particular, 
\begin{equation*}
  \Phi(2;\l,z)=\frac{1}{\Delta(z)}\; e^{\l z}.
\end{equation*}
If $\Theta=\Pi$, then 
\begin{equation}
  \label{eq:sphrankone}
  \frac{\pedPi{\varphi}(m;\l,z)}{\pedPi{c}^+(m;\rho)}=2^{m/2-1}
\prod_{k=0}^{m/2-1} \frac{(m/2+k)}{(\l^2-k^2)}  \;
 \Big( \frac{1}{\sinh z} \; \frac{d}{dz}\Big)^{m/2} \cosh(\l z)
\end{equation}
is the formula for Harish-Chandra's spherical functions on the Riemannian symmetric 
spaces $\SO_0(1,n)/\SO(n)$ (with $n=m+1$), as found by Takahashi 
(see \cite{Taka}, Formula (21), p. 326).
In particular
\begin{equation*}
   \frac{\pedPi{\varphi}(2;\l,z)}{\pedPi{c}^+(2;\rho)}=\frac{1}{\l}\; \frac{\sinh(\l z)}{\sinh z}
\end{equation*}
in the rank-one complex case.
\end{Ex}

\begin{Ex}[The complex case] 
\label{ex:phithetacomplex}  
In this case $d(\Theta,m)=\abs{\Sigma^+\setminus \rootstheta^+}$.
Since
\begin{equation*}
\prod_{\a \in \Sigma^+} \partial_\a\; e^{w\l(\log a)}= 
\Big(\prod_{\a \in \Sigma^+} \inner{w\l}{\a}\Big)  e^{w\l(\log a)}=
\e(w) \Big(\prod_{\a \in \Sigma^+}  \inner{\l}{\a}\Big) e^{w\l(\log a)}, 
\end{equation*}
we obtain 
\begin{equation}
  \label{eq:phithetacomplex}
  \pedtheta{\varphi}(2;\l,a)=(-1)^{\abs{\Sigma^+\setminus \rootstheta^+}}
                \frac{1}{\pi(\l)}\; \frac{1}{\Delta(a)} \sum_{w\in \pedtheta{W}} \e(w) e^{w\l(\log a)}.
\end{equation}
In particular,
\begin{equation}\label{eq:sphcomplex}
  \frac{\pedPi{\varphi}(2;\l,a)}{\pedPi{c}^+(2;\rho)}=
                \frac{\pi(\rho)}{\pi(\l)}\; 
  \frac{1}{\Delta(a)} \,\sum_{w\in W} \e(w) e^{w\l(\log a)}
\end{equation}
gives Harish-Chandra's formula for the spherical functions on Riemannian
symmetric spaces with a complex structure, and (with $W_0:=\ped{\Pi_0}{W}$)  
\begin{equation}\label{eq:sphcomplexNCC}
  \frac{\ped{\Pi_0}{\varphi}(2;\l,a)}{\ped{\Pi_0}{c}^+(2;\rho) \ped{\Pi_0}{c}^-(2;\rho) }=
                \frac{\pi(\rho)}{\pi(\l)}\; 
  \frac{1}{\Delta(a)} \sum_{w\in W_0} \e(w) e^{w\l(\log a)}
\end{equation} 
gives the formula for the spherical functions on NCC symmetric spaces with complex structure
as in \cite{FHO}.  
\end{Ex}

\section{Some remarks on the Paley-Wiener space}

Before proceeding with the the proof of Theorem  \ref{thm:pw}, 
let us add some remarks on the definition of the Paley-Wiener space. 
Let $m \in \mathcal M^+$ be an even multiplicity function, 
$\Theta \subset \Pi$ a set of positive simple roots, and $C$ a compact, convex and 
$\pedtheta{W}$-invariant subset of $\pedtheta{\frak a}$ 
 Recall from Definition \ref{def:pwspacetheta} that the Paley-Wiener space
$\pwthetamC$ is the space of all 
$\pedtheta{W}$-invariant meromorphic functions $g:\frakacs \to \C$ 
satisfying the following properties:
\begin{enumerate}
\renewcommand{\labelenumi}{\theenumi.}
\item 
The function
$\pedtheta{e}^-(m;\l) g(\l)$ is entire of exponential type $C$ 
and rapidly decreasing, that is 
for every $N \in \N$  
there is a constant $C_N\geq 0$ such that
$$
\abs{\pedtheta{e}^-(m;\l) g(\l)} \leq C_N (1+\abs{\l})^{-N} e^{q_C(\Re\l)}
$$
for all $\l \in \frakacs$.
\item 
The function 
\begin{equation*}
\Pavtheta g(\l):=\sum_{w \in \pedtheta{W}\setminus W} g(w\l)
\end{equation*}
 extends to an entire function on $\frakacs$. 
\end{enumerate} 

\begin{Prop} \label{prop:PWspace}
Condition 2 above is 
automatically satified in the following cases:
\begin{enumerate}
\thmlist
\item
Euclidean case ($m=0$, $\Theta$ arbitrary);
\item
Heckman-Odpam case ($\Theta=\Pi$, $m$ arbitrary);
\item
Complex case ($m=2$, $\Theta$ arbitrary).
\end{enumerate}
\end{Prop}
\pf
In the first two cases we have $\pedtheta{e}^-(m;\l)=1$, hence all functions in
$\pedtheta{\rm PW}(m;C)$ are entire.
In the complex case $\pedtheta{e}^-(2;\l)=\prod_{\a \in \Sigma^+\setminus \rootstheta^+} \l_\a$.
Condition 1 in the definition of the Paley-Wiener space implies
that every $g \in \pwthetamC$ has at most a first order pole along each hyperplane 
$\l_\a=0$. The same is then true for $\Pavtheta g$. But $\Pavtheta g$ is $W$-invariant, so 
it cannot have first order poles on hyperplanes $\l_\a=0$ with $\a\in \Sigma^+$. 
They must be therefore removable singularities, 
i.e. Condition \ref{item:pwspacethetadue} holds.
\qed

\begin{Rem}
Condition \ref{item:pwspacethetadue} in Definition \ref{def:pwspacetheta} does not 
follow in general from Condition \ref{item:pwspacethetauno}. For instance, in the rank-one
case with $m=4$ and the usual identification of $\frak a$ with $\R$,
we have $\pedtheta{e}^-(4;\l)=\l(\l-1)(\l+1)$. If $h \in \PW(C)$, where $C$ is an arbitrary 
compact interval in $(0,+\infty)$, then the meromorphic function 
$g(\l)=\big(\frac{a}{\l-1}+\frac{b}{\l+1}\big) h(\l)$ satisfies the first condition, 
but not the second  when $a\neq -b$. 
\end{Rem}

Let $\ell$ denote the length function on $W$ with respect to the simple reflections
$r_\beta$ with $\beta \in \Pi$. For $\Theta \subset \Pi$ set
\begin{equation}\label{eq:apthetaW}
\aptheta{W}:=\{w \in W: \text{$\ell(wr_\a)>\ell(w)$ for all $\a \in \Theta$}\}.
\end{equation}
Then every element $w \in W$ can be uniquely written as 
$w=uv$ with $u \in \aptheta{W}$ and $v \in \pedtheta{W}$ (see e.g. \cite{Hu}, p. 19). 
The element $u$ is characterized as the unique element of smallest length in the coset 
$w\pedtheta{W}$. Hence
a set of coset representatives for $\pedtheta{W}\setminus W$ consists of 
the elements of $\{u^{-1}: u \in \aptheta{W}\}$. In particular we obtain
\begin{equation}\label{eq:Pavthetarep}
\Pavtheta g(\l)=\sum_{u \in \aptheta{W}} g(u^{-1}\l).
\end{equation}

\begin{Lemma} \label{lemma:apthetaWonatheta}
Let be $\aptheta{W}$ as in (\ref{eq:apthetaW}).  
Then for $u,u' \in \aptheta{W}$ with $u\neq u'$ one has 
$u \pedtheta{\frak a} \cap u' \pedtheta{\frak a} =\emptyset$. 
\end{Lemma}
\pf
Since $W$ acts simply transitively on the Weyl chambers, there are $\abs{W}$ Weyl chambers. 
A nonempty intersection of $u \pedtheta{\frak a}$ and $u' \pedtheta{\frak a}$ would contain 
at least a Weyl chamber. This is not possible since $\abs{W}=\abs{\aptheta{W}} 
\abs{\pedtheta{W}}$.
\qed

The following lemma, which is an easy modification of Lemma 5.13 in \cite{He3},
will be applied several times in the sequel.

\begin{Lemma} \label{lemma:malgrange}
Let $C$ be a compact convex subset of $\frak a\equiv\R^l$, and 
let $g:\frakacs\equiv \C^l \to \C$ belong to the Paley-Wiener space $\PW(C)$.  
Suppose $q$ is a polynomial so that $h:=g/q$ is entire. Then $h \in \PW(C)$. 
\end{Lemma}

In the following we denote by $\conv (W( C))$ the closed convex hull of the $W$-orbit
of a subset $C$ of $\frak a$. Moreover, we normalize the  Lebesgue measure on $i\frak a^*$
so that the inverse transform to $\mathcal F_A$ is given by
\begin{equation*}
 \mathcal F_A^{-1}g(a):=\int_{i\frak a^*} g(\l) e^{-\l(\log a)} \; d\l.
\end{equation*}

\begin{Prop} \label{prop:Pavthetainjective}
Let $C \subset \pedtheta{\frak a}$ be compact, convex and $\pedtheta{W}$-invariant.
Suppose $g \in \pwthetamC$.
Then $\Pavtheta g \in \PW\big(\conv (W (C))\big)^W$.
Moreover, the map  $\Pavtheta: \pwthetamC \to \PW\big(\conv (W( C))\big)^W$
is linear and injective.
\end{Prop}
\pf
The first statement is an immediate consequence of Lemma \ref{lemma:malgrange} 
with $f:=\Pavtheta g$ and 
\begin{equation} \label{eq:q}
 q(\l):= \prod_{\a \in \Sigma} \; \prod_{k=-\ma/2+1}^{\ma/2-1} (\l_{\a}-k).
\end{equation}
Indeed $q$ is $W$-invariant and  
$q(\l)\Pavtheta g(\l)=\sum_{w \in \pedtheta{W}\setminus W} q(\l)g(w\l)$,
with $q(\l)g(w\l) \in \PW(wC)$ for all $w \in W$.

To prove the injectivity of $\Pavtheta$, suppose $g \in \pwthetamC$ satisfies 
$\Pavtheta g=0$. With $q$ as above, one has
$q g \in \PW(C)^{\pedtheta{W}}$, and $\Pavtheta(qg)=q\Pavtheta g=0$.
Hence for all $a \in A$
\begin{align}
  0&=\mathcal F_A^{-1}\big(\Pavtheta(qg)\big)(a) \notag \\
    &=\sum_{w\in \pedtheta{W}\setminus W} 
                        \int_{i\frak a^*} q(w\l)g(w\l) e^{-\l(\log a)} \; d\l  \notag\\
    &=\sum_{w\in \pedtheta{W}\setminus W} 
                        \int_{i\frak a^*} q(\l)g(\l) e^{-\l(\log wa)} \; d\l  \notag\\
    &=\sum_{w\in \pedtheta{W}\setminus W} \mathcal F_A^{-1}(qg)(wa) \notag\\
    &=\sum_{u \in \aptheta{W}} \mathcal F_A^{-1}(qg)(u^{-1}a). 
              \label{eq:sumzero}
\end{align}
By the classical Paley-Wiener theorem 
\begin{equation*}
  \supp \big[\mathcal F_A^{-1}(qg) \circ u^{-1}\big] 
\subset u\exp C \subset u\pedtheta{A}.
\end{equation*}
Lemma \ref{lemma:apthetaWonatheta} and (\ref{eq:sumzero}) therefore imply 
$\mathcal F_A^{-1}(qg) \circ u^{-1}=0$
for all $u \in \aptheta{W}$. In particular,  $\mathcal F_A^{-1}(qg)=0$. Thus $qg=0$,
which proves $g=0$ because $g$ is meromorphic.
\qed

\begin{Cor}
Let $q$ be as in (\ref{eq:q}), and let 
$\pedtheta{\iota}:\pedtheta{A} \hookrightarrow A$ be the inclusion map. Then
\begin{equation*}
  \Big(\frac{1}{q} \; \mathcal F_A \circ \pedtheta{\iota} \circ 
\mathcal F_A^{-1} \circ q\big) 
\circ \Pavtheta ={\id}_{\pwthetamC}.
\end{equation*}
\end{Cor}
\pf
As in the proof of Proposition \ref{prop:Pavthetainjective}, we have for $g \in 
\pwthetamC$
\begin{align*}
\pedtheta{\iota}\; \mathcal F_A^{-1} (q \Pavtheta g)&=
\pedtheta{\iota}\; \mathcal F_A^{-1} \Pavtheta(q g)\\
   &=\pedtheta{\iota}\Big( 
\sum_{u \in \aptheta{W}} \big(\mathcal F_A^{-1}(q g) \circ u^{-1}\big)\Big)\\
&=\mathcal F_A^{-1}(q g).
\end{align*}
\qed

\section{The $\Theta$-spherical transform}

Recall from (\ref{eq:thetasphtr}) that the $\Theta$-spherical transform 
of a sufficiently regular $\pedtheta{W}$-invariant function 
$f:\pedtheta{A}\to\C$ is the $\pedtheta{W}$-invariant 
function $\pedtheta{\mathcal F}f(m)$ on $\frakacs$ defined by
\begin{equation*} 
\pedtheta{\mathcal F}f(m;\l):=
\frac{1}{\abs{\pedtheta{W}}} \, \int_{\pedtheta{A}} 
f(a) \, \pedtheta{\varphi}(m;\l,a) \Delta(m;a)\; da
\end{equation*}
with $\Delta$ as in (\ref{eq:delta}).

 \begin{Lemma} \label{lemma:sphtrentire}
If $f \in \ccitheta$, then  $\pedtheta{e}^-(m;\l)\pedtheta{\mathcal F}f(m;\l)$ is a 
$\pedtheta{W}$-invariant  entire function on $\frakacs$. 
 \end{Lemma}
\pf 
This follows immediately from Theorem \ref{thm:sphentire}. \qed

\begin{Cor} \label{cor:functionaleqtr}
Every function $f \in \ccitheta$ can be uniquely extended to a $W$-invariant function 
$\pedPi{f} \in C_c^\infty(A)^W$. Moreover
\begin{equation*}
  \big(\pedPi{\mathcal F}\pedPi{f}\big)(m;\l)=
(-1)^{d(\Theta,m)} \sum_{\pedtheta{W}\setminus W} 
  \big(\pedtheta{\mathcal F}f)(m;w\l)=
(-1)^{d(\Theta,m)}  \big(\Pavtheta  \pedtheta{\mathcal F}f)(m;\l),
\end{equation*}
where $d(\Theta,m)$ is as in (\ref{eq:dmtheta}).
\end{Cor}
\pf
Immediate consequence of Lemmas \ref{lemma:apthetaWonatheta} and 
 \ref{lemma:functionaleq}. \qed

The inversion formula for the NCC spaces was first proven in
\cite{O}. The general case was treated in
in \cite{P2} (see also \cite{PHab}). 
Observing that in the even multiplicity case 
$$\pi(\l)\ped\Pi{e}^+(m;\l)=
\frac{1}{\abs{\pedPi{c}^+(m;\l)}^2},$$ 
we can state it as follows.

 \begin{Thm}[\cite{P2}, Theorem 4.5; see also \cite{PHab}, Theorem 2.4.5] 
\label{thm:inversion}
Let $m \in \mathcal M^+$ be an even multiplicity function, and 
let $f \in \cctheta$. Then there is a constant $k>0$ (depending only on 
the normalization of the measures) so that for all $f \in \cctheta$ 
the following inversion formula holds:
For all $a \in \pedtheta{A}$
\begin{align*}
    f(a)&= (-1)^{d(\Theta,m)}  k\; \frac{\abs{W}}{\abs{\pedtheta{W}}} \; \int_{i\frak a^*} 
\big(\pedtheta{\mathcal F}f\big)(m;\l) 
\pedPi{\varphi}(m;-\l,a) \; \frac{d\l}{\abs{\pedPi{c}^+(m;\l)}^2}\\
 &=(-1)^{d(\Theta,m)}  k\; \frac{\abs{W}}{\abs{\pedtheta{W}}} \; \int_{i\frak a^*} 
\big(\pedtheta{\mathcal F}f\big)(m;\l) 
\pi(\l) \pedPi{e}^+(m;\l)\pedPi{\varphi}(m;-\l,a) \; d\l
\end{align*}
 \end{Thm}

\section{Transforms of compactly supported smooth functions}

In this section we begin the proof of the Paley-Wiener theorem  
(Theorem \ref{thm:pw}) by showing that the $\Theta$-spherical transform 
maps $C^\infty_c(C)^{\pedtheta{W}}$ into the 
Paley-Wiener space $\pwthetamC$. The key property is the fact that  
$D_m=\Delta(m)D_+(m)=D_-(m)^*$ is a $W$-invariant differential operator on $A$ with
smooth coefficients.

\begin{Thm} \label{lemma:pwinto}
Let $C \subset \pedtheta{A}$ be compact, convex and $\pedtheta{W}$-invariant.
Then the $\Theta$-spherical transform maps 
$C^\infty_c(C)^{\pedtheta{W}}$ injectively into $\pwthetamC$.
\end{Thm}
\pf
Suppose $f \in C^\infty_c(C)^{\pedtheta{W}}$. Because of
Theorem \ref{thm:Dm} we have 
\begin{align*}
\pi(\l)\pedtheta{e}^+(m;\l) \pedtheta{e}^-(m;\l) \pedtheta{\mathcal F}f(m;\l)&=
\frac{1}{\abs{\pedtheta{W}}} \, \int_{\pedtheta{A}} 
f(a) \, \pi(\l)\pedPi{e}^+(m;\l) 
\pedtheta{\varphi}(m;\l,a) \Delta(m;a)\; da\\
&=\frac{1}{\abs{\pedtheta{W}}} \, 
\int_{\pedtheta{A}} \,
f(a) \, \sum_{w \in \pedtheta{W}} D_m e^{w\l(\log a)}\; da\\
&=\frac{1}{\abs{\pedtheta{W}}} \,
\sum_{w \in \pedtheta{W}}
\int_{\pedtheta{A}} f(a) \,D_-(m)^*e^{w\l(\log a)}\; da\\
 &=\frac{1}{\abs{\pedtheta{W}}} \, 
\sum_{w \in \pedtheta{W}}
\int_{\pedtheta{A}} [D_-(m)f(a)] \,e^{\l(\log w^{-1}a)}\; da\\
&=
\int_{\pedtheta{A}} [D_-(m)f(a)] \,e^{\l(\log a)}\; da\\
&=
\big[\mathcal F_A (D_-(m)f)\big](\l).
\end{align*}
In the above equalities we have used the fact that $D_-(m)$ is a $W$-invariant 
differential operator with smooth coefficients on $A$,
which implies that $D_-(m)f(a)$ is a $\pedtheta{W}$-invariant smooth 
function on $A$ with compact support in $\exp C \subset \pedtheta{A}$.
The classical Paley-Wiener theorem for the Fourier transform yields then
that 
$$
\pedPi{\pi}(m;\l)\pedtheta{e}^+(m;\l) \pedtheta{e}^-(m;\l)\pedtheta{\mathcal F}f(m;\l)
$$
belongs to $\PW(C)$. By Lemma \ref{lemma:sphtrentire} the function 
$\pedtheta{e}^-(m;\l)\pedtheta{\mathcal F}f(m;\l)$ is entire and $\pedtheta{W}$-invariant. 
It therefore an element of 
$\PW(C)^{\pedtheta{W}}$ by Lemma \ref{lemma:malgrange}. 
Furthermore, Corollary  \ref{cor:functionaleqtr} gives
$\Pavtheta  \pedtheta{\mathcal F}f=(-1)^{d(\Theta,m)} \pedPi{\mathcal F}\pedPi{f}$,
which is entire again by Lemma  \ref{lemma:sphtrentire}. 
This proves that $\pedtheta{\mathcal F}f(m) \in \pwthetamC$. 

The injectivity of $\pedtheta{\mathcal F}$ on $C^\infty_c(C)^{\pedtheta{W}}$ 
follows from the inversion formula in Theorem \ref{thm:inversion}.
\qed

\section{Wave packets}

\begin{Def}
Let $m\in \mathcal M^+$ be a fixed even multiplicity function. The  
{\em wave-packet} of $g:\frakacs \to \C$ is the function 
$\mathcal I g=\mathcal I g(m):A \to \C$ defined by
\begin{align}\label{eq:wave}
\big(\mathcal I g\big)(a)
&:=\int_{i\frak a^*} g(\l) \pedPi{\varphi}(m;-\l,a) \; \frac{d\l}{\abs{\pedPi{c}^+(m;\l)}^2} \notag \\
& =\int_{i\frak a^*} g(\l) \pi(\l)\pedPi{e}^+(m;\l) \pedPi{\varphi}(m;-\l,a) \; d\l, 
\end{align}
provided the integrals converge. In this case, the {\em $\Theta$-wave-packet}
of $g$ is the function on $\pedtheta{A}$ obtained by restriction
of $\mathcal I g$ to $\pedtheta{A}$, that is 
\begin{equation}
  \label{eq:thetawave}
\pedtheta{\mathcal I} g = \mathcal I g \circ \pedtheta{\iota},   
\end{equation}
where $\pedtheta{\iota}:\pedtheta{A} \hookrightarrow A$ is the inclusion map.
\end{Def}

\begin{Rem} \label{rem:invarianceI}
Suppose the integrals (\ref{eq:wave}) converge for $g:A\to\C$. The $W$-invariance of 
$\pedPi{\varphi}(m;-\l,a)$ in $a \in A$ implies that 
$\mathcal I g$ is $W$-invariant. Furthermore, 
by $W$-invariance in the $\l$-variable 
of  $\pedPi{\varphi}(m;-\l,a)$ and $\abs{\pedPi{c}^+(m;\l)}^2$, we have
$\mathcal I g=\frac{\abs{\pedtheta{W}}}{\abs{W}}\; \mathcal I \Pavtheta g$. 
\end{Rem}

\begin{Lemma} \label{lemma:Igcinfty}
Let $m \in \mathcal M^+$ be a fixed even multiplicity function.
Let $C$ be a compact convex subset of $\frak a$.
Suppose that $g:\frakacs \to \C$ satisfies  $\pi(\l) \pedPi{e}^+(m;\l) g(\l) \in \PW(C)$.
Then $\mathcal Ig$ is a well-defined $W$-invariant smooth function on $A$.

In particular, if $C$ is $\pedtheta{W}$-invariant and $g \in \pwthetamC$, then 
$\pedtheta{\mathcal I}g$ is a well-defined $\pedtheta{W}$-invariant smooth function on 
$\pedtheta{A}$.
\end{Lemma}
\pf 
The assumption guarantees that for all $N \in \N$ there is a constant $C_N>0$ for which
$\abs{\pi(\l) \pedPi{e}^+(m;\l) g(\l)} \leq C_N (1+\abs{\l})^{-N}$ for all $\l \in i\frak a^*$.
Because of (\ref{eq:estsphPider}), we can therefore differentate (\ref{eq:wave}) 
under integral sign.
\qed

\begin{Lemma} \label{lemma:suppIWinv}
Let $C\subset \pedtheta{\frak a}$ be compact, convex and $\pedtheta{W}$-invariant,
and let $g \in \pwthetamC$. Then $\supp \mathcal I g \subset \conv (W (C))$.
\end{Lemma}
\pf 
By Remark \ref{rem:invarianceI} we have 
$\mathcal Ig=\frac{\abs{\pedtheta{W}}}{\abs{W}} \mathcal I  \Pavtheta g$ with
$\Pavtheta g \in \PW(\conv (W( C)))^W$. 
Formulas (\ref{eq:formula}) and (\ref{eq:wave}) give for all $h \in \PW(\conv (W( C)))^W$ 
and $a \in A$ 
\begin{align*}
  \Delta(m;a) (\pedPi{\mathcal I}h)(a)&
=\int_{i\frak a^*} h(\l)\pi(\l)\pedPi{e}^+(m;\l)\pedPi{\varphi}(m;-\l,a)\; d\l \\
&=\int_{i\frak a^*} h(\l)\sum_{w\in W} D_m e^{w\l(\log a)}\; d\l \\
&=\sum_{w\in W} D_m \int_{i\frak a^*} h(\l) e^{w\l(\log a)} \; d\l\\
&=\abs{W} \big(D_m \mathcal F_A^{-1} h\big)(a).
\end{align*}
The classical Paley-Wiener theorem implies that $\mathcal F_A^{-1} h$, and hence
$\Delta(m) \pedPi{\mathcal I}h$, has support in $\conv (W (C))$. Since 
 $\pedPi{\mathcal I}h$ is smooth, this implies the claim.
\qed

\begin{Rem}
An alternative proof of Lemma  \ref {lemma:suppIWinv}, which indeed holds for the
$\Theta$-spherical transform corresponding to an arbitrary multiplicity function $m \in 
\mathcal M^+$, can be obtained via the equality 
$\mathcal Ig=\frac{\abs{\pedtheta{W}}}{\abs{W}} \mathcal I  \Pavtheta g$ and the 
Paley-Wiener theorem for the Opdam transform (see Theorem \ref{thm:pwopdam} below).
The proof given above requires however only elementary tools.
\end{Rem}

\section{The case $\Theta=\Pi$}
\label{section:PWopdam}

In the case $\Theta=\Pi$ our spherical transform coincides (up to a constant multiple 
depending on the multiplicites and on the normalization of the measures)
with the Opdam transform.  For the latter transform, 
the following version of the Paley-Wiener theorem 
was proven for arbitrary $m \in \mathcal M^+$ in \cite{OpdActa}. 
\footnote{For arbitrary multiplicity functions 
$m \in \mathcal M^+$, the $\Theta$-spherical transform
is defined by (\ref{eq:thetasphtr}) with 
$\Delta(m):=\prod_{\a\in \Sigma^+} \abs{e^\a-e^{-\a}}^{\ma}$.}

\begin{Thm}[\cite{OpdActa}, Theorems 8.6 and 9.13(4); see also \cite{Opd00}, p. 49]
\label{thm:pwopdam}
Let  $m \in \mathcal M^+$ be arbitrarily fixed.
Let  $H \in \frak a$ and $C(H):=\conv (W( H))$. Then $\pedPi{\mathcal F}$ maps 
$C^\infty_c \big(C(H)\big)^W$ bijectively onto 
$\pedPi{\PW}\big(m;C(H)\big)\equiv \PW\big(C(H)\big)^W$. 
The inverse of $\pedPi{\mathcal F}$  is $k \pedPi{\mathcal I}$, where 
$k$ is the normalizing constant appearing in Theorem \ref{thm:inversion}.
\end{Thm}

\begin{Cor}  \label{cor:pwopdam}
Let  $m \in \mathcal M^+$ be arbitrary, and let $C$ be a compact, convex and $W$-invariant 
subset of $\frak a$. Then $\pedPi{\mathcal F}$ maps 
$C^\infty_c (C)^W$ bijectively onto 
$\pedPi{\PW}(m;C)\equiv \PW(C)^W$.  The inverse transform is 
$k\pedPi{\mathcal I}$, where 
$k$ is the normalizing constant appearing in Theorem \ref{thm:inversion}.
\end{Cor}
\pf
Since $C \subset C(H)$ for suitable $H \in \frak a$ and  $\pedPi{\mathcal F}$ is bijective with
inverse $k\pedPi{\mathcal I}$ on $C^\infty_c \big(C(H)\big)^W$, we only need to show
that $\pedPi{\mathcal F}$ maps 
$C^\infty_c (C)^W$ onto $\PW(C)^W$. One can argue as
in Theorem 8.6 (1) of \cite{OpdActa}, that is using the
estimates (\ref{eq:estsphPider}) together with the fact that $\pedPi{\varphi}(m;\l)$ satisfies
the hypergeometric system of differential equations (\ref{eq:hypereq}).
See also Proposition 8.3 in \cite{BS}.

Conversely, suppose $g \in \PW(C)^W$.  The  classical Paley-Wiener 
theorem ensures that $\mathcal F_A^{-1}g \in C^\infty_c (C)^W$. 
Let $\e>0$, and let $U_\e(C)$ denote
the closed $\e$-neighborhood of $C$. Then we can find $H_1,\dots,H_n \in \frak a$
so that $C \subset \cup_{j=1}^n C(H_j)^0 \subset U_\e(C)$. Here $C(H)^0$ denotes the 
interior of $C(H)$. Let $\{\chi_1,\dots,\chi_n\}$ be a smooth partition of unity  
associated with the open covering $\{C(H_1)^0,\dots,C(H_n)^0\}$ of $C$. Hence
$\chi_j\geq 0$ with $\supp \chi_j \subset  C(H_j)^0$ for all $j=1,\dots,n$, and 
$\sum_{j=1}^n \chi_j=1$ on $C$. Since $\mathcal F_A^{-1}g$ is supported in $C$,
we have
\begin{equation*}
 \mathcal F_A^{-1}g= \big(\mathcal F_A^{-1}g\big)\big(\sum_{j=1}^n \chi_j\big)=
                      \sum_{j=1}^n  \big(\mathcal F_A^{-1}g\big) \chi_j.
\end{equation*}
Therefore  $g= \sum_{j=1}^n g_j$ with
\begin{equation*}
 g_j:=\mathcal F_A\Big(\big(\mathcal F_A^{-1}g\big) \chi_j\Big) \in \PW\big(C(H_j)\cap C\big) 
\subset  \PW\big(U_\e(C)\big),
\end{equation*}
again by the classical Paley-Wiener theorem.
Observe that $U_\e(C)$ is compact, convex and $W$-invariant. 
So, by of Theorem \ref{thm:pwopdam}, 
$\pedPi{\mathcal I}g= \sum_{j=1}^n\pedPi{\mathcal I} g_j$ vanishes on $A\setminus
\exp U_\e(C)$. Since 
$\e>0$ is arbitrary, we thus conclude that $\supp \pedPi{\mathcal I}g \subset C$.
\qed

The differential operator $D_m$ allows a simple direct proof of Corollary \ref{cor:pwopdam}
when $m\in \mathcal M^+$ is even.

\medskip
\label{secondproof}
\noindent {\em Proof of Corollary \ref{cor:pwopdam} for $m \in \mathcal M^+$ even}.\;
Theorem \ref{lemma:pwinto} implies that $\pedPi{\mathcal F}$ maps $C_c^\infty(C)^W$ 
injectively into $\PW(C)^W$. Lemma \ref{lemma:suppIWinv} proves that 
$\pedPi{\mathcal I}$ maps $\PW(C)^W$ into $C_c^\infty(C)^W$. Finally, to show that 
$\pedPi{\mathcal F}$ and $k\pedPi{\mathcal I}$ are inverses to each other, see the proof
of Theorem \ref{thm:pw} in Section \ref{section:endproof}.
\qed

\section{The complex case} 
\label{section:PWcomplex}

Before proceeding with the proof of the surjectivity of $\pedtheta{\mathcal F}$ in the general 
case, let us first treat the elementary situation corresponding to the complex case.

\begin{Thm} \label{thm:complexPW}
Let $\Theta \subset \Pi$ be arbitrary, but suppose $\ma=2$ for all $\a \in \Sigma$.
Let $C \subset \pedtheta{\frak a}$ be compact, convex and $\pedtheta{W}$-invariant.
Then $\pedtheta{\mathcal I}$ maps $\pedtheta{\PW}(2;C)$ into $C^\infty_c(C)^{\pedtheta{W}}$.
\end{Thm}
\pf 
Because of the definition of $\pedtheta{\mathcal I}$ and 
Lemma \ref{lemma:Igcinfty}, we only need to prove that 
$\supp \mathcal I g \subset \exp W ( C)$.  Indeed $W (C) \cap \pedtheta{\frak a}=C$.

Recall from Example \ref{ex:phithetacomplex} that 
$\Delta(a) \pi(-\l) \pedPi{\varphi}(2;-\l,a)=\sum_{w\in W} \e(w) e^{-w\l(\log a)}$.
Therefore
\begin{align*}
  \Delta(a)(\mathcal I g)(a)&= 
\int_{i\frak a^*} g(\l) \pi(\l) \; \sum_{w\in W} \e(w) e^{-\l(\log w^{-1} a)} \; d\l\\
&=\sum_{w\in W} \e(w) \int_{i\frak a^*} g(\l) \pi(\l) e^{-\l(\log w^{-1} a)} \; d\l\\
&=\sum_{w\in W} \e(w) \left[ \mathcal F_A^{-1} (g \cdot \pi) \circ w^{-1} \right](a). 
\end{align*}
Since $\pi(\l) g(\l)$ is an entire function of exponential type $C$ and rapidly decreasing, 
the classical Paley-Wiener theorem implies that  
$\supp  \mathcal F_A^{-1} (g \cdot \pi) \subset \exp C$. Thus 
$\supp  \Delta\; \mathcal I g \subset \exp W(C)$, which prove the result since
$\mathcal I g$ is smooth.
\qed

\section{Condition A} \label{section:condA}

The proof of the surjectivity of the $\Theta$-spherical transform will depend on the following
condition on the pair $(m,\Theta)$.

\medskip
\noindent {\bf Condition A.\,} {\em Either 
\begin{equation} \label{eq:conditionAuno}
\text{$\ma\leq 2$ for all $\a \in \Sigma^+ \setminus \rootstheta^+$,} \tag{${\rm A}_1$}
\end{equation}
or there exists a $\b\in\Pi\setminus\Theta$ such that
\begin{equation} \label{eq:conditionAdue}
\text{$\Pi \setminus \Theta=\{\b\}$ and
$\inner{\b}{\a}\geq 0$ for 
all $\a \in \Sigma^+\setminus \rootstheta^+$.} \tag{${\rm A}_2$}
\end{equation}} 
\medskip

\begin{Rem}
\begin{enumerate}
\item
The requirement $\ma\leq 2$ for all $\a \in \Sigma^+ \setminus \rootstheta^+$
means that there are at most singular $\l$-hyperplanes of the form 
$\l_\a=0$ with $\a\in \Sigma^+ \setminus \rootstheta^+$. When all $m_\a$  are equal,
which happens for instance in the geometric case corresponding to 
$K_\e$ symmetric spaces with even multiplicities, this 
requirement is equivalent to one of the following conditions:    
  \begin{enumerate}
\thmlist
  \item $\Theta=\Pi$,
  \item $m=0$ (Euclidean case),
  \item $m=2$ (complex case).
\end{enumerate}
\item
In an irreducible reduced system $\Sigma$ there are at most two root lengths. Moreover, if
$\abs{\a}=\abs{\b}$, then there is $w \in W$ with $w\a=\b$ 
(see e.g. \cite{Bou}, Proposition 11, p. 164) . By $W$-invariance, multiplicity 
functions can therefore assume only two values on $\Sigma$. 

Suppose $\Sigma$ admits two root lengths, so $\Sigma$ is of type $B_l$ ($l\geq 2$), 
$C_l$ ($l \geq 2$), $F_4$ or $G_2$. 
It is then straightforward to verify that for $\abs{\Pi \setminus \Theta}=1$ the set
$\Sigma^+ \setminus \rootstheta^+$ always contains two root lengths. 
The condition $\ma\leq 2$ for all $\a \in \Sigma^+ \setminus \rootstheta^+$ 
is therefore equivalent to  $\ma\leq 2$ for all $\a \in \Sigma^+$.  
We infer that, besides (a), (b) and (c) above,
only the two following additional cases can occur:
   \begin{enumerate}
\item[(d)]  $\ma=0$ for all short $\a\in \Sigma^+$ and $\ma=2$ for all long $\a\in \Sigma^+$,
\item[(e)]  $\ma=0$ for all long $\a\in \Sigma^+$ and $\ma=2$ for all short $\a\in \Sigma^+$.
\end{enumerate}
\item
The dual cone of $\pedtheta{\frak a}$ is the $\pedtheta{W}$-invariant cone
\begin{equation}
 \label{eq:athetadual}
 \pedtheta{\frak a}^*:=\{\l \in \frak a^*:
\text{$\l(H)\geq 0$ for all $H \in \pedtheta{\frak a}$}\} 
   =\sum_{\a \in \Sigma^+ \setminus \rootstheta^+} \R_0^+ \a,
\end{equation}
where $\R_0^+$ denotes the set of nonnegative real numbers.
Condition ${\rm A}_2$ requests that $\pedtheta{\frak a}^*$
must not be ``too wide'' in $\frak a^*$, or equivalently, that $\pedtheta{\frak a}$
is ``wide enough'' in $\frak a$. It is always satisfied in the following cases:
\begin{enumerate}
\thmlist
\item $\Sigma$ of rank one;
\item $\Sigma$ corresponds to a NCC space (not necessarily with even multiplicities) 
       and $\Pi\setminus \Theta=\{\gamma_1\}$,
       where $\gamma_1$ is the unique simple noncompact root.
  \end{enumerate}
The latter statement depends on the fact that  $\gamma_1+\a$ is never a root when 
$\a$ is a noncompact positive root because $(\gamma_1+\a)(X_0)=2$ for $X_0$ as in 
(\ref{eq:signatureNCC}). Thus $\inner{\gamma_1}{\a}\geq 0$ by \cite{Bou}, 
Corollary to Theorem 1, p. 162.

The complete list of the cases in which 
Condion ${\rm A}_2$ is satisfied can be obtained using classification of the irreducible
reduced root systems together with the property that $\inner{\a}{\b} \geq 0$ 
for two non-proportional roots for which $\a+\b$ is not a root 
(see e.g. \cite{Bou}, Ch. VI, Corollary to Theorem 1). 
For non-exceptional root systems, we get the following list, given in the notation of 
\cite{Bou}, Ch. VI, \S 4. If $\Sigma$ is of type $A_l$ ($l\geq 1$), then Condition ${\rm A}_2$ 
holds for $\Pi \setminus \Theta=\{\b\}$ for every choice of $\b\in \Pi$. The same is true 
for $\Sigma$ of type $B_2=C_2$ or $D_3$ or $D_4$. 
If $\Sigma$ is of type $B_l$ (with $l>2$) or $C_l$ (with $l>2$) or $D_l$ (with $l>4$), 
then Condition ${\rm A}_2$ holds only for $\Pi \setminus \Theta=\{\a_l\}$. 
Observe, though, that the root system $G_2$ does not satisfy Condition ${\rm A}_2$ 
for any choice of the simple root $\beta$. 
\item
Condition A is never satisfied by $\Theta=\emptyset$ when 
$\dim \frak a>1$. Hence our Paley-Wiener theorem applies to the transform associated
with Harish-Chandra series only in the even-multiplicity rank-one case. Notice however 
that the Harish-Chandra series agree with geometrically defined spherical functions only
in the rank-one case. 
\item 
It is important to observe that, as a consequence of (1) and (3) above together with the list
in Appendix \ref{app:list}, 
all pairs $(\Theta,m)$ corresponding to $K_\e$ symmetric spaces with even multiplicities 
satisfy Condition A. Thus Condition A applies to 
all geometrical situations in which the theory of the $\Theta$-spherical functions might 
be relevant. 
\end{enumerate}
\end{Rem}

In Sections \ref{section:PWopdam} and \ref{section:PWcomplex} 
we have already given 
direct proofs of the Paley-Wiener theorem for $\Theta=\Pi$ 
and in the complex case.
The proof in the latter case can be easily modified to apply also 
to the situation in which all 
multiplicities are either $0$ or $2$. 
Indeed, the shift operator $D_2$ to be considered in this 
case is obtained as in Example \ref{ex:shiftcomplex} with the product over 
$\Sigma^+$ replaced by a product over the positive roots with multiplicity 2.
Thus all cases in which Condition ${\rm A}_1$ holds can be 
easily treated directly.
For the proof of the surjectivity of the $\Theta$-spherical 
transform onto the space 
$\pedtheta{\PW}(m;C)$ we shall therefore 
restrict ourselves only to the case in which Condition 
${\rm A}_2$ is satisfied. 
This simplifies some technical aspects of the argument. 
Nevertheless a unified proof for the pairs 
$(\Theta,m)$ satisfying Condition A is possible, 
and is based on the observation in 
Corollary \ref{cor:shiftnonsing}. 
We shall indicate in Remark \ref{rem:generalproof} 
how to proceed in the general case.
 
\section{Away from the walls of $\pedtheta{\frak a}$}
\label{section:walls}

In the following we assume $\Theta \subsetneq  \Pi$.  Set
\begin{align}
  \label{eq:CTheta}
  \pedtheta{C}&:=\{H \in \frak a: \text{$\a(H)=0$ for all $\a \in \Theta$ and 
                                   $\a(H)>0$ for all $\a \in \Pi \setminus\Theta$}\} \notag\\
                      &=\{H \in \frak a: \text{$\a(H)=0$ for all $\a \in \rootstheta$ and 
                                   $\a(H)>0$ for all $\a \in \Sigma^+ \setminus\rootstheta^+$}\}. 
\end{align}
Then $\pedtheta{C}$ is a closed convex subset of $\overline{\frak a^+}$ with 
$\dim \big(\Span_\R \pedtheta{C}\big)=\abs{\Pi\setminus\Theta}$. 
Moreover, $\overline{\frak a^+}=\overset{\cdot}{\cup}_{\Theta \subset \Pi} \pedtheta{C}$. 
We refer to \cite{Hu}, pp. 25-26, for the proof of these statements.  
The group $\pedtheta{W}$ stabilizes $\pedtheta{C}$. More precisely, the following lemma
holds.

\begin{Lemma}[\cite{Hu}, Proposition 1.15]\label{lemma:WeylonCTheta}
The following properties are equivalent for $w \in W$:
\begin{enumerate}
\thmlist
\item $w \in \pedtheta{W}$,
\item $w\pedtheta{C}=\pedtheta{C}$,
\item $wH=H$ for all $H \in \pedtheta{C}$.
\end{enumerate}
\end{Lemma}

Observe that $\pedPi{C}=\{0\}$ and $\pedempty{C}=\frak a^+$. 
Let $\R^+$ denote the set of positive real numbers.
If $\abs{\Pi\setminus\Theta}=1$, then $\pedtheta{C}=\R^+ X^0$, for some 
$X^0 \in \pedtheta{\frak a}$. In particular, if $\Sigma$ is the
root system of a NCC symmetric space and $\Theta=\Pi_0$ is the set of compact
positive simple roots, then $\pedtheta{C}=\R^+ X^0$, where $X^0 \in \frak a$ 
is the cone generating element for the causal structure. Finally, observe
from (\ref{eq:CTheta}) that $\pedtheta{C} \subset \pedtheta{\frak a}$.

We shall be interested in the closed cones in $\pedtheta{\frak a}$
\begin{equation}
  \label{eq:CrX0}
C(r,X^0):=  rX^0 + \overline{\pedtheta{\frak a}}\,,
\end{equation}
where $X^0 \in \pedtheta{C}$ and $r>0$.

\begin{Lemma} \label{lemma:anotinC}
Let $r>0$ and $X^0 \in \pedtheta{C}$. Then $C(r,X^0)$ is a closed 
$\pedtheta{W}$-invariant cone in $\pedtheta{\frak a}$.
If $a \in A^+ \setminus \exp C(r,X^0)$, then
there exists $\beta \in \Pi \setminus \Theta$ with $\beta(\log a-r X^0)<0$.
\end{Lemma}
\pf
The first statement follows immediately from (\ref{eq:CrX0}) and Lemma 
\ref{lemma:WeylonCTheta}.

Suppose now $a \in A^+ \setminus \exp C(r,X^0)$. Then 
$$\log a-r X^0 
\notin \overline{\pedtheta{\frak a}}=\{H \in \frak a: \text{$\a(H)\geq 0$ for all 
$\a \in \Sigma^+ \setminus \rootstheta^+$}\}$$
Hence, there exists a $\beta_1 \in \Sigma^+ \setminus \rootstheta^+$ with 
$\beta_1(\log a-r X^0)<0$. Write $\beta_1=\sum_{\a\in \Pi\setminus \Theta} n_\a \a+
\sum_{\gamma\in \Theta} n_\gamma \gamma$ with $n_\a,n_\gamma\geq 0$. 
Then $(\Pi \setminus\Theta)_1:=\{\a \in \Pi \setminus \Theta:n_\a>0\}$ is nonempty.
Set $$\beta_2:=\sum_{\a\in (\Pi\setminus \Theta)_1} n_\a \a =
\beta_1-\sum_{\gamma\in \Theta} n_\gamma \gamma.$$
Then 
\begin{align*}
  \beta_2(\log a-r X^0)&=\beta_1(\log a-r X^0)-\sum_{\gamma\in \Theta} n_\gamma 
[\gamma(\log a)-r \gamma(X^0)]\\
&=\beta_1(\log a-r X^0)-\sum_{\gamma\in \Theta} n_\gamma \gamma(\log a) <0.
\end{align*}
Since $n_\a>0$ for all $\a\in (\Pi\setminus \Theta)_1$, there must be 
$\beta \in (\Pi\setminus \Theta)_1$ satisfying $\beta(\log a-r X^0)<0$, as requested.
\qed

Recall the definition of the dual cone $ \pedtheta{\frak a}^*$ of $\pedtheta{\frak a}$ 
from (\ref{eq:athetadual}). 
Its interior is  $\big(\pedtheta{\frak a}^*\big)^0=\{\l \in \frak a^*: \text{$\l(H)>0$ for all 
$H \in \pedtheta{\frak a}$}\}.$

\begin{Lemma} \label{lemma:intathetastar}
Suppose $\abs{\Pi \setminus \Theta}=1$. 
Then $\big(\pedtheta{\frak a}^*\big)^0\neq \emptyset$.
\end{Lemma}
\pf
Since $\abs{\Pi\setminus\Theta}=1$ there is $X^0\neq 0$ so that
$\pedtheta{C}=\R_+ X^0$. We claim that 
$\inner{X^0}{H}>0$ for all $H \in \pedtheta{\frak a}$.   
Indeed, suppose first $H \in \overline{\pedtheta{\frak a}}$. Then 
there is $H_1 \in \overline{\frak a^+}$ and $w \in \pedtheta{W}$ so that
$H=wH_1$. By Lemma \ref{lemma:WeylonCTheta} and Lemma A in \cite{Eber}, p. 197, 
we have $\inner{X^0}{H}=\inner{X^0}{H_1} \geq 0$.
Hence $\overline{\pedtheta{\frak a}} \subset \{H \in \frak a: \inner{X^0}{H}\geq 0\}$, which 
implies  $\pedtheta{\frak a} \subset \{H \in \frak a: \inner{X^0}{H}> 0\}$. 

Let $\l_0 \in \frak a$ correspond to $X^0$ under the identification of 
$\frak a$ and $\frak a^*$ under the inner product, i.e. $X^0=A_{\l_0}$. Then 
$0 \neq \l_0 \in  \big( \pedtheta{\frak a}^*\big)^0$.
\qed

For an even multiplicity function $m \in \mathcal M^+$ we set 
\begin{equation}
  \label{eq:athetadualm}
  \pedtheta{\frak a}^*(m)=\{\l \in \pedtheta{\frak a}^*: 
\text{$\la\geq \ma/2$ for all $\a \in \Sigma^+ \setminus \rootstheta^+$}\}. 
\end{equation}
The geometrical properties of $\pedtheta{\frak a}^*(m)$ are collected 
in the following lemma.

\begin{Lemma}
The set  $\pedtheta{\frak a}^*(m)$ is closed, convex and  $\pedtheta{W}$-invariant.  
If $\mu\in \pedtheta{\frak a}^*(m)$ and $t\geq 1$, then also $t\mu \in \pedtheta{\frak a}^*(m)$. 
Furthermore, if $\abs{\Pi \setminus \Theta}=1$, then $\pedtheta{\frak a}^*(m)$ has
non-empty interior.
\end{Lemma}
\pf
The only property needing some proof is that $\pedtheta{\frak a}^*(m)$ has non-empty 
interior when $\abs{\Pi \setminus \Theta}=1$.
Set $$\big( \pedtheta{\frak a}^*\big)^\dagger:=
\{\l \in \frak a^*: \text{$\inner{\l}{\a}>0$ for all 
$\a \in\Sigma^+ \setminus \rootstheta^+$}\}.$$
Then it suffices to show that $\big( \pedtheta{\frak a}^*\big)^\dagger 
\cap \big(\pedtheta{\frak a}^*\big)^0 \not=\emptyset$. Indeed, if this is
the case, then there is $0 \neq \mu \in
 \big( \pedtheta{\frak a}^*\big)^\dagger 
\cap (\pedtheta{\frak a}^*)^0$.
If $t \geq \max_{\a\in\Sigma^+\setminus \rootstheta^+} \frac{\ma}{2} \; 
\frac{\inner{\a}{\a}}{\inner{\mu}{\a}}$, then $t\mu \in \big(\pedtheta{\frak a}^*\big)^0$
and $(t\mu)_\a=t \, \frac{\inner{\mu}{\a}}{\inner{\a}{\a}} \geq\frac{\ma}{2}$.

Recall the notation $A_\l$ for the element in $\frak a$ identified with $\l \in \frak a^*$
via the inner product. In particular, as in the proof of Lemma \ref{lemma:intathetastar},
let $A_{\l_0}=X^0$.  Observe that
\begin{align*}
 \big( \pedtheta{\frak a}^*\big)^\dagger&=
   \{\l \in \frak a^*:\text{$\inner{\l}{\nu}>0$ for all $\nu \in \pedtheta{\frak a}^*$}\}\\
  &=\{\l \in \frak a^*: A_\l \in (\pedtheta{\frak a}^{**})^0=\pedtheta{\frak a}\}.
\end{align*}
Thus, because of Lemma \ref{lemma:intathetastar}, we conclude
$\R^+ \l_0 
\subset  \big( \pedtheta{\frak a}^*\big)^\dagger \cap (\pedtheta{\frak a}^*)^0$.
\qed

The set $\pedtheta{\frak a}^*(m)$ is introduced because it is a ``large'' closed
subset of $\pedtheta{\frak a}^*$ which is ``away'' from the possible singularities of 
every $g \in \pedtheta{\PW}(m;C)$. 
This is made precise by the following lemma.

\begin{Lemma} \label{lemma:estgshift}
Let $C\subset \frak a$ be compact and convex, and let
$g:\frakacs\to \C$ satisfy $\pedtheta{e}^-(m;\l)g(\l) \in PW(C)$. 
Then $g$ is holomorphic on a neighborhood of the 
convex set $i\frak a^*-\pedtheta{\frak a}^*(m)$.
Furthermore, for every $N \in \N$, there is a 
constant $C_N>0$ such that for all 
$\l \in i\frak a^*$ and $\mu \in -\pedtheta{\frak a}^*(m)$
\begin{equation}
  \label{eq:estgaway}
  \abs{g(\l+\mu)} \leq C_N(1+\abs{\l})^{-N} e^{q_C(\mu)}. 
\end{equation}
\end{Lemma}
\pf 
The function $g$ is holomorphic on the open set 
$$i\frak a^* -\{\l \in \frak a^*:\text{$\l_\a >\ma/2-1$ for all
$\a \in \Sigma^+ \setminus \rootstheta^+$}\}.$$
To prove the estimate, observe that,  
if $\l \in i\frak a^*$ and $\mu\in -\pedtheta{\frak a}^*(m)$, then 
\begin{equation*}
\abs{\pedtheta{e}^-(m;\l+\mu)}\geq 
\prod_{\a \in \Sigma^+ \setminus \rootstheta^+} \;
\prod_{k=-\ma/2+1}^{\ma/2-1} \abs{\mu_\a-k} \geq 1.
\end{equation*}
\qed

\begin{Lemma} \label{lemma:intassum}
Let $C\subset \pedtheta{\frak a}$ be compact, convex and $\pedtheta{W}$-invariant.
Assume that $g:\frakacs\to \C$ satisfy $\pedtheta{e}^-(m;\l)g(\l) \in PW(C)^{\pedtheta{W}}$. 
Then, for all $\mu \in -\pedtheta{\frak a}^*(m)$ and $a \in A$,
\begin{equation*}
  \Delta(m;a) \int_{i\frak a^*} g(\l) \pedPi{\varphi}(m;-\l,a) \; \frac{d\l}{\abs{\pedPi{c}^+(m;\l)}^2}
 =\sum_{w\in W} D_m \int_{i\frak a^*} g(\l+\mu) e^{-w(\l+\mu)(\log a)} \; d\l\,.
\end{equation*}
Moreover, for all $a \in A$ and $w \in W$, the integral 
\begin{equation} \label{eq:integral}
   \int_{i\frak a^*} g(\l+\mu) e^{-w(\l+\mu)(\log a)} \; d\l
\end{equation}
is independent of $\mu \in -\pedtheta{\frak a}^*(m)$. If $w \in \pedtheta{W}$, then 
(\ref{eq:integral}) is $\pedtheta{W}$-invariant. In particular, 
\begin{equation*} 
   \int_{i\frak a^*} g(\l+\mu) e^{-w(\l+\mu)(\log a)} \; d\l=
   \int_{i\frak a^*} g(\l+\mu) e^{-(\l+\mu)(\log a)} \; d\l.
\end{equation*}
  \end{Lemma}
\pf
Recall formula (\ref{eq:formulabis}) for the function $\pedPi{\varphi}$.
Since 
$$g(\l) \pedPi{\varphi}(m;-\l,a) \abs{\pedPi{c}^+(m;\l)}^{-2}=
g(\l)\pi(\l)\pedPi{e}^+(m;\l) \pedPi{\varphi}(m;-\l,a)$$ 
is entire and rapidly decreasing, we can apply
Cauchy's theorem and get for all $\mu \in -\pedtheta{\frak a}^*(m)$
\begin{align*}
   &\Delta(m;a) \int_{i\frak a^*} 
g(\l) \pedPi{\varphi}(m;-\l,a) \; \frac{d\l}{\abs{\pedPi{c}^+(m;\l)}^2}\\
&=\Delta(m;a) \int_{i\frak a^*} g(\l+\mu)\pi(\l+\mu) \pedPi{e}^+(m;\l+\mu) 
\pedPi{\varphi}(m;-\l-\mu,a) \; d\l\\
&=\int_{i\frak a^*} g(\l+\mu) \sum_{w\in W} D_m e^{-w(\l+\mu)(\log a)} \; d\l\\
&=\sum_{w\in W} D_m \int_{i\frak a^*}  g(\l+\mu) e^{-w(\l+\mu)(\log a)} \; d\l\,.
\end{align*}
The last equality is justified by the fact that $D_m$ is a differenital operator with
smooth coefficients on $A$ and by Lemma \ref{lemma:estgshift}.

Set 
$$F_w(a,\mu):=\int_{i\frak a^*} g(\l+\mu) e^{-w(\l+\mu)(\log a)} \; d\l\,.$$ 
It is independent of $\mu \in -\pedtheta{\frak a}^*(m)$ by Cauchy's theorem, since 
its integrand is holomorphic in a simply connected open neighborhood of 
$i\frak a^* -\pedtheta{\frak a}^*(m)$ on which an estimate of the form 
(\ref{eq:estgaway}) holds. 
For $w,v \in \pedtheta{W}$ 
\begin{align*}
  F_w(va,\mu)&=\int_{i\frak a^*} g(\l+\mu) e^{-v^{-1}w(\l+\mu)(\log a)} \; d\l\\
&=\int_{i\frak a^*} g(w^{-1}vw\l+\mu) e^{-(w\l+v^{-1}w\mu)(\log a)} \; d\l\\
&=\int_{i\frak a^*} g(\l+w^{-1}v^{-1}w\mu) e^{-(w\l+v^{-1}w\mu)(\log a)} \; d\l \qquad
   \text{(by $\pedtheta{W}$-invariance of $g$)}\\
&=F_w(a,w^{-1}v^{-1}w\mu)\\
&=F_w(a,\mu)  \qquad \text{(because $w^{-1}v^{-1}w\mu \in-\pedtheta{\frak a}^*(m)$)}.
\end{align*}
Since $F_w(a,\mu)=F_e(w^{-1}a,\mu)$, it follows in particular that 
$F_w(a,\mu)=F_e(a,\mu)$ for $w \in \pedtheta{W}$.
\qed

The possibility of separating points in $\pedtheta{\frak a}$ from points outside $\pedtheta{\frak a}$
 is guaranteed by the following lemma. 
Recall from (\ref{eq:apthetaW}) the definition of  $\aptheta{W}$ and length function 
$\ell$.

\begin{Lemma} \label{lemma:apthetaminus}
If $u \in \aptheta{W} \setminus \{e\}$, then there is $\b \in \Pi \setminus \Theta$ with 
$u\b \in -\Sigma^+$.
\end{Lemma}
\pf
If $u\neq e$, then $\ell(u)>1$. Since $\ell(u)$ represents the number of elements in 
$\Sigma^+$ mapped by $u$ into $-\Sigma^+$, there must be $\beta \in \Pi$ so that 
$u\b\in -\Sigma^+$. But this implies $\ell(ur_\b)<\ell(u)$ (see e.g. \cite{Hu}, Lemma 1.6).
Thus $u\b \notin \rootstheta$ because of the definition of $\aptheta{W}$.
\qed 

In the following we assume that $\Theta$ satisfies Condition ${\rm A}_2$. 
This assumption is crucial because it will allow us either 
to separate points in $\pedtheta{\frak a}$ from point outside it, or to move the contour 
of integration without encountering singular hyperplanes.

\begin{Lemma} \label{lemma:betaplusastarminus}
Suppose that Condition ${\rm A}_2$ is satisfied by $\Theta$, and let $\beta \in 
\Pi \setminus \Theta$.  If $\mu_0 \in \pedtheta{\frak a}^*(m)$, 
then $t\beta+\mu_0 \in \pedtheta{\frak a}^*(m)$ for all $t\geq 0$.
\end{Lemma}
\pf 
Since $\pedtheta{\frak a}^*$ is a closed cone, the assumption 
$\beta, \mu_0 \in \pedtheta{\frak a}^*$ implies
$t\beta+\mu_0 \in \pedtheta{\frak a}^*$ for all $t\geq 0$. 
Consider now $\a \in \Sigma^+\setminus \rootstheta^+$. 
Then, assuming Condition ${\rm A}_2$, we have 
$$(t\beta + \mu_0)_\a=t \b_\a + (\mu_0)_\a\geq \ma/2.$$
\qed

\begin{Lemma} \label{lemma:intzero}
Suppose Condition ${\rm A}_2$ is satisfied. Let $a \in A^+$, 
$\mu \in -\pedtheta{\frak a}^*(m)$, $w \notin \pedtheta{W}$, and let
$g:\frakacs\to \C$ satisfy $\pedtheta{e}^-(m;\l)g(\l) \in PW(C)^{\pedtheta{W}}$.  
Then 
\begin{equation}
  \label{eq:intzero}
 \int_{i\frak a^*} g(\l+\mu) e^{-w(\l+\mu)(\log a)}\; d\l=0. 
\end{equation}
\end{Lemma}
\pf
Write $w=uv$ with $v \in \pedtheta{W}$ and $u \in \aptheta{W}$. 
If $w \notin \pedtheta{W}$, then $u \neq e$. By Lemma \ref{lemma:apthetaminus},
there is $\beta \in \Pi \setminus \Theta$ so that $\beta u \in -\Sigma^+$. 
Let $\mu_0 \in \pedtheta{\frak a}^*(m)$. Then 
$t\beta+\mu_0 \in \pedtheta{\frak a}^*(m)$ for all $t\geq 1$
by Lemma \ref{lemma:betaplusastarminus}. Set $\mu_t:=-(t\beta+\mu_0).$
By $\pedtheta{W}$-invariance of $g$ and $\pedtheta{\frak a}^*(m)$,
 and since (\ref{eq:integral}) is independent 
of $\mu \in -\pedtheta{\frak a}^*(m)$, we have
\begin{align*}
  \int_{i\frak a^*} g(\l+\mu) e^{-w(\l+\mu)(\log a)}\; d\l&=
   \int_{i\frak a^*} g(\l+\mu_t) e^{-w(\l+\mu_t)(\log a)}\; d\l \\
  &= \int_{i\frak a^*} g(v\l+v\mu_t) e^{-u(v\l+v\mu_t)(\log a)}\; d\l \\
  &= \int_{i\frak a^*} g(\l+v\mu_t) e^{-u(\l+v\mu_t)(\log a)}\; d\l \\
  &= \int_{i\frak a^*} g(\l+\mu_t) e^{-u(\l+\mu_t)(\log a)}\; d\l.
\end{align*}
Consider the last integral. By Lemma \ref{lemma:estgshift} we have for all $N \in \N$:
\begin{equation*}
   \abs{g(\l+\mu_t)} \leq C_N(1+\abs{\l})^{-N} e^{q_C(\mu_t)}=
       C_N e^{q_C(-\mu_0)} (1+\abs{\l})^{-N} e^{tq_C(-\b)}
\end{equation*}
with $q_C(-\b)\leq 0$ because $\b \in \Pi \setminus\Theta$ and $C \subset \pedtheta{\frak a}$.
Moreover
\begin{equation*}
\abs{e^{-u(\l+\mu_t)(\log a)}}=\abs{e^{-u\mu_t(\log a)}}=e^{u\mu_0(\log a)}
e^{tu\b(\log a)}.
\end{equation*}
with $u\b(\log a)<0$ because $u\b \in -\Sigma^+$ and $a \in A^+$. 
The result then follows by taking the limit $t\to+\infty$.
\qed 

\begin{Thm} \label{thm:newinv}
Suppose that Condition ${\rm A}_2$ is satisfied. 
Let $C\subset \pedtheta{\frak a}$ be compact and convex, and let
$g:\frakacs\to \C$ satisfy $\pedtheta{e}^-(m;\l)g(\l) \in \PW(C)^{\pedtheta{W}}$.
Suppose furthermore that $\mu \in -\pedtheta{\frak a}^*(m)$.
Then for  all $a \in \pedtheta{A}$
\begin{equation} \label{eq:newinv}
    \Delta(m;a) \int_{i\frak a^*} g(\l) \pedPi{\varphi}(m;-\l,a) \; \frac{d\l}{\abs{\pedPi{c}^+(m;\l)}^2}
= \abs{\pedtheta{W}} D_m \int_{i\frak a^*} g(\l+\mu) e^{-(\l+\mu)(\log a)}\; d\l\,.
\end{equation}
  \end{Thm}
\pf 
Let $a \in A^+$. By Lemmas \ref{lemma:intassum} and  \ref{lemma:intzero}
\begin{align*}
\Delta(m;a) \int_{i\frak a^*} g(\l) \pedPi{\varphi}(m;-\l,a) \; \frac{d\l}{\abs{\pedPi{c}^+(m;\l)}^2}
&=\sum_{w\in W} D_m \int_{i\frak a^*} g(\l+\mu) e^{-w(\l+\mu)(\log a)} \; d\l\\
&=\sum_{w\in \pedtheta{W}} D_m \int_{i\frak a^*} g(\l+\mu) e^{-w(\l+\mu)(\log a)} \; d\l\\
&=\abs{\pedtheta{W}} D_m \int_{i\frak a^*} g(\l+\mu) e^{-(\l+\mu)(\log a)}\; d\l\,.
\end{align*}
Both sides of (\ref{eq:newinv}) are $\pedtheta{W}$-invariant and continuous.
Hence they agree on the entire $\pedtheta{A}$.
\qed 

As an immediate consequence we obtain an alternative expression for the
inversion formula for the $\Theta$-spherical transform. It is 
a modification by the differential operator $D_m$ of the
inversion formula for the Euclidean Fourier transform. 
\begin{Cor} \label{cor:newinv}
Suppose that Condition ${\rm A}_2$ is satisfied. Let $f\in C^\infty_c(C)^{\pedtheta{W}}$. Then  
for  all $a \in \pedtheta{A}$ and  $\mu \in -\pedtheta{\frak a}^*(m)$.
\begin{equation} \label{eq:newinvf}
    \Delta(m;a) f(a)
= (-1)^{d(\Theta,m)}  k\; \abs{W}^2 \;
D_m \int_{i\frak a^*} (\pedtheta{\mathcal F}f)(m;\l+\mu) e^{-(\l+\mu)(\log a)}\; d\l\,,
\end{equation}
where $k$ is the normalizing constant in Theorem  \ref{thm:inversion} and 
$d(\Theta,m)$ is as in (\ref{eq:dmtheta}).
  \end{Cor}
\pf Immediate from Theorems \ref{thm:inversion}  and \ref{thm:newinv}.
\qed 

Recall the notation $C(r,X^0)=rX^0+\overline{\pedtheta{\frak a}}\,$ from (\ref{eq:CrX0}).

\begin{Lemma} \label{lemma:suppIgCrX}
Suppose Condition ${\rm A}_2$ is satisfied. Let $r>0$ and $X^0 \in \pedtheta{C}$.
Let $C\subset \pedtheta{\frak a}$ be compact and convex, and let
$g:\frakacs\to \C$ satisfy $\pedtheta{e}^-(m;\l)g(\l) \in PW(C)^{\pedtheta{W}}$.
 Assume furthermore that $C \subset C(r,X^0)$. Then 
$\supp \pedtheta{\mathcal I} g \subset C(r,X^0)$.
\end{Lemma}
\pf
Let $a \in A^+ \setminus \exp C(r,X^0)$, and let 
$\beta \in \Pi\setminus \Theta$ be as in Lemma
\ref{lemma:anotinC}. For a fixed $\mu_0 \in \pedtheta{\frak a}^*(m)$ and  $t\geq 1$ we set
$\mu_t:=-(t\beta+\mu_0)$. As in Lemma \ref{lemma:intzero}, we have
$\mu_t \in -\pedtheta{\frak a}^*(m)$ and, for all $N \in \N$,
\begin{equation*}
   \abs{g(\l+\mu_t)}e^{-(\l+\mu_t)(\log a)} \leq 
         C_N e^{q_C(-\mu_0)} e^{\mu_0(\log a)} 
         (1+\abs{\l})^{-N} e^{tq_C(-\b)}e^{t\b(\log a)}.
\end{equation*}
Notice that $C\subset C(r,X^0)$ implies $C=rX^0+C'$, where 
$C':=C-rX^0 \subset \overline{\pedtheta{\frak a}^+}$ is compact,
convex, $\pedtheta{W}$-invariant. Hence
\begin{equation*}
  q_C(-\b)=\sup_{H \in C} [-\beta(H)]=-\b(rX^0)+ \sup_{H \in C'} [-\beta(H)]=
                  -\b(rX^0)+ q_{C'}(-\b)                         
\end{equation*}
with $q_{C'}(-\b)\leq 0$ because $\b \in \Pi \setminus\Theta$.
Thus 
\begin{equation*}
   \abs{g(\l+\mu_t)}e^{-(\l+\mu_t)(\log a)} \leq C_N e^{\mu_0)(\log a)}  
         (1+\abs{\l})^{-N} e^{t\b(\log a-rX^0)},
\end{equation*}
which converges to $0$ as $t\to +\infty$ by Lemma \ref{lemma:anotinC}.
Together with Theorem \ref{thm:newinv}, this proves that for 
$a \in A^+ \setminus \exp C(r,X^0)$
we have $\Delta(m;a)(\mathcal Ig)(a)=0$. 
Since $\mathcal Ig$ is $\pedtheta{W}$-invariant and 
smooth on $A$ by Lemma \ref{lemma:Igcinfty}, we conclude that 
$\mathcal Ig(a)=0$ for $a \in A \setminus \cup_{w \in W} \exp w\big(C(r,X^0)\big)=
A \setminus \cup_{u \in \aptheta{W}} \exp u\big(C(r,X^0)\big)$.
Because of Lemma \ref{lemma:apthetaWonatheta}, the sets $u(C(r,X^0))$ are
pairwise disjoint. Thus $\supp \pedtheta{\mathcal I}g \subset C(r,X^0)$.
\qed

\begin{Rem} \label{rem:generalproof}
 One could use Corollary \ref{cor:division} and replace the set
$\pedtheta{\frak a}^*(m)$ with
\begin{equation*}
  \{\l \in \pedtheta{\frak a}^*: 
\text{$\la\geq \ma/2$ for all $\a \in \Sigma^+ \setminus 
\rootstheta^+$ with $\ma>2$}\}. 
\end{equation*}
This would allows us a proof of the previous lemmas 
which applies also to the complex case (which generally does not satisfy condition $A_2$). 
\end{Rem}

The results obtained so far in this section apply in particular to every element $g \in 
 \pedtheta{\PW}(m;C)$, but  the additional condition that 
$\Pavtheta g$ extends to be entire on $\frakacs$ has been not used.
This condition is needed in the proof of the Paley-Wiener theorem 
to apply Lemma \ref{lemma:suppIgCrX}, in particular to know that $\pedtheta{\mathcal I}g$
has compact support.

At this point we can already prove that $\pedtheta{\mathcal I}$ maps
$\pedtheta{\PW}(m;C)$ into $C^\infty_c(C)^{\pedtheta{W}}$, but only for
a special class of compact convex and  $\pedtheta{W}$-invariant exausting subsets $C$ of 
$\pedtheta{\frak a}$. 

Recall for $H\in \frak a$ the notation $C(H):=\conv (W( H))$.

\begin{Cor}
Let $C \subset \pedtheta{\frak a}$ be a compact, convex 
and $\pedtheta{W}$-invariant subset of the form $C= C(r,X^0) \cap C(H)$ with 
$r>0, X^0 \in \pedtheta{C}$ and $H \in \frak a$. 
Suppose $g \in \pedtheta{\PW}(m;C)$. 
Assume Condition ${\rm A}_2$ is satisfied.
Then $\supp \pedtheta{\mathcal I}g \subset \exp C$. 
\end{Cor}
\pf   
Lemmas \ref{lemma:Igcinfty}  and \ref{lemma:suppIgCrX} ensure that 
$\pedtheta{\mathcal I}g$ is smooth and compactly 
supported inside $C(H)$.  
Moreover, Lemma \ref{lemma:suppIgCrX} guarantees that  $\supp \pedtheta{\mathcal I}
\subset C(r,X^0)$.
\qed

\section{Support properties}

By means of Condition $A_2$ we could prove in Section \ref{section:walls} that 
$\pedtheta{\mathcal I}g$ is compactly supported inside $\pedtheta{A}$ for every
$g \in \pedtheta{\PW}(m;C)$, where $C$ is an arbitrary compact convex and 
$\pedtheta{W}$-invariant subset of $\pedtheta{\frak a}$. In this section we apply 
a support theorem proven in \cite{OP4} to show that the support of $\pedtheta{\mathcal I}g$
is indeed contained in $\exp C$.

Recall the $W$-invariant polynomial  $q$ introduced in  (\ref{eq:q}).
It satisfies $q(-\l)=q(\l)$ for all $\l \in \frakacs$, 
and $qg \in \PW(C)^{\pedtheta{W}}$ for all 
$g \in \pedtheta{\PW}(m;C)$.
Since $q\in\polya^W$, we can consider the corresponding operator 
$D(m;q) \in \D(\frak a,\Sigma,m)$, that is 
\begin{equation} \label{eq:Dmeexplicit}
  D(m;q) =\Upsilon \left( \prod_{\a \in \Sigma} \; \prod_{k=-\ma/2+1}^{\ma/2-1}
 \big( T(m;H_{\a}/2)-k\big) \right)
\end{equation}
(see Definition \ref{def:cherednik} and (\ref{eq:diffops})). 
Hence
\begin{equation}
  \label{eq:dmevarphi}
D(m;q) \pedtheta{\varphi}(m;\l,a)=q(\l) \pedtheta{\varphi}(m;\l,a) 
\end{equation}
for all $a \in A^+$ and all $\Theta \subset \Pi$. Moreover, 
there is $k \in \N$ such that
\begin{equation} \label{eq:De}
D_q:=\Delta^k D(m;q)
\end{equation}
is a $W$-invariant differential operator on $A$ with analytic coefficients. 

\begin{Lemma} \label{lemma:keyequality}
Suppose $g \in \pedtheta{\PW}(m,C)$. 
Then, for all $a \in \pedtheta{A}$,
\begin{equation} 
  D_q\big(\pedtheta{\mathcal I}g\big)(a)=
\abs{\pedtheta{W}}\, \Delta(k-m;a) D_m \mathcal F^{-1}_A(qg)(a)
\end{equation}
with $\Delta(k-m;a) D_m \mathcal F^{-1}_A(qg)  \in C^\infty_c(C)$. 
\end{Lemma}
\pf
Because of (\ref{eq:dmevarphi}) and Theorem  \ref{thm:newinv}, for all $a \in A$ we have
\begin{align*}
  \Delta(m;a)  D_q \big(\pedtheta{\mathcal I}g\big)(a)&=
           \Delta(m;a) \Delta^k(a)  \int_{i\frak a^*} g(\l)q(\l) \pedPi{\varphi}(m;-\l,a) \;
\frac{d\l}{\abs{\pedPi{c}^+(m;-\l)}^2}\\
         &=\abs{\pedtheta{W}} \Delta^k(a) 
               D_m \int_{i\frak a^*} q(\l+\mu)g(\l+\mu)e^{-(\l+\mu)(\log a)} \;d\l.
\end{align*}
Since $qg \in \PW(C)$, we can shift the contour of integration and obtain
\begin{align} \label{eq:partialeqDe}
   \Delta(m;a)  D_q \big(\pedtheta{\mathcal I}g\big)(a)&=
   \abs{\pedtheta{W}} \Delta^k(a)  D_m \int_{i\frak a^*} q(\l)g(\l) e^{-\l(\log a)} \;d\l \notag\\
   &=\abs{\pedtheta{W}} \Delta^k(a) D_m \mathcal F_A^{-1}(qg) (a).
\end{align}
The classical Paley-Wiener theorem guarantees that 
$\supp \pedtheta{\mathcal F}^{-1}(qg) \subset C$. 
Since $D_q\big(\pedtheta{\mathcal I}g\big)$ is smooth, we conclude that 
$\Delta(m;a)$ must divide the right-hand side of (\ref{eq:partialeqDe}). Thus
$$
\Delta(m-k) D_m \mathcal F^{-1}_A(qg) =
\frac{\Delta(m)D_m \mathcal F^{-1}_A(qg)}{\Delta^k(a)}
$$
is smooth and supported in $C$. This proves the lemma.
\qed

Lemma \ref{lemma:keyequality} shows that $D_q\big(\pedtheta{\mathcal I}g\big)$
is supported inside $\exp C$. We can therefore conclude 
that $\pedtheta{\mathcal I}$ maps $\pedtheta{\PW}(m;C)$
to $C_c^\infty(C)^{\pedtheta{W}}$ if we prove that the support of  $\pedtheta{\mathcal I}g$
is contained in $\exp C$ when the support of 
$D_q\big(\pedtheta{\mathcal I}g\big)$ is in $\exp C$.
This is guaranteed by the following proposition, which depends on the 
particular form of the principal symbol of the differential operator $D_q$.

\begin{Prop} \label{prop:key}
Let  $D_q$ be the differential operator from (\ref{eq:De}). Then
\begin{equation*}
  \supp f \subset \exp C \iff \supp D_q f \subset \exp C
\end{equation*}
for every $f \in C_c^\infty(A)^{\pedtheta{W}}$ and for every 
 compact, convex and $\pedtheta{W}$-invariant subset 
$C  \subset \pedtheta{\frak a}$ with nonempty interior.
\end{Prop}
\pf 
The implication of $\supp D_q f \subset \exp C$ from 
$\supp f \subset \exp C$ is obvious. 
The converse inclusion follows from Theorem 1.2 in \cite{OP4}.
Indeed, the highest homogeneous part of $q(\l)$ is 
$$q_h(\l):=\prod_{\a \in \Sigma} \l_{w\a}^{\ma-1},$$
and Formula (\ref{eq:Dmeexplicit}) 
implies that the highest homogeneous part of $D(m;q)$ is 
\begin{equation*} 
  \prod_{\a \in \Sigma} \big(\partial(H_{w\a}/2)\big)^{\ma-1}.
\end{equation*}
The principal symbol of $D_q$ at $(H,\l) \in \frak a \times \frak a^*$ is therefore
\begin{equation}
  \label{eq:HHP}
  \sigma(D_q)(H,\l)=
 \Delta^k(H) \prod_{\a \in \Sigma} \l_{w\a}^{\ma-1}.
\end{equation}
\qed

\begin{Thm} \label{thm:onto}
Suppose $\Theta$ satisfies Condition ${\rm A}_2$. 
Let $C  \subset \pedtheta{\frak a}$ be an arbitrary compact, convex 
and $\pedtheta{W}$-invariant subset. Then $\supp \pedtheta{\mathcal I}g \subset C$
for all $g \in \pedtheta{\PW}(m;C)$.   
\end{Thm}
\pf 
Immediate consequence of Lemma \ref{lemma:keyequality} and  Proposition 
\ref{prop:key}.
\qed

\section{Conclusion of the proof of the Paley-Wiener theorem}
\label{section:endproof}

\noindent {\em Conclusion of the proof of Theorem  \ref{thm:pw}.\; }
Because of Theorems \ref{lemma:pwinto}
and  \ref{thm:onto}, we only need to prove that
$\pedtheta{\mathcal F}: C^\infty_c(C)^{\pedtheta{W}} \to \pedtheta{\PW}(m;C)$ and 
$\pedtheta{\mathcal I}: \pedtheta{\PW}(m;C) \to C^\infty_c(C)^{\pedtheta{W}}$ 
are (up to a fixed constant depending on $\Theta$ and $m$) inverse to each other. 

Let $g \in \pedtheta{\PW}(m;C)$. We have then proven that 
$\pedtheta{\mathcal F} \pedtheta{\mathcal I} g \in  \pedtheta{\PW}(m;C)$.
Observe that $\pedPi{\mathcal I}g=\mathcal I g$ is the $W$-invariant extension
of $\pedtheta{\mathcal I}g$ to $A$. 
Corollaries \ref{cor:functionaleqtr} and \ref{cor:pwopdam} together with Remark 
\ref{rem:invarianceI} therefore yield
\begin{equation*}
\big(\Pavtheta  \pedtheta{\mathcal F}\pedtheta{\mathcal I} g)(m;\l)=
(-1)^{d(\Theta,m)}  \big(\pedPi{\mathcal F}\pedPi{\mathcal I} g\big)(m;\l)=
\frac{1}{k(\Theta,m)} \Pavtheta g(m;\l),
\end{equation*}
where we have set
$$k(\Theta,m):=(-1)^{d(\Theta,m)}k\;  \frac{\abs{\pedtheta{W}}}{\abs{W}}.$$  
The injectivity of $\Pavtheta$ from Proposition \ref{prop:Pavthetainjective} gives then
\begin{equation*}
 \pedtheta{\mathcal F}\pedtheta{\mathcal I} g=\frac{1}{k(\Theta,m)}\; g
 \end{equation*}
for all $g \in \pedtheta{\PW}(m;C)$.
Conversely, the inversion formula of  Theorem \ref{thm:inversion} yields
\begin{equation*}
k(\Theta,m)\; \pedtheta{\mathcal I}\pedtheta{\mathcal F} f= f
 \end{equation*}
for $f \in C^\infty_c(C)^{\pedtheta{W}}$. 
This concludes the proof of Theorem \ref{thm:pw}.
\qed

\appendix\section{Estimates for the Harish-Chandra series}
\label{appendix:estPhi}

In this appendix we prove estimates for the derivatives 
of the Harish-Chandra series in the spectral parameter. 
As a corollary, they will provide a proof for Lemma \ref{lemma:partialPhi}.
In the case when no differentiation occurs, these estimates 
have been proven in \cite{P1}. The general case follows by
a straightforward modification of the arguments from 
the special case. 
We therefore only outline the proof and refer the reader to Section 4 
in \cite{P1} for the details.

In the following we adopt the notation of Sections \ref{section:prelim}
and \ref{section:HC}, but we allow $\Sigma$ to be an arbitrary
(not necessarily reduced) root system. We also assume that $\ma=0$
is equivalent to $\a\notin\Sigma$.
As in Section \ref{section:shift}, we denote by 
$(\l_1,\dots, \l_l)$ the complex coordinates in $\frakacs$ 
associated with an orthonormal  basis $(\xi_1,\dots,\xi_l)$ of $\frak a^*$. 
Moreover, if $I=(\iota_1,\dots,\iota_l) \in \N_0^\l $ is a 
multi-index and $H \in \liecomplex{a}$, we set 
\begin{equation*}
  \xi^I(H):=\big(\xi_1(H)\big)^{\iota_1} \dots \big(\xi_l(H)\big)^{\iota_l}.
\end{equation*}
The product rule for the derivation is then given by
\begin{equation*}
  \partial^I_\l(fg)=\sum_{J+K=I} \frac{I!}{J! K!} \, 
(\partial^J_\l f) (\partial^K_\l g).
\end{equation*}

\begin{Lemma} \label{lemma:multiindex}
Suppose $I=(\iota_1,\dots,\iota_l)$, $p \in \polya$, 
$m \in \N$,  $\l \in \frakacs$
and $a \in A$. 
\begin{enumerate}
\thmlist
\item \label{item:multiindexder}
$\partial^I_\l e^{\l(\log a)}=\xi^I(\log a) e^{\l(\log a)}$.
\item  \label{item:multindexholo}
Let $a$ be a meromorphic function on a domain  $D\subset \frakacs$. Suppose
$p(\l)a(\l)$ is holomorphic on $D$. Then also
\begin{equation*}
  p(\l)^{\abs{I}+1} \partial^I_\l a(\l)
\end{equation*}
is holomorphic in $D$.
\end{enumerate}
\nqed
\end{Lemma}

Estimates for the Harish-Chandra series and its derivatives will be obtained 
as in \cite{G1} and \cite{P1} using the modified Harish-Chandra series
\begin{equation*}
   \Psi(m;\l,a):=\Delta\left(m/2;a\right) \Phi(m;\l,a) 
\end{equation*}
with $\Delta(m)$ defined by (\ref{eq:Deltam}).
In particular, for every multi-index $I$, we have
\begin{equation}
  \label{eq:modHCseriesder}
  \partial_\l^I \Psi(m;\l,a):=\Delta\left(m/2;a\right) 
\partial_\l^I\Phi(m;\l,a). 
\end{equation}

For $S>0$, set $\frak a^+(S):=\{H \in \frak a: \text{$\a_j(H)>S$ 
for all $j=1,\dots, l$}\}$ and
$A^+(S):=\exp \frak a^+(S)$.
Suppose $\l \in \frakacs$ satisfies 
$\inner{\mu}{\mu-2\l}\neq 0$ for all $\mu \in  2\Lambda\setminus \{0\}$. 
Then, on $A^+$, we have the series expansions (\ref{eq:HCexp}) and
\begin{equation} \label{eq:seriesdelta}
   \Delta\left(m/2;a\right)=
e^{\rho(m)(\log a)} \prod_{\a\in \Sigma^+} 
\left(1-e^{-2\a(\log a)}\right)^{\ma/2}
=e^{\rho(m)(\log a)} \sum_{\mu \in 2\Lambda} b_\mu(m) e^{-\mu(\log a)}
\end{equation}
with $b_0(m)=1$. 
The Cauchy product of (\ref{eq:seriesdelta}) and  (\ref{eq:HCexp}) yields the
series expansion 
\begin{equation}\label{eq:seriesPsi} 
\Psi(m;\l,a)=e^{\l(\log a)} \sum_{\mu \in 2\Lambda} 
a_\mu(m;\l) e^{-\mu(\log a)},
\qquad a \in A^+,
\end{equation}
with
\begin{equation}\label{eq:amu}
a_\mu(m;\l):=\sum_{\stackrel{\nu,\eta \in 2\Lambda}{\nu+\eta=\mu}}  
b_\nu(m)\Gamma_\eta(m;\l)
\end{equation}
and $a_0(m;\l)=1$. The series converges absolutely in $A^+$ 
and uniformly in $A^+(S)$. 

The coefficients $a_\mu(m;\l)$ satisfy the recursion relations  
\begin{equation}\label{eq:recursiona}
a_\mu(m;\l) \inner{\mu-2\l}{\mu}=
\sum_{\a\in \Sigma^+} \ma(2-\ma-2\mduea)  \inner{\a}{\a}
\sum_{\stackrel{k\in \N}{\mu-2k\a\in 2\Lambda}} k a_{\mu-2k\a}(m;\l) 
\end{equation}
for $\mu \in 2\Lambda \setminus \{0\}$, with initial condition $a_0(m;\l)=1$.
See equation (37) in \cite{P1}.
It follows in particular that for $\mu \in 2\Lambda \setminus \{0\}$
each coefficient $a_\mu(m;\l)$ is a rational function of 
$\l \in \frakacs$ with at most simple poles 
along the hyperplanes
\begin{equation*}
  \mathcal H_\eta:=\{\l\in \frakacs:\inner{\eta-2\l}{\eta}=0\}
\end{equation*}
for some $\eta \in 2\Lambda \setminus \{0\}$ with $\eta \leq \mu$. 

In the following $R$ will always denote a \emph{finite} positive real number. 
We define
\begin{align}
  \frakacs(R)&:=\{\l \in \frakacs: \Re\inner{\l}{\a} < R 
                 \quad \text{for all $\a \in \Sigma^+$}\},
       \label{eq:frakacsR}\\
  \pedR{\mathcal X}&:=\{\eta \in 2\Lambda \setminus\{0\}: 
\mathcal H_\eta \cap \frakacs(2R)\neq 0\},
       \label{eq:Xr}\\
  \pedR{p}(\l)&:=\prod_{\eta \in \pedR{\mathcal X}} 
                          \inner{\eta-2\lambda}{\eta}.
      \label{eq:pR}
\end{align}
By Lemma 4.4 in \cite{P1}, 
the set $\pedR{\mathcal X}$ is finite, and there is $\omega>0$
such that  $\pedR{\mathcal X}=\emptyset$ for $0<R \leq \omega$. 
We convene that 
$\pedR{p}(\l)\equiv 1$ when $\pedR{\mathcal X}=\emptyset$ and
denote by $\deg \pedR{p}$ the degree of the polynomial $\pedR{p}$. 

\begin{Lemma} \label{lemma:series}
  Let $R,S>0$ be arbitrarily fixed. Then the series
  \begin{equation} \label{eq:series}
    \pedR{p}(\l) \sum_{\mu \in 2\Lambda} a_{\mu}(m;\l)e^{-\mu(\log a)}
  \end{equation}
converges uniformly in $(\l,a) \in\overline{V} 
\times A^+(S)T$, where $V$ is any open
subset with compact closure $\overline{V} \subset \frakacs(2R)$. 
Consequently, 
\begin{equation*}
  \pedR{p}(\l)\Psi(m;\l,a)= \pedR{p}(\l) e^{\l(\log a)} 
 \sum_{\mu \in 2\Lambda} a_{\mu}(m;\l)e^{-\mu(\log a)}
\end{equation*}
and, for every multi-index $I$,
\begin{equation*}
  \pedR{p}(\l)^{\abs{I}+1} \partial^I_\l\Psi(m;\l,a)
\end{equation*}
are holomorphic functions of $(\l,a) \in \frakacs(2R) \times A^+(S)U$, 
where $U\subset T$
is a neighborhood of $e$ on which $e^{\l(\log a)}$ is holomorphic 
for all $\l \in \frakacs$.
\end{Lemma}
\pf
The first part of the lemma depends on the fact that 
$\pedR{p}(\l)\Gamma_\mu(m;\l)$ 
is holomorphic in $\frakacs(2R)$ for all $\eta \in 2\Lambda$, together with 
(\ref{eq:amu}) and an easy modification of Lemma 2.1 in \cite{Opd88b}. 
The final
statement follows from  Lemma \ref{lemma:multiindex}(\ref{item:multindexholo}).
\qed

Because of Lemma \ref{lemma:series}, the series (\ref{eq:seriesPsi})  
converges uniformly on compact subsets of 
$\Big(\frakacs\setminus \big(\cup_{\mu \in 2\Lambda\setminus \{0\}} 
\mathcal H_\mu\big)\Big) \times
A^+U$. Termwise differentiation and 
Lemma \ref{lemma:multiindex}(\ref{item:multiindexder}) yield there
\begin{align}
  \partial^I_\l \Psi(m;\l,a)&= \sum_{\mu \in 2\Lambda} 
\partial_\l^I \Big(e^{\l(\log a)} a_\mu(m;\l)\Big) e^{\mu(\log a)} \notag \\
&=e^{\l(\log a)} \sum_{\mu \in 2\Lambda} \sum_{J+K=I} \frac{I!}{J!K!} \;
\xi^J(\log a)\, 
\partial_\l^K a_\mu(m;\l) e^{-\mu(\log a)}. \label{eq:PsiseriesI}
\end{align}
Finding estimates for $\partial^I_\l \Psi(m;\l,a)$ on $A^+$ 
is therefore reduced
to estimating the $\partial_\l^K a_\mu(m;\l)$.

\begin{Lemma} \label{lemma:polesamu}
Suppose $R>0$ and $\mu \in 2\Lambda$. Then the following holds
 for every multi-index $I$. 
\begin{enumerate}
\thmlist
\item
$\pedR{p}^{\abs{I}+1}(\l)\partial^I_\l a_\mu(m;\l)$ is holomorphic 
in $\l\in \frakacs(2R)$;
\item 
For all $\nu < \mu$, the function 
$ 
\; \displaystyle{\frac{\pedR{p}(\l)^{\abs{I}+1} 
\partial^I_\l a_\nu(m;\l)}{\inner{\mu-2\l}{\mu}}} \;
$
 is holomorphic in $\l \in \frakacs(2R)$.
\end{enumerate}
\end{Lemma}
\pf
This follows from Lemma 4.5 in \cite{P1} and  
Lemma \ref{lemma:multiindex}(\ref{item:multindexholo}).\qed

Let $e_j \in \N_0^l$ denote the multi-index with all coordinates 
zero except for the 
$j$-th equal to $1$. Differentiation of  (\ref{eq:recursiona})
yields for multi-indices $I=(\iota_1,\dots,\iota_l)\neq 0$  
the recursion relations
 \begin{multline}\label{eq:recursionaI}
\null \hskip -.3truecm \inner{\mu-2\l}{\mu} \partial_\l^I a_\mu(m;\l)
=2\sum_{\stackrel{j \in \{1,\dots,l\}}{I-e_j \in \N_0^l}} 
\inner{\xi_j}{\mu} \iota_j \partial^{I-e_j}_\l a_\mu(m;\l)\\ 
+\sum_{\a\in \Sigma^+} \ma(2-\ma-2\mduea)  \inner{\a}{\a}\!\!
\sum_{\stackrel{k\in \N}{\mu-2k\a\in 2\Lambda}} k 
\partial_\l^I a_{\mu-2k\a}(m;\l) 
\end{multline} 
for $\mu \in 2\Lambda \setminus\{0\}$, with initial condition 
 $\partial_\l^I a_0(m;\l)=0$.

The procedure for estimating the derivatives 
$\partial_\l^Ia_\mu(m;\l)$ is based, as
in the case $I=0$ of \cite{P1}, on the comparison of the recursion relations
 (\ref{eq:recursionaI}) with the recursion relations satisfied by the 
coefficients $d_\mu(m;c)$ of the series expansion of the function
\begin{equation}
  \label{eq:deltac}
  \Delta_c(m):=e^{c\rho(m)} \Delta(-cm/2)=
\prod_{\a\in \Sigma^+} \left(1-e^{-2\a}\right)^{-c\ma/2},
\end{equation}
where $c\in [0,1)$ will be suitably chosen.

\begin{Lemma} \label{lemma:deltac}
Suppose $m \in \mathcal M^+$ and $c \in \R$.
\begin{enumerate}
\thmlist
\item {\rm (\cite{P1}, Lemma 4.8)}
The function $\Delta_c(m)$ defined in (\ref{eq:deltac})
admits the series expansion 
\begin{equation} \label{eq:seriesdeltac}
\Delta_c(m)=\sum_{\mu \in 2\Lambda} d_\mu(m;c) e^{-\mu}
\end{equation}
which converges absolutely in $A^+$ and uniformly in $A^+(S)$ for all $S>0$.
We have $d_0(m;c)=1$, and,
if we assume $c\in (0,\infty)$, then $d_\mu(m;c) > 0$ for all 
$\mu \in 2\Lambda$. 
\item {\rm (\cite{P1}, Lemma 4.10)}
Suppose $H \in \frak a$ and $c \in [0,\infty)$.
Then the  coefficients $d_\mu(m;c)$ of the series (\ref{eq:seriesdeltac}) 
satisfy the recurrence relations
\begin{multline} \label{eq:recursiond}
  \big(\inner{\mu}{\mu}+\mu(H)\big)d_\mu(m;c)= 
\sum_{\a \in \Sigma^+} \sum_{\stackrel{k\in\N}{\mu-2k\a \in2\Lambda}} 
\left[ 2c\ma \left(\frac{\a(H)}{2}-c\inner{\rho(m)}{\a}\right) +\right.\\
+ \Big.k c\ma(c\ma+2c\mduea +2) \inner{\a}{\a} \Big] 
d_{\mu-2k\a}(m;c)
\end{multline}
for $\mu \in 2\Lambda \setminus \{0\}$, and $d_0(m;c)=1$.
\nqed
\end{enumerate}
\end{Lemma}

The constant $c$ is chosen according to the following lemma.

\begin{Lemma} 
\label{lemma:cr}
{\rm (\cite{P1}, Lemma 4.11)}
There exist constants $0\leq c<1$ and $r>1$ such that the inequality
  \begin{equation*}
    c(c\ma+2c\mduea+2) \geq r\abs{\ma+2\mduea-2}
  \end{equation*}
 holds for all $\a \in \Sigma^+$.  
\nqed
\end{Lemma}

\begin{Rem}
Assume, as above, that $\ma>0$ for $\a in\Sigma$. Then
the system of inequalities in Lemma \ref{lemma:cr} 
admits the solution $c=0$ if and only if $\Sigma$ is reduced and 
$\ma=2$ for all $\a \in \Sigma^+$.
\end{Rem}

Let $c \in [0,1)$ obtained according to Lemma \ref{lemma:cr} for some
$r>1$. The comparison of the recursion relations (\ref{eq:recursionaI})
and (\ref{eq:recursiond}) is possible for a choice of $H \in \frak a$ so 
that $\a(H) \geq \max\{2c\inner{\rho(m)}{\a},0\}$ for all $\a \in \Sigma^+$.
This condition ensures that $\mu(H)\geq 0$ for all $\mu \in 2\Lambda$
and that the first summands on the right-hand side of (\ref{eq:recursiond})
are nonnegative.

Arguments similar to those in Lemmas 4.7, 4.14 and 4.16 in \cite{P1}, 
applied to the recursion relations 
(\ref{eq:recursionaI}) inductively on $\abs{I}$, lead to the 
required estimates 
 for the derivatives $\partial^I_\l a_\mu(m;\mu)$ and hence to 
estimates for $\partial^I_\l \Psi(m;\l,a)$ and $\partial^I_\l \Phi(m;\l,a)$.

\begin{Lemma}\label{lemma:estaI}
Suppose $m \in \mathcal M^+$.
Let $R>0$ and let $\pedR{p}$ be the polynomial defined in (\ref{eq:pR}). 
Suppose $c\in[0,1)$
is chosen as in Lemma \ref{lemma:cr} for some $r>1$. 
Then, for every multi-index $I$, 
there is a constant $K_{R,c,m,I}>0$ such that 
 \begin{equation*}
\abs{\pedR{p}(\l)^{\abs{I}+1} \partial_\l^I a_\mu(m;\l)} 
\leq K_{R,c,m,I} d_\mu(m;c) (1+\abs{\l})^{(\abs{I}+1)\deg \pedR{p}}  
 \end{equation*}
for all $\l \in \frakacs(R)$ and all $\mu \in 2\Lambda$. 
\nqed
\end{Lemma}
 
\begin{Thm}\label{thm:estPhiI}
Let $m \in \mathcal M^+$, $R>0$ and let $I$ be a multi-index. 
Let $\pedR{p}$ denote the polynomial 
defined in (\ref{eq:pR}). 
Then there exist $c \in [0,1)$ 
(depending only on the multiplicity function $m$) and  $C_{R,c,m,I}>0$
such that 
\begin{equation}
  \label{eq:estPsi}
 \abs{\pedR{p}(\l)^{(\abs{I}+1)} \partial_\l^I\Psi(m;\l,a)} 
\leq  C_{R,c,m,I} \Delta(m;a)^{-c/2} (1+\abs{\log a})^{\abs{I}}
(1+\abs{\l})^{(\abs{I}+1)\deg \pedR{p}}
e^{(c\rho(m)+\Re\l)(\log a)} 
\end{equation}
and
 \begin{equation}
  \label{eq:estPhiI}
\abs{\pedR{p}(\l)^{(\abs{I}+1)}  \Delta(m;a)^{(c+1)/2} \partial_\l^I
\Phi(m;\l,a)} \leq  
C_{R,c,m,I}  (1+\abs{\log a})^{\abs{I}}(1+\abs{\l})^{(\abs{I}+1)\deg \pedR{p}}
e^{(c\rho(m)+\Re\l)(\log a)}
\end{equation}
for all $a \in A^+$ and $\l \in \frakacs(R)$.

There is $\omega >0$ so that we can choose $\pedR{p}\equiv 1$ for 
 $0<R \leq \omega$. In this case,
 \begin{equation}
  \label{eq:estPhiIRsmall}
\abs{\Delta(m;a)^{(c+1)/2} \partial_\l^I
\Phi(m;\l,a)} \leq  
C_{R,c,m,I}  (1+\abs{\log a})^{\abs{I}}e^{(c\rho(m)+\Re\l)(\log a)}
\end{equation}
for all $a \in A^+$ and $\l \in \frakacs(R)$.
\end{Thm}
\pf
Notice first that for $J \leq I$
\begin{equation*}
  \abs{\xi^J(\log a)} \leq \abs{\log a}^{\abs{J}} \leq
  (1+\abs{\log a})^{\abs{I}}.
\end{equation*}
Hence, because of (\ref{eq:PsiseriesI}), (\ref{eq:seriesdeltac})
 and Lemma \ref{lemma:estaI},
\begin{align*}
\abs{\pedR{p}(\l)^{\abs{I}+1} \partial^I_\l \Psi(m;\l,a)}&= 
e^{\Re\l(\log a)} \sum_{\mu \in 2\Lambda} \sum_{J+K=I} \frac{I!}{J!K!} \;
\abs{\xi^J(\log a)} 
\abs{\pedR{p}(\l)^{\abs{I}+1}\partial_\l^K a_\mu(m;\l)} e^{-\mu(\log a)}\\
&\leq C_{R,c,m,I} \sum_{\mu \in 2\Lambda} d_\mu(m;c) e^{-\mu(\log a)} 
(1+\abs{\log a})^{\abs{I}}(1+\abs{\l})^{(\abs{I}+1)\deg \pedR{p}}  
e^{\Re\l(\log a)}\\
&\leq C_{R,c,m,I} 
\Delta_c(m;a)(1+\abs{\log a})^{\abs{I}}
(1+\abs{\l})^{(\abs{I}+1)\deg \pedR{p}}   e^{\Re\l(\log a)}.
\end{align*}
The first inequality therefore follows from (\ref{eq:deltac}).
The last statement is a consequence of (\ref{eq:modHCseriesder}) and  
of the remark after (\ref{eq:pR}).
\qed

Lemma \ref{lemma:partialPhi} is a special case of the following
corollary to Theorem \ref{thm:estPhiI}. 

\begin{Cor} \label{cor:estonA}
Let $m \in \mathcal M^+$ and let $\omega$ be as in Theorem \ref{thm:estPhiI}.
For every fixed $\l \in \frakacs(\omega)$ and every multi-index 
$I$ the function 
$\Delta(m;a) \partial_\l^I \Phi(m;\l,a)$ extends continuously 
on $\overline{A^+}$ 
by setting it equal to zero on the boundary $\partial(A^+)$ of $\overline{A^+}$.
\nqed
\end{Cor}

\section{Proof of Corollary \ref{cor:pieceofDm}}
\label{app:proofcor}

The ring $\polyAc$ is an integral domain and a UFD 
(see e.g. \cite{Bou}, Ch. VI, \S3, Lemma 1(i)). 
This in particular implies that irreducible and prime elements in $\polyAc$ coincide
and that every finite set of elements of $\polyAc$ possesses a greatest common divisor.

\begin{Lemma} 
  \label{lemma:relprimebis}
  Let $\Sigma$ be a reduced root system in $\frak a^*$, and let $\Delta=\Delta(1)$ denote the
   Weyl denominator.
  \begin{enumerate}
  \thmlist
  \item 
For every $\a \in \Sigma^+$ the elements $1+e^{-2\a}$ and $1-e^{-2\a}$ are relatively
prime in $\polyAc$.
   \item
$\prod_{\a\in \Sigma^+} \partial(A_\a)(\Delta) \in \polyAc$, and 
$\prod_{\a\in \Sigma^+} \partial(A_\a)(\Delta)$ and $1-e^{-2\beta}$ are relatively prime 
for all $\beta \in \Sigma^+$.
  \end{enumerate}
\end{Lemma}
\pf
Suppose $u \in \polyAc$ divides $1+e^{-2\a}$ and $1-e^{-2\a}$. Then there are 
$f,g \in \polyAc$ so that $1+e^{-2\a}=fu$ and $1-e^{-2\a}=gu$. Thus
$2=(1+e^{-2\a})+(1-e^{-2\a})=(f+g)u$, which implies that $u$ is a unit in $\polyAc$.
This proves (a). 

For (b), let $\beta_1,\dots,\beta_n$ be an enumeration of $\Sigma^+$, and set
$A_j=A_{\b_j}$ for $j=1,\dots,n$. Then
\begin{equation*}
  \partial(A_j)(\Delta)= e^{\rho(2)} \; 
        \sum_{\a\in \Sigma^+} \inner{\b_j}{\a} (1+e^{-2\a}) 
   \prod_{\gamma \in \Sigma^+ \setminus \{\a\}} (1-e^{-2\gamma}) \in \polyAc. 
\end{equation*}
Hence $\prod_{\a\in \Sigma^+} \partial(A_\a)(\Delta) \in \polyAc$, and more precisely, 
we have
\begin{align*}
  \prod_{\a\in \Sigma^+} \partial(A_\a)(\Delta) &=\prod_{j=1}^n \partial(A_j)(\Delta)\\
    &=e^{n\rho(2)} \prod_{j=1}^n \sum_{k_j=1}^n \inner{\b_j}{\b_{k_j}} 
        (1+e^{-2\b_{k_j}}) \prod_{\gamma_{k_j} \in \Sigma^+ \setminus \{\b_{k_j}\}} 
       (1-e^{-2\gamma_{k_j}})\\
    &=e^{n\rho(2)} \sum_{k_1=1}^n \cdots \sum_{k_n=1}^n \;\;
      \prod_{j=1}^n \Big[\inner{\b_j}{\b_{k_j}}    (1+e^{-2\b_{k_j}})
      \prod_{\gamma_{k_j} \in \Sigma^+ \setminus \{\b_{k_j}\}} (1-e^{-2\gamma_{k_j}}) \Big].
 \end{align*}
Let $\b\in \Sigma^+$ be fixed. The only summand in which  $1-e^{-2\b}$ does not
appear as a factor corresponds to $\b_{k_1}=\dots=\b_{k_n}=\beta$. 
Collecting the remaining terms together, we get for some $f_\b \in \polyAc$
\begin{equation*}
 \prod_{\a\in \Sigma^+} \partial(A_\a)(\Delta)=(1-e^{-2\b})f_\beta+ 
 e^{n\rho(2)}  \Big( \prod_{j=1}^n \inner{\b_j}{\b} \Big) (1+e^{-2\b})^n 
    \prod_{\gamma\in \Sigma^+ \setminus \{\b\}} (1-e^{-2\gamma})^n.
\end{equation*}
By Lemma \ref{lemma:relprime} and Part (a), the elements $1-e^{-2\b}$ and 
$ (1+e^{-2\b})^n 
    \prod_{\gamma\in \Sigma^+ \setminus \{\b\}} (1-e^{-2\gamma})^n$ 
are relatively prime. Thus also $1-e^{-2\b}$ and 
$\prod_{\a\in \Sigma^+} \partial(A_\a)(\Delta)$ must be relatively prime. 
\qed

\begin{Lemma}
  \label{lemma:divisionDelta}
Let $\omega \in \polyAc$, and suppose $\Delta$ does not divide $\omega$  in $\polyAc$. 
Then $\Delta$ does not divide $ \omega \; \prod_{\a\in \Sigma^+} \partial(A_\a)(\Delta) $ in 
 $\polyAc$. 
\end{Lemma}
\pf
We prove that, if  $\Delta$ divides $\omega \; \prod_{\a\in \Sigma^+} \partial(A_\a)(\Delta)$ in 
 $\polyAc$, then $\Delta$ must divide $\omega$. 
Let $u\in\polyAc$ be a prime dividing $\Delta=e^{\rho(2)} \prod_{\b \in \Sigma^+} (1-e^{-2\b})$. 
Then $u$ divides $1-e^{-2\b}$ for some
$\b \in \Sigma^+$. Since  $1-e^{-2\b}$ and 
$\prod_{\a\in \Sigma^+} \partial(A_\a)(\Delta)$ are relatively prime, so are also 
$u$ and $\prod_{\a\in \Sigma^+} \partial(A_\a)(\Delta)$. But we are assuming that
$\Delta$, hence $u$, divides $\omega \; \prod_{\a\in \Sigma^+} \partial(A_\a)(\Delta)$.
So $u$ divides $\omega$. Since $\polyAc$ is a UFD and $u$ is an arbitrary prime
dividing $\Delta$, we conclude that $\Delta$ divides $\omega$. 
\qed

\begin{Lemma}
  \label{lemma:singularcoeff}
Let $\Delta^k \sum_I \omega_I \otimes \partial_a^I$ be the representation of an 
element of $\polyAc \otimes \polya$ with $\omega_I \in \polyAc$ and $k \in \Z$ maximal.
If $k \in -\N$, then the differential operator 
\begin{equation*}
 \Big( \prod_{\a\in \Sigma^+} \partial(A_\a)(\Delta)\Big) \circ 
 \Big( \Delta^k \sum_I \omega_I \otimes \partial_a^I \Big)
\end{equation*}
has singular coefficients on $\complex{A}$.
\end{Lemma}
\pf
As in the proof of Lemma \ref{lemma:relprimebis} (b), let us write 
$\prod_{\a\in \Sigma^+} \partial(A_\a)=\prod_{j=1}^n \partial(A_j)$. 
We first prove by induction on $h\leq n$ that 
\begin{multline}
  \label{eq:composition}
\prod_{j=1}^h \partial(A_j) \circ \Delta^k \omega_I =
 k(k-1)\cdots (k-h+1) \; \Delta^{k-h} \prod_{j=1}^h \partial(A_j)(\Delta) \; \omega_I \\
     + \Delta^{k-h+1} \Big[\omega_{I,h} + \sum_{J(I,h)} \omega_{J(I,h)} \otimes 
\partial(A_{J(I,h)}) \Big]
\end{multline}
with $\omega_{I,h},  \omega_{J(I,h)} \in \polyAc$, and where
$\partial(A_{J(I,h)}):=\partial(A_1)^{j_1}\cdots \partial(A_h)^{j_h}$ if 
$J(I,h)=(j_1,\dots,j_h)$.
Indeed, for $h=1$ we have
\begin{equation*}
  \partial(A_1) \circ \Delta^k \omega_I=k\Delta^{k-1} \partial(A_1)(\Delta)\omega_I
+ \Delta^k\big[ \partial(A_1)(\omega_I)+ \omega_I \otimes \partial(A_1)\big].
\end{equation*}
Suppose inductively that (\ref{eq:composition}) holds for $h$. Then
\begin{align*}
\prod_{j=1}^{h+1} \partial(A_j) \circ \Delta^k \omega_I &=
  \partial(A_{h+1}) 
\Big[k(k-1)\cdots (k-h+1) \Delta^{k-h} \prod_{j=1}^h \partial(A_j)(\Delta) \omega_I\Big]\\
& \qquad + \partial(A_{h+1}) \Big[ \Delta^{k-h+1} \omega_{I,h}
+ \Delta^{k-h+1} \sum_{J(I,h)} \omega_{J(I,h)} \otimes \partial(A_{J(I,h)}) \Big]
\displaybreak[0] \\
&=k(k-1)\cdots (k-h+1)(k-h) \; \Delta^{k-h-1} 
 \prod_{j=1}^{h+1} \partial(A_j)(\Delta) \; \omega_I \\
& \qquad + k(k-1)\cdots (k-h+1)\; \Delta^{k-h} \partial(A_{h+1}) 
        \Big[\prod_{j=1}^h \partial(A_{h+1})(\Delta) \;\omega_I\Big] \\
&\qquad + (k-h+1)\; \Delta^{k-h}\;\partial(A_{h+1})(\Delta)\;\omega_{I,h}+  
                        \Delta^{k-h+1} \;\partial(A_{h+1})(\omega_{I,h}) \\
&\qquad + (k-h+1)\; \Delta^{k-h} \;\partial(A_{h+1})(\Delta) 
          \sum_{J(I,h)} \omega_{J(I,h)} \otimes \partial(A_{J(I,h)})\\
& \qquad  +\Delta^{k-h+1}\sum_{J(I,h)} \partial(A_{h+1})(\Delta) \;
                              \omega_{J(I,h)}\otimes \partial(A_{J(I,h)})\\
& \qquad   + \Delta^{k-h+1} \sum_{J(I,h)} \omega_{J(I,h)} 
                                    \otimes \partial(A_{h+1}) \partial(A_{J(I,h)})\\
&= k(k-1)\cdots (k-h) \; \Delta^{k-h-1}  
        \prod_{j=1}^{h+1} \partial(A_j)(\Delta) \; \omega_I \\
& \qquad
   +  \Delta^{k-h}\Big[ \omega_{I,h+1} +  
       \sum_{J(I,h+1)} \omega_{J(I,h+1)} \otimes \partial(A_{J(I,h+1)}) \Big],
\end{align*}
where
\begin{multline*}
 \omega_{I,h+1}:= k(k-1)\cdots (k-h+1)\;\Delta^{k-h}\; \partial(A_{h+1}) 
        \Big[\prod_{j=1}^h \partial(A_{h+1})(\Delta) \omega_I\Big] \\
+ (k-h+1)\Delta^{k-h}\;\partial(A_{h+1})(\Delta)\;\omega_{I,h}+  
                        \Delta^{k-h+1}\; \partial(A_{h+1})(\omega_{I,h})  \in \polyAc
\end{multline*}
and
\begin{multline*}
 \sum_{J(I,h+1)} \omega_{J(I,h+1)} \otimes \partial(A_{J(I,h+1)}):=
  (k-h+1) \partial(A_{h+1})(\Delta) 
          \sum_{J(I,h)} \omega_{J(I,h)} \otimes \partial(A_{J(I,h)})\\
+ \Delta \sum_{J(I,h)} \partial(A_{h+1})(\Delta)\omega_{J(I,h)}\otimes \partial(A_{J(I,h)})\\
  + \Delta \sum_{J(I,h)} \omega_{J(I,h)} \otimes \partial(A_{h+1}) \partial(A_{J(I,h)}) 
\in \polyAc \otimes \polya.
\end{multline*}
Since $k$ is maximal, there is a multiindex $I_0$ such that $\Delta$ does not divide
$\omega_{I_0}$. 
 Equation (\ref{eq:composition}) for $h=n$ shows that the coefficient of $\partial_a^{I_0}$ in 
$ \prod_{\a\in \Sigma^+} \partial(A_\a)(\Delta) \circ 
\Delta^k \sum_I \omega_I \otimes \partial_a^I$ is 
equal to the sum of 
\begin{equation*}
 k(k-1)\cdots (k-n+1) \Delta^{k-n} \prod_{\a\in\Sigma^+} 
        \partial(A_\a)(\Delta) \; \omega_{I_0}
     + \Delta^{k-n+1} \omega_{I_0,n}
 \end{equation*}
and of those coefficients of 
$\Delta^{k-n+1} \sum_{I\neq I_0} \sum_{J(I,n)} \omega_{J(I,n)} \otimes 
\partial(A_{J(I,n)})\partial_a^I$
coming from the terms in $\partial(A_{J(I,n)})\partial_a^I$ which are equal to $\partial_a^{I_0}$.
Consequently, the coefficient of $\partial_a^{I_0}$ is of the form
\begin{equation*}
  k(k-1)\cdots (k-n+1) \Delta^{k-n} \prod_{\a\in\Sigma^+} 
        \partial(A_\a)(\Delta)\; \omega_{I_0} +  \Delta^{k-n+1} \sigma_0
\end{equation*}
with $\sigma_0 \in \polyAc$. By Lemma \ref{lemma:divisionDelta}, $\Delta$ does not divide
$\prod_{\a\in \Sigma^+} \partial(A_\a) \omega_{I_0}$. 
Thus the singularity of $\Delta^{k-n}$ (with $k-n<0$) is neither cancelled by 
$\prod_{\a\in \Sigma^+} \partial(A_\a) \; \omega_{I_0}$ nor by $\Delta^{k-n+1} \sigma_0$
(in which a singularity of lower order appears).
\qed 

\noindent {\em Proof of Corollary \ref{cor:pieceofDm}.\;}
Set $G:=\Delta(m)G_+(m-2;2)\circ \Delta(-1)$. Using the definition of $D_m$ in (\ref{eq:Dm})
together with (\ref{eq:DplusHC}) and Example \ref{ex:shiftcomplex}, we obtain
\begin{align*}
\Delta(m)D_m&=\Delta(m)D_+(m)\\
                   &=\Delta(m)G_+(m-2;2)\circ D_+(2)\\
  &=\Delta(m)G_+(m-2;2)\circ 
                   \widetilde\sigma \Delta^{-1} \prod_{\a\in\Sigma^+} \partial(A_\a)\\
  &=\widetilde{\sigma} \Big( G \circ \prod_{\a\in \Sigma^+} \partial(A_\a)\Big).
\end{align*}
Hence 
\begin{align*}
  \big(\Delta(m)D_m\big)^*&=\widetilde\sigma (-1)^{\abs{\Sigma^+}} 
                           \prod_{\a\in \Sigma^+} \partial(A_\a) \circ G^*\\
        &=\sigma  \prod_{\a\in \Sigma^+} \partial(A_\a) \circ G^*.
\end{align*}
By Theorem \ref{thm:Dm}, $\Delta(m)D_m$, and hence $\big(\Delta(m)D_m\big)^*$, belongs to
$\polyAc\otimes \polya$. Since $G \in \C_\Delta[\complex{A}]\otimes \polya$, so does
$G^*$. Thus $G^*=\Delta^k \sum_I \omega_I \otimes \partial_a^I$ with 
$\omega_I \in \polyAc$ and $k \in \Z$ maximal. Since $k \in -\N$ 
would imply $\big(\Delta(m)D_m\big)^* \notin \polyAc\otimes \polya$ by 
Lemma \ref{lemma:singularcoeff}, we conclude
that $k \in \N_0$, i.e. $G^* \in  \polyAc\otimes \polya$. Thus $G \in \polyAc\otimes \polya$. 
\qed

\section{$K_\e$-symmetric spaces with even multiplicities}
\label{app:list}

In this appendix we report the infinitesimal classification of 
$K_\e$-symmetric spaces with even multiplicities by listing
the $K_\e$-symmetric pairs $(\frak g,\frak h)$ 
with even multiplicities for which $\frak g$ is simple
and noncompact. The list has been extracted 
from the classification due to Oshima and Sekiguchi \cite{OS}.
It is presented in three tables respectively collecting 
(for the even multiplicity case) the Riemannian symmetric pairs (Table 1),
the non-compactly causal (NCC) symmetric pairs (Table 2) and the other 
$K_\e$ symmetric pairs (Table 3). A non-Riemannian $K_\e$-symmetric pair
is said to be of type $K_\e I$ if its signature $\e$ comes from a gradation
of first kind according to \cite{Kane}. Otherwise it is said to be 
of type $K_\e II$. The symmetric pairs of type $K_\e I$ coincide with the 
NCC symmetric pairs. Table 3 therefore collects all symmetric pairs 
with even multiplicities of type $K_\e II$. 

The restricted root system $\Sigma$ of a $K_\e$-symmetric pair 
with even multiplicities has at most two root lengths. The classification
below shows that all multiplicities $m_\a$ of $\Sigma$ are equal, 
and moreover that they are all equal to $2$ for symmetric pairs of type
$K_\e II$.   
The restricted root system and multiplicities of a $K_\e$-symmetric pair 
$(\frak g, \frak h)$ coincide with
those of the corresponding Riemannian dual symmetric pair 
$(\frak g,\frak k)$. They are explicitly reported in Tables 
2 and 3 for the reader's convenience.

If $\Sigma$ is of type $X_n$ (with $X_n \in \{A_n,B_n,C_n,\dots\}$), then 
the index $n$ denotes the real rank of $\frak g$. The range for $n$ 
is chosen to avoid overlappings due to isomorphisms of symmetric spaces.
These isomorphisms arise from isomorphisms of the lower dimensional 
complex Lie algebras. We refer to \cite{He1}, Ch. X, \S 6, for more 
information. After each table below we report the relevant
symmetric pair isomorphisms.

\begin{table}[h]
\label{table:Riem}
\setlength{\extrarowheight}{4pt}
{\begin{tabular}{|>{$}c<{$}|>{$}c<{$}|>{$}c<{$}|>{$}c<{$}|>{$}c<{$}|>{$}c<{$}|}
\hline
\frak{g} &\frak{h}=\frak{k}&\Sigma&m_\alpha
\rule[-0.2cm]{0cm}{0.0cm}
&\\
\hline
\hline
\frak{sl}(n,\mathbb{C})& \frak{su}(n)
&A_{n-1}&2&n\geq 2
\rule[-0.2cm]{0cm}{0.0cm}
\cr\hline
\frak{so}(2n+1,\mathbb{C})&\frak{so}(2n+1)
&B_n&2&n\geq 2
\rule[-0.2cm]{0cm}{0.0cm}
\cr\hline
\frak{sp}(n,\mathbb{C})&\frak{sp}(n)
&C_n&2&n\geq 3
\rule[-0.2cm]{0cm}{0.0cm}
\cr\hline
\frak{so}(2n,\mathbb{C})&\frak{so}(2n)
&D_n&2&n\geq 4
\rule[-0.2cm]{0cm}{0.0cm}
\cr\hline
(\frak{e}_6)_{\smC}&\frak{e}_6&E_6&2
\rule[-0.2cm]{0cm}{0.0cm}
&\cr\hline
(\frak{e}_7)_{\smC}&\frak{e}_7&E_7&2
\rule[-0.2cm]{0cm}{0.0cm}
&\cr\hline
(\frak{e}_8)_{\smC}&\frak{e}_8&E_8&2
\rule[-0.2cm]{0cm}{0.0cm}
&\cr\hline
(\frak{f}_4)_{\smC}&\frak{f}_4&F_4&2
\rule[-0.2cm]{0cm}{0.0cm}
&\cr\hline
(\frak{g}_2)_{\smC}&\frak{g}_2&G_2&2
\rule[-0.2cm]{0cm}{0.0cm}
&\cr\hline
\frak{su}^*(2n)&\frak{sp}(n)&A_{n-1}&4&n\geq 2
\rule[-0.2cm]{0cm}{0.0cm}
\cr\hline
\frak{e}_{6(-26)}&\frak{f}_{4(-20)}&A_2&8
\rule[-0.2cm]{0cm}{0.0cm}
&\cr\hline
\frak{so}(2n+1,1)&\frak{so}(2n+1)&A_1&2n&n\geq 3
\rule[-0.2cm]{0cm}{0.0cm}
\cr\hline
\end{tabular}
\bigskip
\caption{Riemannian symmetric pairs with even multiplicities.} }
\end{table}

Special isomorphisms of Riemannian symmetric spaces with even multiplicities
are
\begin{alignat*}{2}
 \frak{so}(3,\C)=\frak{sp}(1,\C) &\approx \frak{sl}(2,\C), 
  &\quad
 \frak{so}(3)=\frak{sp}(1) &\approx \frak{su}(2)\,;\\
 \frak{sp}(2,\C) &\approx \frak{so}(5,\C),
  &\quad
\frak{sp}(2) &\approx \frak{so}(5)\,;\\
 \frak{so}(6,\C) &\approx \frak{sl}(4,\C),
  &\quad
\frak{so}(6) &\approx \frak{su}(4)\,;\\
 \frak{so}(3,1) &\approx \frak{sl}(2,\C),
  &\quad
\frak{so}(3) &\approx \frak{su}(2)\,;\\
 \frak{so}(5,1) &\approx \frak{su}^*(4),
  &\quad
\frak{so}(5) &\approx \frak{sp}(2)\,.
\end{alignat*}

The Lie algebra $\frak{so}(2,\C)$ is not semisimple.  
Observe  also that 
$\frak{so}(4,\C)\cong \frak{sl}(2,\C) \times \frak{sl}(2,\C)$ is not simple.
The structure of its homogeneous spaces can be therefore deduced from 
the structure of the homogeneous spaces of $\frak{sl}(2,\C)$.
 
\bigskip 
In the next table we list all the non-compactly causal, or $K_\e I$, 
symmetric pairs with even multiplicities. The third column reports the
subalgebra of $\frak g$ fixed by $\theta\theta_\e$, where $\theta_e$ is the 
involution associated with the $K_\e$-pair $(\frak g,\frak h)$.

\begin{table}[h]\label{table:NCC}
\setlength{\extrarowheight}{4pt}
\begin{tabular}{|>{$}c<{$}|>{$}c<{$}|>{$}c<{$}|>{$}c<{$}|>{$}c<{$}|>{$}c<{$}|}
\hline
\frak{g} &\frak{h}& \frak{g}^{\theta\theta_\e}&\Sigma& m_\alpha 
\rule[-0.2cm]{0cm}{0.0cm}
& 
 \\
\hline
\hline
\frak{sl}(n,\mathbb{C})&\frak{su}(n-j,j)
&\frak{sl}(n-j,\mathbb{C})\times\frak{sl}(j,\mathbb{C})\times\mathbb{C}
&A_{n-1}&2
&n\geq 2,\, 1\leq j\leq [n/2]
\rule[-0.2cm]{0cm}{0.0cm}
\cr\hline
\frak{so}(2n+1,\mathbb{C})&\frak{so}(2n-1,2)
&\frak{so}(2n-1,\mathbb{C})\times \mathbb{C}&B_n&2
&n\geq 2
\rule[-0.2cm]{0cm}{0.0cm}
\cr\hline
\frak{sp}(n,\mathbb{C})&\frak{sp}(n,\mathbb{R})
&\frak{gl}(n,\mathbb{C})&C_n&2
&n\geq 3
\rule[-0.2cm]{0cm}{0.0cm}
\cr\hline
\frak{so}(2n,\mathbb{C})&\frak{so}(2n-2,2)
&\frak{so}(2n-2,\mathbb{C})\times \mathbb{C}&D_n&2
&n\geq 4
\rule[-0.2cm]{0cm}{0.0cm}
\cr\hline
\frak{so}(2n,\mathbb{C})&\frak{so}^*(2n)
&\frak{gl}(n,\mathbb{C})&D_n&2
&n\geq 4
\rule[-0.2cm]{0cm}{0.0cm}
\cr\hline
(\frak{e}_6)_{\smC}&\frak{e}_{6(-14)}
&\frak{so}(10,\mathbb{C})\times \mathbb{C}&E_6&2
\rule[-0.2cm]{0cm}{0.0cm}
& 
\cr\hline
(\frak{e}_7)_{\smC}&\frak{e}_{7(-25)}
&(\frak{e}_6)_{\smC}\times\mathbb{C}&E_7&2
\rule[-0.2cm]{0cm}{0.0cm}
& 
\cr\hline
\frak{su}^*(2n)&\frak{sp}(n-j,j)
&\frak{su}^*(2(n-j))\times \frak{su}^*(2j)\times\mathbb{R}&A_{n-1}&4
&n\geq 2,\, 1\leq j \leq [n/2]
\rule[-0.2cm]{0cm}{0.0cm}
\cr\hline
\frak{e}_{6(-26)}&\frak{f}_{4(-20)}
&\frak{so}(9,1)\times \mathbb{R}&A_2&8&\cr\hline
\frak{so}(2n+1,1)&\frak{so}(2n,1)
&\frak{so}(2n+1)\times \mathbb{R}&A_1&2n
&n\geq 3
\rule[-0.2cm]{0cm}{0.0cm}
\cr\hline
\end{tabular}
\bigskip
\caption{Non-compactly causal symmetric pairs with even multiplicities.}  
\end{table}

\addtolength{\textheight}{40pt}

Special isomorphisms of NCC symmetric pairs with even multiplicities are
\begin{alignat*}{2}
 \frak{so}(3,\C)=\frak{sp}(1,\C) &\approx \frak{sl}(2,\C), 
  &\quad
 \frak{so}(1,2)\approx\frak{sp}(1,\R) &\approx \frak{su}(1,1)\,;\\
 \frak{sp}(2,\C) &\approx \frak{so}(5,\C),
  &\quad
\frak{sp}(2,\R) &\approx \frak{so}(3,2)\,;\\
 \frak{so}(6,\C) &\approx \frak{sl}(4,\C),
  &\quad
\frak{so}(4,2) &\approx \frak{su}(2,2)\,;\\
\frak{so}(6,\C) &\approx \frak{sl}(4,\C),
  &\quad
\frak{so}^*(6) &\approx \frak{su}(3,1)\,;\\
 \frak{so}(3,1) &\approx \frak{sl}(2,\C),
  &\quad
\frak{so}(2,1) &\approx \frak{su}(1,1)\,;\\
 \frak{so}(5,1) &\approx \frak{su}^*(4),
  &\quad
\frak{so}(4,1) &\approx \frak{sp}(1,1)\,.
\end{alignat*}
 
\smallskip 
The last table contains all the other $K_\e$-symmetric pairs, i.e. those of
type $K_\e II$, with even multiplicities.

\begin{table}[h] \label{table:nonNCC}
\setlength{\extrarowheight}{4pt}
{\begin{tabular}{|>{$}c<{$}|>{$}c<{$}|>{$}c<{$}|>{$}c<{$}|>{$}c<{$}|>{$}c<{$}|}
\hline
\frak{g} &\frak{h}& \frak{g}^{\theta\theta_\e}&\Sigma& m_\alpha 
\rule[-0.2cm]{0cm}{0.0cm}
& \\
\hline
\hline
\frak{so}(2n+1,\mathbb{C})&\frak{so}(2(n-j)+1,2j)
&\frak{so}(2(n-j)+1,\mathbb{C})\times
\frak{so}(2j,\mathbb{C})&B_n&2
&n \geq 2, \, 2\leq j\leq n
\rule[-0.2cm]{0cm}{0.0cm}
\cr\hline
\frak{sp}(n,\mathbb{C})&\frak{sp}(n-j,j)
&\frak{sp}(n-j,\mathbb{C})\times \frak{sp}(j,\mathbb{C})&C_n&2
&n \geq 3,\, 1\leq j \leq [n/2]
\rule[-0.2cm]{0cm}{0.0cm}
\cr\hline
\frak{so}(2n,\mathbb{C})&\frak{so}(2(n-j),2j)
&\frak{so}(2(n-j),\mathbb{C})\times\frak{so}(2j,\mathbb{C})&D_n&2
&n\geq 4,\, 2 \leq j \leq [n/2]
\rule[-0.2cm]{0cm}{0.0cm}
\cr\hline
(\frak{e}_6)_{\mathbb{C}}&\frak{e}_{6(2)}
&\frak{sl}(6,\mathbb{C})\times\frak{sl}(2,\mathbb{C})&E_6&2&
\cr\hline
(\frak{e}_7)_{\mathbb{C}}&\frak{e}_{7(7)}
&\frak{sl}(8,\mathbb{C})&E_7&2
\rule[-0.2cm]{0cm}{0.0cm}
&
\cr\hline
(\frak{e}_7)_{\mathbb{C}}&\frak{e}_{7(-5)}
&\frak{so}(12,\mathbb{C})\times\frak{sl}(2,\mathbb{C})&E_7&2
\rule[-0.2cm]{0cm}{0.0cm}
&
\cr\hline
(\frak{e}_8)_{\mathbb{C}}&\frak{e}_{8(8)}
&\frak{so}(16,\mathbb{C})&E_8&2
\rule[-0.2cm]{0cm}{0.0cm}
&
\cr\hline
(\frak{e}_8)_{\mathbb{C}}&\frak{e}_{8(-24)}
&e_{7\mathbb{C}}\times\frak{sl}(2,\mathbb{C})&E_8&2
\rule[-0.2cm]{0cm}{0.0cm}
&
\cr\hline
(\frak{f}_4)_{\mathbb{C}}&\frak{f}_{4(4)}
&\frak{sp}(3,\mathbb{C})\times\frak{sl}(2,\mathbb{C})&F_4&2
\rule[-0.2cm]{0cm}{0.0cm}
&
\cr\hline
(\frak{f}_4)_{\mathbb{C}}&\frak{f}_{4(-20)}
&\frak{so}(9,\mathbb{C}) &F_4&2&
\cr\hline
(\frak{g}_2)_{\smC}&\frak{g}_{2(2)}
&\frak{sl}(2,\mathbb{C})\times\frak{sl}(2,\mathbb{C})&G_2&2
\rule[-0.2cm]{0cm}{0.0cm}
&
\cr\hline
\end{tabular}
\bigskip
\caption{Other $K_\e$-symmetric pairs with even multiplicities.} }
\end{table}

\addtolength{\textheight}{-40pt}

A special isomorphism of $K_\e II$ symmetric pairs with even multiplicities
is
\begin{alignat*}{2}
 \frak{sp}(2,\C) &\approx \frak{so}(5,\C),
  &\quad
\frak{sp}(1,1) &\approx \frak{so}(1,4)\,.
\end{alignat*}


\end{document}